\documentclass[a4paper,12pt]{article}
\usepackage{amsmath}
\usepackage{amsfonts}
\usepackage{amssymb}
\usepackage{latexsym}
\usepackage{epsfig}
\usepackage{graphicx}
\usepackage{oldgerm}
\usepackage{theorem}

\def\C{\mathbb{C}}
\def\R{\mathbb{R}}
\def\N{\mathbb{N}}
\def\W{\mathbb{W}}
\def\Z{\mathbb{Z}}
\def\I{\mathbb{I}}
\def\J{\mathbb{J}}

\def\P{{\mathbb P}}
\def\E{{\mathbb E}}
\def\hP{\hat{\mathbb P}}
\def\hE{\hat{\mathbb E}}

\def\mbK{\mathbb{K}}

\def\bK{{\bf K}}
\def\bP{{\bf P}}
\def\bE{{\bf E}}

\def\rP{{\rm P}}
\def\rE{{\rm E}}
\def\rC{{\rm C}}

\def\1{{\bf 1}}
\def\0{{\bf 0}}
\def\x{\mib{x}}
\def\y{\mib{y}}
\def\z{\mib{z}}
\def\u{\mib{u}}
\def\v{\mib{v}}

\def\bB{\mib{B}}
\def\bR{\mib{R}}
\def\f{\mib{f}}
\def\t{\mib{t}}
\def\vchi{\mib{\chi}}

\def\X{\mib{X}}
\def\V{\mib{V}}
\def\bW{\mib{W}}
\def\bZ{\mib{Z}}

\def\sfM{{\sf M}}

\def\hV{\hat{V}}

\def\hX{\hat{\mib{X}}}

\def\hXi{\hat{\Xi}}

\def\mM{\mathfrak{M}}
\def\hmM{\hat{\mathfrak{M}}}

\def\mX{\mathfrak{X}}

\def\supp{{\rm supp}\,}
\def\cG{{\cal G}}

\def\cM{{\cal M}}

\def\cF{{\cal F}}

\def\sc{{\sf c}}
\def\cC{{\cal C}}
\def\cS{{\cal S}}

\def\Det{\mathop{\mathrm{Det}}}
\def\bM{{\bf M}}
\def\bJ{{\bf J}}

\def\uj{\underline{j}}
\def\oj{\overline{j}}

\def\cT{{\cal T}}
\def\cD{{\cal D}}

\def\cH{{\cal H}}
\def\cF{{\cal F}}
\def\cL{{\cal L}}
\def\cP{{\cal P}}

\def\B{\mib{B}}
\def\T{\mib{T}}
\def\F{\mib{F}}
\def\k{\mib{k}}
\def\e{\mib{e}}
\def\vnu{\mib{\nu}}
\def\vmu{\mib{\mu}}
\def\vrho{\mib{\rho}}

\def\mbT{\mathbb{T}}

\def\hnu{\widehat{\nu}}

\def\d={\stackrel{\rm d}{=}}

\theorembodyfont{\itshape}

\newtheorem{thm}{Theorem}[section]
\newtheorem{lem}[thm]{Lemma}
\newtheorem{cor}[thm]{Corollary}
\newtheorem{prop}[thm]{Proposition}

\newtheorem{df}[thm]{Definition}

\newcommand{\mib}[1]{\mbox{\boldmath $#1$}}
\newcommand{\SSC}[1]{\section{#1}\setcounter{equation}{0}}
\newcommand{\qed}{\hbox{\rule[-2pt]{3pt}{6pt}}}

\begin{document}
\title{\bf 
Determinantal Martingales \\
and Interacting Particle Systems
\footnote{
This is a note
prepared for the lectures at 
Faculty of Mathematics, Kyushu University,
24-28 June 2013.
}}
\author{
Makoto Katori
\footnote{
Department of Physics,
Faculty of Science and Engineering,
Chuo University, 
Kasuga, Bunkyo-ku, Tokyo 112-8551, Japan;
e-mail: katori@phys.chuo-u.ac.jp
}}
\date{9 July 2013}
\pagestyle{plain}
\maketitle
\begin{abstract}
Determinantal process is a dynamical extension of
a determinantal point process such that
any spatio-temporal correlation function is given by
a determinant specified by a single continuous function
called the correlation kernel.
Noncolliding diffusion processes are important examples
of determinantal processes.
In the present lecture, we introduce determinantal martingales
and show that if the interacting particle system (IPS)
has determinantal-martingale representation,
then it becomes a determinantal process.
From this point of view, the reason why noncolliding diffusion
processes and noncolliding random walk 
are determinantal is simply explained 
by the fact that the harmonic transform with the Vandermonde
determinant provides a proper determinantal martingale.
Recently O'Connell introduced an interesting IPS, which can be regarded
as a stochastic version of a quantum Toda lattice.
It is a geometric lifting of the noncolliding Brownian motion
and is not determinantal,
but Borodin and Corwin discovered a determinant formula
for a special observable for it.
We also discuss this new topic from the present
view-point of determinantal martingale.
\end{abstract}

\footnotesize 
\tableofcontents
\vspace{3mm}
\normalsize

\clearpage
\SSC{Introduction \label{sec:Introduction}}
\subsection{Determinantal martingale \label{sec:det_mar}}

Let $V(t), t \in \cT$ be a Markov process in a state space
$S \subset \R$, where the set of time is continuous
$\cT=[0, \infty)$ or discrete $\cT=\N_0 \equiv \{0,1,2, \dots\}$.
The probability space is denoted as $(\Omega, \cF, \rP)$,
where the expectation is written as $\rE$.
We introduce a filtration $\{\cF(t) : t \in \cT \}$
defined by $\cF(t)=\sigma(V(s), s \in [0, t] \cap \cT)$.
When $V(t), t \in [0, \infty)$ is a continuous time Markov process,
provided it has right-continuous sample paths, the
transition probability density is given by
\begin{equation}
p(t, y|x)= \frac{\partial}{\partial y}
\rP(V(t) \leq y| V(0)=x), \quad x, y \in S,
\label{eqn:tpd1}
\end{equation}
and when $V(t), t \in \N_0$ is a Markov chain,
the transition probability is given by
\begin{equation}
p(t,y|x)=\rP(V(t)=y|V(0)=x), \quad x, y \in S.
\label{eqn:tpd2}
\end{equation}

The process $M(t), t \in \cT$ on $(\Omega, \cF, \rP)$
adapted to $\{\cF(t) : t \in \cT \}$,
is said to be a {\it martingale} if, for every $s < t, s, t \in \cT$,
\begin{equation}
\rE[M(t)|\cF(s)]=M(s) \quad \mbox{a.s.} \quad \rP.
\label{eqn:martingale1}
\end{equation}
Let $\N = \{1, 2, \dots\}$.
For a Markov process $V(t), t \in [0, \infty)$, 
if there exists a nondecreasing sequence 
$\{T_n : n \in \N \}$
of stopping times of $\{\cF(t) : t \in [0, \infty) \}$,
such that
$M_n(t) \equiv M(t \wedge T_n), t \in \cT$ is a martingale
for every $n \in \N$ and
$\displaystyle{\rP \left[\lim_{n \to \infty} T_n = \infty \right]=1}$,
then we say $M(t), t \in \cT$ is a {\it local martingale}.

For $N \in \N$, we put
\begin{equation}
\W_N=\{\x=(x_1, \dots, x_N) \in S^{N} : x_1 < x_2 < \cdots < x_N \}.
\label{eqn:Weyl1}
\end{equation}
For $\u=(u_1, \dots, u_N) \in \W_N$,
we define a measure $\xi$ by a sum of
point masses concentrated on $u_j$'s, $1 \leq j \leq N$,
\begin{equation}
\xi(\cdot)=\sum_{j=1}^N \delta_{u_j}(\cdot).
\label{eqn:xiAA1}
\end{equation}
Depending on $\xi$, we 
assume that there is a one-parameter
family of maps
\begin{equation}
\cM_{\xi}^u(\cdot, \cdot) : 
\cT \times \R \mapsto \R
\label{eqn:map1}
\end{equation}
with a parameter $u \in \C$, 
such that
\begin{eqnarray}
&(i)& \quad
\cM_{\xi}^u(\cdot, V(\cdot)) 
\, \, \mbox{is a local martingale},
\nonumber\\
&(ii)& \quad
\cM_{\xi}^{u_k}(\cdot, x) : 1 \leq k \leq N 
\quad \mbox{are linearly independent functions of $x$},
\nonumber\\
&(iii)& \quad
\cM_{\xi}^{u_k}(0,u_j)=\delta_{jk}, \quad
u_j, u_k \in \supp \xi =\{u_1, \dots, u_N\}.
\label{eqn:map2}
\end{eqnarray}

Let $\{V_j(t), t \in \cT : 1 \leq j \leq N \}$
be a collection of $N$ independent copies 
of $V(t), t \in \cT$.
We consider the $N$-component vector-valued process
$\V(t)=(V_1(t), \dots, V_N(t))$, $t \in \cT$,
for which the initial values are fixed to be
$V_j(0)=u_j \in S, 1 \leq j \leq N$
and
the probability space is denoted by
$(\Omega, \cF, \rP_{\u})$ with expectation
$\rE_{\u}, \u=(u_1, \dots, u_N)$. 

For $n \in \N$, let
$\I_{n} = \{1,2, \dots, n\}$.
Let  $\x=(x_1, \dots, x_N) \in S^N$
and $1 \leq N' \leq N$.
The cardinality of a finite set $A$
is denoted by $\sharp A$.
We write $\J \subset \I_N, \sharp \J=N'$,
if $\J=\{j_1, \dots, j_{N'}\},
1 \leq j_1 < \dots < j_{N'} \leq N$,
and put $\x_{\J}=(x_{j_1}, \dots, x_{j_{N'}})$.
In particular, we write
$\x_{N'}=\x_{\I_{N'}}, 1 \leq N' \leq N$.
(By definition $\x_N=\x$.)
Suppose $\u \in \W_N$ and 
$\xi(\cdot)=\sum_{j=1}^N \delta_{u_j}(\cdot)$.
For $\J \subset \I_N, 1 \leq \sharp \J \leq N$,
consider a determinant of local martingales
\begin{equation}
\cD_{\xi}(t, \V_{\J}(t)) 
=\det_{j,k \in \J} \Big[ 
\cM_{\xi}^{u_k}(t, V_j(t)) \Big],
\quad t \in \cT.
\label{eqn:D1}
\end{equation}
We call (\ref{eqn:D1}) 
a {\it determinantal martingale} \cite{Kat13}.

Let $t \in \cT, t \leq T \in \cT$.
In this lecture we study the following expectation
for an $\cF(t)$-measurable function $F$ of $\V(\cdot)$,
which is symmetric at each time, 
weighted by the determinantal martingale,
\begin{equation}
\rE_{\u} \Big[
F(\V(\cdot)) \cD_{\xi}(T, \V(T)) \Big],
\quad \u \in S^N.
\label{eqn:exp1}
\end{equation}
By the assumptions (\ref{eqn:map2}) for the map (\ref{eqn:map1}), 
we can prove the following.
\begin{lem}
\label{thm:reducibility}
Assume that
$\xi(\cdot)=\sum_{j=1}^N \delta_{u_j}(\cdot)$
with $\u \in \W_N$.
Let $1 \leq N' \leq N$.
For $t \in \cT, t \leq T < \infty$ and a measurable function
$F_{N'}$ on $S^{N'}$,
\begin{eqnarray}
&& \sum_{\J \subset \I_N, 
\sharp \J=N'}
\rE_{\u} \left[
F_{N'}(\V_{\J}(t))
\cD_{\xi}(T, \V(T)) \right]
\nonumber\\
&& \quad
= \int_{\W_{N'}} \xi^{\otimes N'} (d\v)
\rE_{\v} \left[
F_{N'}(\V_{N'}(t))
\cD_{\xi}(T, \V_{N'}(T)) \right].
\label{eqn:reducibility}
\end{eqnarray}
\end{lem}
This shows the
{\it reducibility} of the determinantal martingale
in the sense that,
if we observe a symmetric function depending
on $N'$ variables, $N' \leq N$,
then the size of determinantal
martingale can be reduced from $N$ to $N'$.
Proof is given in Section \ref{sec:proof_red}.

For an integer $M \in \N$, consider 
a sequence of times
$0 < t_1 < \cdots < t_M \leq T \in \cT$,
$t_m \in \cT, 1 \leq m \leq M$, 
and a sequence of measurable functions
$\vchi=(\chi_{t_1},\dots,\chi_{t_M})$.
Then given an integral kernel
\begin{equation}
\bK(s,x;t,y); 
\quad 
(s,x), (t,y) \in \cT \times S,
\label{eqn:kernel1}
\end{equation}
the {\it Fredholm determinant} is defined as
\begin{eqnarray}
&& \mathop{{\rm Det}}_
{\substack{
(s,t)\in \{t_1, \dots, t_M\}, \\
(x,y)\in S^2}
}
 \Big[\delta_{st} \delta_x(y)
+ \bK(s,x;t,y) \chi_{t}(y) \Big]
\nonumber\\
&& 
=\sum_
{\substack
{N_m \geq 0, \\ 1 \leq m \leq M} }
\int_{\prod_{m=1}^{M} \W_{N_{m}}}
\prod_{m=1}^{M} \left\{ d \x_{N_m}^{(m)}
\prod_{j=1}^{N_{m}} 
\chi_{t_m} \Big(x_{j}^{(m)} \Big) \right\}
\det_{\substack
{1 \leq j \leq N_{m}, 1 \leq k \leq N_{n}, \\
1 \leq m, n \leq M}
}
\Bigg[
\bK(t_m, x_{j}^{(m)}; t_n, x_{k}^{(n)} )
\Bigg],
\nonumber\\
\label{eqn:Fredholm1}
\end{eqnarray}
where $\x^{(m)}_{N_m}$ denotes
$(x^{(m)}_1, \dots, x^{(m)}_{N_m})$
and
$d \x^{(m)}_{N_m}= \prod_{j=1}^{N_m} dx^{(m)}_j$,
$1 \leq m \leq M$.
Let $\1(\omega)$ be the indicator of $\omega$;
$\1(\omega)=1$ if $\omega$ is satisfied,
and $\1(\omega)=0$ otherwise.
The reducibility of determinantal martingales 
(Lemma \ref{thm:reducibility}) implies the
following identity.
\begin{lem}
\label{thm:Fredholm}
Let $\u \in \W_N$ and
$\xi=\sum_{j=1}^N \delta_{u_j}$.
For $t \in \cT, 0 \leq t \leq T \in \cT$,
\begin{eqnarray}
&& \rE_{\u} \left[
\prod_{m=1}^M \prod_{j=1}^N
\{1+\chi_{t_m}(V_j(t_m)) \}
\cD_{\xi}(T, \V(T)) \right]
\nonumber\\
&& \qquad =
\mathop{{\rm Det}}_
{\substack{
(s,t)\in \{t_1, \dots, t_M\}, \\
(x,y)\in S^2}
}
 \Big[\delta_{st} \delta_x(y)
+ \mbK_{\xi}(s,x;t,y) \chi_{t}(y) \Big],
\label{eqn:Fredholm2}
\end{eqnarray}
where 
\begin{equation}
\mbK_{\xi}(s,x;t,y)
=\int_{S} \xi(dv) p(s, x|v) \cM_{\xi}^{v}(t,y)
- \1(s>t) p(s-t,x|y).
\label{eqn:K1}
\end{equation}
\end{lem}
Proof is given in Section \ref{sec:proof_Fred}.

\subsection{Determinantal-martingale representations (DMR) \label{sec:DM_rep}}

Let $\mM$ be the space of nonnegative integer-valued Radon measures 
on $S$.
Any element $\xi$ of $\mM$ can be represented as
$\xi(\cdot) = \sum_{j \in \I}\delta_{x_j}(\cdot)$
with a countable index set $\I$, in which
a sequence of points in $S$, $\x =(x_j)_{j \in \I}$ 
satisfying $\xi(K)=\sharp\{x_j, x_j \in K \} < \infty$ 
for any compact subset $K \subset S$.
In this lecture, we consider 
interacting particle systems 
as $\mM$-valued processes and write them as
\begin{equation}
\Xi(t, \cdot)=\sum_{j=1}^N \delta_{X_j(t)}(\cdot),
\quad t \in \cT.
\label{eqn:Xi1}
\end{equation}
The probability law of $\Xi(t, \cdot)$
starting from a fixed configuration $\xi \in \mM$
is denoted by $\P_{\xi}$ 
and the process
specified by the initial configuration
is expressed by
$(\Xi(t), \P_{\xi})$.
The expectations w.r.t.$\P_{\xi}$
is denoted by $\E_{\xi}$.
We introduce a filtration $\{{\cal F}(t) \}_{t \in \cT}$
defined by ${\cal F}(t) = \sigma (\Xi(s), s \in \cT \cap [0,t])$. 
Let $\rC_0(S)$ be the set of all continuous
real-valued functions with compact supports on $S$.
We set 
\begin{equation}
\mM_{0}= \{ \xi \in \mM : 
\xi(\{x\}) \leq 1 \mbox { for any }  x \in S \},
\label{eqn:mM0}
\end{equation}
which denotes a collection of configurations
without any multiple points.

For any integer $M \in \N$,
a sequence of times
$\t=(t_1,\dots,t_M) \in \cT^M$ with 
$0 < t_1 < \cdots < t_M \leq T \in \cT$,
$t_m \in \cT, 1 \leq m \leq M$, 
and a sequence of functions
$\f=(f_{t_1},\dots,f_{t_M}) \in \rC_0(S)^M$,
the {\it moment generating function} of multitime distribution
of $(\Xi(t), \P_{\xi})$ is defined by
\begin{equation}
\Psi_{\xi}^{\t}[\f]
\equiv \E_{\xi} \left[ \exp \left\{ \sum_{m=1}^{M} 
\int_{S} f_{t_m}(x) \Xi(t_m, dx) \right\} \right].
\label{eqn:GF1}
\end{equation}
It is expanded w.r.t.
$\chi_{t_m}(\cdot)=e^{f_{t_m}(\cdot)}-1,
1 \leq m \leq M$ as
\begin{equation}
\Psi_{\xi}^{\t}[\f]
=\sum_
{\substack
{N_m \geq 0, \\ 1 \leq m \leq M} }
\int_{\prod_{m=1}^{M} \W_{N_{m}}}
\prod_{m=1}^{M} \left\{ d \x_{N_m}^{(m)}
\prod_{j=1}^{N_{m}} 
\chi_{t_m} \Big(x_{j}^{(m)} \Big) \right\}
\rho_{\xi} 
\Big( t_{1}, \x^{(1)}_{N_1}; \dots ; t_{M}, \x^{(M)}_{N_M} \Big),
\label{eqn:GFB1}
\end{equation}
and it defines the {\it spatio-temporal correlation functions}
$\rho_{\xi}(\cdot)$ for the process $(\Xi(t), \P_{\xi})$.

We introduce the following definitions.

\begin{df}
\label{thm:determinantal}
Let $t \in \cT$, $0 < t \leq T \in \cT$.
For an $\cF(t)$-measurable bounded function $F$,
if $\E_{\xi}[F(\Xi(\cdot))]$ is expressed by
a Fredholm determinant,
$\E_{\xi}[F(\Xi(\cdot))]$ is said to be 
Fredholm determinantal (F-determinantal, for short).
If the moment generating function (\ref{eqn:GF1}) 
is Fredholm-determinantal, we say
the process $(\Xi(t), \P_{\xi})$ is determinantal. 
In this case, 
all spatio-temporal correlation functions
are given by determinants as 
\begin{equation}
\rho_{\xi} \Big(t_1,\x^{(1)}_{N_1}; \dots;t_M,\x^{(M)}_{N_M} \Big) 
=\det_{\substack
{1 \leq j \leq N_{m}, 1 \leq k \leq N_{n}, \\
1 \leq m, n \leq M}
}
\Bigg[
\mbK_{\xi}(t_m, x_{j}^{(m)}; t_n, x_{k}^{(n)} )
\Bigg],
\label{eqn:rho1}
\end{equation}
$0 < t_1 < \cdots < t_M \leq T \in \cT$,
$t_m \in \cT, 1 \leq m \leq M$, 
$1 \leq N_m \leq N$,
$\x^{(m)}_{N_m} \in S^{N_m}, 1 \leq m \leq M \in \N$,
and the integral kernel $\mbK_{\xi}$ 
depending on $\xi$ is called
the correlation kernel.
\end{df}
\begin{df}
\label{thm:DM_rep}
Let $t \in \cT, 0 < t \leq T \in \cT$.
If there exists $\cM_{\xi}^{u}(\cdot, \cdot), u \in S$
defining $\cD_{\xi}(\cdot, \cdot)$ by (\ref{eqn:D1}) such that
the following equality holds 
for an $\cF(t)$-measurable bounded function $F$,
\begin{equation}
\E_{\xi}[F(\Xi(\cdot))]
=\rE_{\u} \left[
F \left( \sum_{j=1}^N \delta_{V_j(\cdot)} \right) 
\cD_{\xi}(T, \V(T)) \right],
\label{eqn:DM1}
\end{equation}
then we say $(\Xi(t), \P_{\xi})$ has a determinantal-martingale
representation (DMR, for short) for $F$.
If $(\Xi(t), \P_{\xi})$ has DMR for any $\cF(t)$-measurable bounded function,
$t \in \cT, 0 < t \leq T \in \cT$, 
it is said to be having DMR.
\end{df}

Lemma \ref{thm:Fredholm} gives the following statement.
\begin{prop}
\label{thm:DM_det}
If $(\Xi(t), \P_{\xi})$ has DMR 
for $F$, then
$\E_{\xi}[F(\Xi(\cdot))]$ is 
F-determinantal.
If $(\Xi(t), \P_{\xi})$ has DMR, then
it is determinantal.
\end{prop}

In the present lectures, we will prove the following.
\begin{description}
\item{1.} \quad
The noncolliding Brownian motion (BM) and
the noncolliding squared Bessel process 
with parameter $\nu >-1$ (BESQ$^{(\nu)}$)
have DMR
for $\xi \in \mM_0, \xi(S) < \infty$
(Theorem \ref{thm:DMR1}).
Then they are both determinantal
(Corollaries \ref{thm:K_BM} and \ref{thm:K_BESQ}).

\item{2.} \quad
The simple and symmetric noncolliding random walk (RW)
has DMR
for $\xi \in \mM, \xi(\Z) \in \N$
(Theorem \ref{thm:DMR_RW}).
Then it is determinantal
(Corollary \ref{thm:main_RW}).

\item{3.} \quad
The O'Connell process has a {\it variation} of DMR for a special quantity,
which will be denoted as
\begin{equation}
\Theta^a(X^a_1(\cdot)-h), \quad a>0, \quad h \in \R.
\label{eqn:Theta1}
\end{equation}
Then its expectation is F-determinantal
(Proposition \ref{thm:DMP_OP1}).
\end{description}

\subsection{Complex-process representations (CPR) \label{sec:CP_rep}}

Let $W(t), t \in \cT$ be a Markov process in $\widetilde{S}$ started at 0
defined independently from $V(t), t \in \cT$ on the 
probability space $(\check{\Omega}, \check{\cF}, \check{\rP}_0)$,
where the expectation is denoted as $\check{\rE}_0$.
We introduce a complex process
\begin{equation}
Z(t)=V(t)+ i W(t), \quad t \in \cT,
\label{eqn:Z1}
\end{equation}
where $i=\sqrt{-1}$.
We consider a possibility that there exists 
a one-parameter family of functions $\varphi^{u}_{\xi} : \C \to \C$
with the parameter $u \in \C$ and $\xi=\sum_{j=1}^N \delta_{u_j},
\u \in \W_N$ such that the equality
\begin{equation}
\cM_{\xi}^{u}(t, V(t))=
\check{\rE}_0 [ \varphi^u_{\xi}(Z(t))],
\quad t \in \cT, t \leq T < \infty
\label{eqn:Z2}
\end{equation}
hold.
In this case, we set a collection of $N$ independent
copies of $W(\cdot)$ and denote the probability space
as $(\check{\Omega}, \check{\cF}, \check{\rP}_{\0})$,
where $\0$ denotes the zero in $S^N$.
Define the space $(\Omega, \cF, \bP_{\u})$
as a product of $(\Omega, \cF, \rP_{\u})$ and
$(\check{\Omega}, \check{\cF}, \check{\rP}_{\0})$,
which is the probability space for the $N$-component
complex vector-valued process
$\bZ(t)=(Z_1(t), \cdots, Z_N(t))$
with $Z_j(t)=V_j(t)+i W_j(t), 1 \leq j \leq N$,
$t \in \cT$.
Then the determinantal martingale (\ref{eqn:D1}) is written as
\begin{eqnarray}
\cD_{\xi}(t, \V_{\J}(t))
&=& \det_{j, k \in \J} \left[\check{\rE}_0[\varphi^{u_k}_{\xi}(Z_j(t))] \right]
\nonumber\\
&=& \check{\rE}_{\0} \left[
\det_{j, k \in \J} [ \varphi^{u_k}_{\xi}(Z_j(t))] \right],
\label{eqn:D2}
\end{eqnarray}
for $\J \in \I_N$ and the DMR
(\ref{eqn:DM1}) is rewritten as
\begin{equation}
\E_{\xi}[F(\Xi(\cdot))]
=\bE_{\u} \left[
F \left( \sum_{j=1}^N \delta_{\Re Z_j(\cdot)} \right) 
\det_{1 \leq j, k \leq N} [ \varphi^{u_k}_{\xi}(Z_j(T))
\right],
\quad t \in \cT, t \leq T < \infty.
\label{eqn:CPR1}
\end{equation}
We call (\ref{eqn:CPR1}) the 
{\it complex-process representation}
(CPR, for short)
for $\E_{\xi}[F(\Xi(\cdot))]$.

We will prove the following.
\begin{description}
\item{4.} \quad
The noncolliding BM and
the noncolliding BES$^{(\nu)}$ with
$\nu=2n+1, n \in \N_0$
have CPRs 
(Corollaries \ref{thm:CPR_noncBM}, \ref{thm:CPR_noncBES}).

\item{5.} \quad
The noncolliding RW
has CPR
(Theorem \ref{thm:DMR_RW}).

\item{6.} \quad
The O'Connell process has {\it a variation} of CPR for 
$\Theta^a(X^a_1(t)-h), a >0, h \in \R$.
(Proposition \ref{thm:CPR_OP1})
\end{description}

\subsection{Infinite particle systems \label{sec:IPS}}

Under some conditions on initial configuration $\xi \in \mM_0$
for interacting particle system $(\Xi(t), \P_{\xi})$,
the map $\cM_{\xi}^{u}(\cdot, \cdot), u \in \C$ will be
well-defined even for infinite particle limits,
$\xi(S)=N \to \infty$.
Assume $M \in \N$, $0 < t_1 < \cdots < t_M \leq T < \infty$,
$\phi_m \in \rC_0(S), 1 \leq m \leq M$
and $G=G(\{x_m\}_{m=1}^M)$
is a polynomial on $S^M$.
For $0 < t \leq T < \infty$, if an $\cF(t)$-measurable
function $F(\Xi(\cdot))$ is represented as
$$
F(\Xi(\cdot))
=G \left( \left\{ \int_{S} \phi_m(x) \Xi(t_m, dx) \right\}_{m=1}^M \right),
$$
we say $F$ is polynomial.
By the reducibility of determinantal martingale
(Lemma \ref{thm:reducibility}),
if the degree of the polynomial $F$ is $n \in \N$,
$\cD_{\xi}(T, \V(T))$ in the DMR
(\ref{eqn:DM1}) can be expressed by using
$\{\cD_{\xi}(T,\V_{\J}(T)) : \J \subset \N$, $\sharp \J \leq n \}$.
In this sense, the present DMRs (and CPRs)
are valid for $(\Xi(t), \P_{\xi})$ also in the case
$\xi(S)=\infty$.
We will discuss the following.

\begin{description}
\item{7.} \quad
Some sufficient conditions for $\xi$ are
given such that the DMRs
(and CPRs)
are valid for the interacting particle systems
$(\Xi(t), \P_{\xi})$ 
with infinite numbers of particles, 
and infinite-dimensional determinantal processes are
well-defined (Section \ref{sec:IPS2}).
\end{description}

\subsection{Systems started at initial configurations 
with multiple points
\label{sec:initial}}

The basic of the present theory is for the deterministic
initial configuration with no multiple point,
$\xi \in \mM_0$.
But we will also discuss the following case.

\begin{description}
\item{8.} \quad
The map $\cM_{\xi}^{u}(\cdot, \cdot), u \in \C$
can be extended for the system
started at configurations with multiple points
(Section \ref{sec:multi}).

\end{description}

\clearpage
\SSC{Proofs of Lemmas \ref{thm:reducibility} and 
\ref{thm:Fredholm} \label{sec:proof1}}

For a finite set $\J$, we write the collection of all permutations
of elements in $\J$ as $\cS(\J)$.
In particular, 
for $\I_n=\{1,2,\dots, n \}, n \in N$, 
we express $\cS(\I_{n})$ simply by $\cS_{n}$.
For $\x=(x_1, \dots, x_N) \in S^n$, $\sigma \in \cS_n$, 
we put $\sigma(\x)=(x_{\sigma(1)}, \dots, x_{\sigma(n)})$.
For an $n \times n$ matrix $B=(B_{jk})_{1 \leq j,k \leq n}$, 
the determinant is defined by
\begin{eqnarray}
\det B &=& \det_{1 \leq j, k \leq n} [B_{jk}]
\nonumber\\
&=& \sum_{\sigma \in \cS_{n}} {\rm sgn}(\sigma)
\prod_{j=1}^n B_{j \sigma(j)}.
\label{eqn:det1}
\end{eqnarray}
Any permutation
$\sigma \in \cS_n$ can be decomposed
into a product of cycles.
Let the number of cycles in the decomposition
be $\ell(\sigma)$ and express $\sigma$ by
\begin{equation}
\sigma = \sc_1 \sc_2 \cdots \sc_{\ell(\sigma)},
\label{eqn:det2a}
\end{equation}
where $\sc_{\lambda}$
denotes a cyclic permutation
\begin{equation}
\sc_{\lambda}= 
(c_\lambda(1) c_\lambda(2) \cdots c_\lambda(q_{\lambda}) ), \quad
1 \leq q_{\lambda} \leq n,  
\quad1 \leq \lambda \leq \ell(\sigma).
\label{eqn:det2b}
\end{equation}
For each $1 \leq \lambda \leq \ell(\sigma)$,
we write the set of entries 
$\{c_{\lambda}(j)\}_{j=1}^{q_{\lambda}}$ of $\sc_{\lambda}$
simply as $\{\sc_{\lambda}\}$,
in which the periodicity
$c_{\lambda}(j+q_{\lambda})=c_{\lambda}(j), 1 \leq j \leq q_{\lambda}$ 
is assumed.
By definition, for each $1 \leq \lambda \leq \ell(\sigma)$,
$c_{\lambda}(j), 1 \leq j \leq q_{\lambda}$
are distinct indices chosen from $\I_n$, 
$\{\sc_{\lambda}\} \cap \{\sc_{\lambda'}\} = \emptyset$
for $1 \leq \lambda \not= \lambda' \leq \ell(\sigma)$, and
$\sum_{\lambda=1}^{\ell(\sigma)} q_{\lambda}=n$.
The determinant (\ref{eqn:det1}) is also given by
\begin{equation}
\det B = \sum_{\sigma \in \cS_n}
(-1)^{n-\ell(\sigma)}
\prod_{\lambda=1}^{\ell(\sigma)}
\prod_{j=1}^{q_{\lambda}} 
B_{c_{\lambda}(j) c_{\lambda}(j+1)}.
\label{eqn:det2}
\end{equation}

\subsection{Proof of Lemma \ref{thm:reducibility} \label{sec:proof_red}}

By definition of determinant (\ref{eqn:det1}),
(\ref{eqn:D1}), and independence of $V_j(\cdot), 1 \leq j \leq N$,
the LHS of (\ref{eqn:reducibility}) is equal to
\begin{eqnarray}
&& \sum_{\J \subset \I_N, \sharp \J=N'}
\rE_{\u} \left[ F_{N'}(\V_{\J}(t))
\det_{j, k \in \I_N}
[\cM_{\xi}^{u_k}(T, V_j(T))] \right]
\nonumber\\
&& \, = \sum_{\J \subset \I_N, \sharp \J=N'}
\rE_{\u} \left[ F_{N'}(\V_{\J}(t))
\sum_{\sigma \in \cS_N} {\rm sgn}(\sigma)
\prod_{j=1}^N \cM_{\xi}^{u_{\sigma(j)}}(T, V_{j}(T)) \right]
\nonumber\\
&& \, = \sum_{\J \subset \I_N, \sharp \J=N'}
\sum_{\sigma \in \cS_N} {\rm sgn}(\sigma)
\rE_{\u} \left[ F_{N'}(\V_{\J}(t))
\prod_{j \in \J} \cM_{\xi}^{u_{\sigma(j)}}(T, V_j(T))
\prod_{k \in \I_N \setminus \J}
\cM_{\xi}^{u_{\sigma(k)}}(T, V_k(T)) \right]
\nonumber\\
&& \, = \sum_{\J \subset \I_N, \sharp \J=N'}
\sum_{\sigma \in \cS_N} {\rm sgn}(\sigma)
\rE_{\u} \left[ F_{N'}(\V_{\J}(t))
\prod_{j \in \J} \cM_{\xi}^{u_{\sigma(j)}}(T, V_j(T)) \right]
\nonumber\\
&& \qquad \qquad \qquad \qquad \qquad \times
\prod_{k \in \I_N \setminus \J}
\rE_{\u} \left[ \cM_{\xi}^{u_{\sigma(k)}}(T, V_k(T)) \right].
\label{eqn:A1}
\end{eqnarray}
By the martingale property of $\cM_{\xi}^u(\cdot, V(\cdot))$
and the condition (iii) of (\ref{eqn:map2}), 
\begin{eqnarray}
\prod_{k \in \I_N \setminus \J}
\rE_{\u} \left[ \cM_{\xi}^{u_{\sigma(k)}}(T, V_k(T)) \right]
&=&
\prod_{k \in \I_N \setminus \J}
\rE_{\u} \left[ \cM_{\xi}^{u_{\sigma(k)}}(0, V_k(0)) \right]
\nonumber\\
&=& \prod_{k \in \I_N \setminus \J} \cM_{\xi}^{u_{\sigma(k)}}(0, u_k) 
\nonumber\\
&=& \prod_{k \in \I_N \setminus \J} \delta_{k \sigma(k)}.
\label{eqn:A2}
\end{eqnarray}
Then (\ref{eqn:A1}) is equal to
\begin{eqnarray}
&&  \sum_{\J \subset \I_N, \sharp \J=N'}
\sum_{\sigma \in \cS(\J)} {\rm sgn}(\sigma)
\rE_{\u} \left[ F_{N'}(\V_{\J}(t))
\prod_{j \in \J} \cM_{\xi}^{u_{\sigma(j)}}(T, V_j(T)) \right]
\nonumber\\
&& \quad =
\sum_{\J \subset \I_N, \sharp \J=N'}
\rE_{\u} \left[ F_{N'}(\V_{\J}(t))
\det_{j, k \in \J} [ \cM_{\xi}^{u_k}(T, V_j(T))] \right]
\nonumber\\
&& \quad =
\int_{\W_{N'}} \xi^{\otimes N'}(d \v)
\rE_{\v} \left[ F_{N'}(\V_{N'}(t))
\det_{j, k \in \I_{N'}} [ \cM_{\xi}^{u_k}(T, V_j(T))] \right],
\label{eqn:A3}
\end{eqnarray}
where equivalence of $V_j(\cdot), 1 \leq j \leq N$ 
in probability law is used.
This is the RHS of (\ref{eqn:reducibility}) and 
the proof is completed. \qed

\subsection{Proof of Lemma \ref{thm:Fredholm} \label{sec:proof_Fred}}

By performing binomial expansion of
$\prod_{m=1}^M \prod_{j=1}^N
\{1+\chi_{t_m}(V_j(t_m)) \}$
and by Lemma \ref{thm:reducibility},
the LHS of (\ref{eqn:Fredholm2}) gives
\begin{equation}
\sum_
{\substack
{N_m \geq 0, \\ 1 \leq m \leq M} }
\sum_{1 \leq p \leq N}
\sum_{\substack{
\sharp \J_m=N_m, \\ 1 \leq m \leq M, \\
\bigcup_{m=1}^{M} \J_m= \I_{p}}}
\int_{\W_{p}} \ 
\xi^{\otimes p}(d\v)
\rE_{\v} \left[
\prod_{m=1}^{M} 
\prod_{j_m\in \J_m} \chi_{t_m}(V_{j_m}(t_m))
\cD_{\xi}(T, \V_p(T))\Big] \right].
\label{eqn:B1}
\end{equation}
On the other hand, the RHS of (\ref{eqn:Fredholm2})
has an expansion according to (\ref{eqn:Fredholm1}).
Then, for proof of Lemma \ref{thm:Fredholm},
it is enough to show 
that the following equality is established
for any $M \in \N, (N_1, \dots, N_M) \in \N^M$
\begin{eqnarray}
&& \int_{\prod_{m=1}^{M} \W_{N_m}} 
\prod_{m=1}^{M} \left\{ d \x_{N_m}^{(m)}
\prod_{j=1}^{N_m} \chi_{t_m} 
\Big(x_{j}^{(m)} \Big) \right\}
\det_{\substack{1 \leq j \leq N_{m}, 1 \leq k \leq N_{n},
\\ 1\le m, n \le M}}
\Bigg[
\mbK_{\xi}(t_m, x_{j}^{(m)}; t_{n}, x_{k}^{(n)} )
\Bigg]
\nonumber\\
&=&
\sum_{1 \leq p \leq N}
\sum_{\substack{
\sharp \J_m=N_m, \\ 1 \leq m \leq M, \\
\bigcup_{m=1}^{M} \J_m= \I_{p}}}
\int_{\W_{p}} \ 
\xi^{\otimes p}(d\v)
\rE_{\v} \left[
\prod_{m=1}^{M} 
\prod_{j_m\in \J_m} \chi_{t_m}(V_{j_m}(t_m))
\det_{j,k \in \I_p}
\Big[\cM_{\xi}^{v_k}(T, V_j(T))\Big] \right].
\nonumber\\
\label{eqn:E=det}
\end{eqnarray}
Here we will prove (\ref{eqn:E=det}) by 
fixing $M \in \N$, $(N_1, \dots, N_M) \in \N^M$.

Let $\I^{(1)}=\I_{N_1}$ and 
$\I^{(m)}=\I_{\sum_{j=1}^{m} N_{j}} 
\setminus \I_{\sum_{j=1}^{m-1} N_{j}}, 2 \leq m \leq M$.
Put $n=\sum_{m=1}^{M} N_m$ and
$\tau_j=\sum_{m=1}^{M} t_m \1(j \in \I^{(m)}),
1 \leq j \leq n$.
Then the integrand in the LHS of (\ref{eqn:E=det})
is simply written as
\begin{equation}
\prod_{j=1}^n \chi_{\tau_j}(x_j)
\det_{1 \leq j, k \leq n} [ 
\mbK_{\xi}(\tau_j, x_j; \tau_k, x_k)],
\label{eqn:det_Z1}
\end{equation}
and the integral 
$\int_{\prod_{m=1}^{M} \W_{N_m}} 
\prod_{m=1}^M d \x_{N_m}^{(m)} (\cdot)$
can be replaced by
$\{\prod_{m=1}^{M} N_m ! \}^{-1}$
$\int_{S^n} d \x \, (\cdot)$.
By the formula (\ref{eqn:det2}) of determinant,
(\ref{eqn:det_Z1}) is expressed as
\begin{eqnarray}
&& \sum_{\sigma \in \cS_n} (-1)^{n-\ell(\sigma)}
\prod_{\lambda=1}^{\ell(\sigma)}
\prod_{j=1}^{q_{\lambda}}
\chi_{\tau_{c_{\lambda}(j)}}(x_{c_{\lambda}(j)})
\mbK_{\xi}(\tau_{c_{\lambda}(j)}, x_{c_{\lambda(j)}};
\tau_{c_{\lambda}(j+1)}, x_{c_{\lambda}(j+1)})
\nonumber\\
&=& 
\sum_{\sigma \in \cS_n} (-1)^{n-\ell(\sigma)}
\prod_{\lambda=1}^{\ell(\sigma)}
\prod_{j=1}^{q_{\lambda}}
\chi_{\tau_{c_{\lambda}(j)}}(x_{c_{\lambda}(j)})
\Big\{
\cG_{\tau_{c_{\lambda}(j)},\tau_{c_{\lambda}(j+1)}}
(x_{c_{\lambda}(j)}, x_{c_{\lambda}(j+1)})
\nonumber\\
&& \qquad \qquad 
-\1(\tau_{c_{\lambda}(j)} > \tau_{c_{\lambda}(j+1)})
p(\tau_{c_{\lambda}(j)} - \tau_{c_{\lambda}(i+1)},
x_{c_{\lambda}(j)} | x_{c_{\lambda}(j+1)} )
\Big\},
\label{eqn:star2}
\end{eqnarray}
where with (\ref{eqn:K1}) we have set 
\begin{eqnarray}
\cG_{s,t}(x,y)
&\equiv& \mbK_{\xi}(s,x;t,y)+\1(s>t) p(s-t,x|y)
\nonumber\\
&=& \int_{S} \xi(dv) p(s, x|v) \cM_{\xi}^{v}(t,y),
\quad (s,t) \in \cT^2, \quad (x,y) \in S^2.
\label{eqn:G1}
\end{eqnarray}

We will perform binomial expansions in (\ref{eqn:star2}).
In order to show the result, 
we introduce the following notations. 
For each cyclic permutation $\sc_{\lambda}$,
we consider a subset of $\{\sc_{\lambda}\}$,
$$
\cC(\sc_{\lambda}) =
\Big\{ c_{\lambda}(j) \in \{\sc_{\lambda}\} : 
\tau_{c_{\lambda}(j)} > \tau_{c_{\lambda}(j+1)} \Big\}.
$$
Choose $\bM_{\lambda}$ such that 
$\{\sc_{\lambda}\} \setminus \cC(\sc_{\lambda})
\subset \bM_{\lambda} \subset \{\sc_{\lambda}\}$,
and define $\bM_{\lambda}^{\rm c}=\{\sc_{\lambda}\} 
\setminus \bM_{\lambda}$.
(Then $\bM_{\lambda}^{\rm c} \subset \cC(\sc_{\lambda})$.)
Therefore if we put
\begin{eqnarray}
G(\sc_{\lambda}, \bM_{\lambda}) &=& 
\int_{S^{\{\sc_{\lambda}\}}} \prod_{j=1}^{q_{\lambda}} \ \bigg\{
dx_{c_{\lambda}(j)}  \  \chi_{{\tau}_{c_{\lambda}(j)}}(x_{c_{\lambda}(j)})
p( \tau_{c_{\lambda}(j)} - \tau_{c_{\lambda}(i+1)},
x_{c_{\lambda}(j)} | x_{c_{\lambda}(j+1)} )
^{\1(c_{\lambda}(j) \in \bM_{\lambda}^{\rm c})}
\nonumber\\
&& \qquad \qquad  \times
\cG_{{\tau}_{c_{\lambda}(j)},{\tau}_{c_{\lambda}(j+1)}}
(x_{c_{\lambda}(j)},x_{c_{\lambda}(j+1)})
^{\1(c_{\lambda}(j) \in \bM_{\lambda})}\bigg\},
\label{eqn:G}
\end{eqnarray}
the LHS of (\ref{eqn:E=det}) is expanded as
\begin{equation}
\frac{1}{\prod_{m=1}^{M} N_m !} 
\sum_{\sigma \in \cS_n} (-1)^{n-\ell(\sigma)}
\prod_{\lambda=1}^{\ell(\sigma)}
\sum_{\substack{\bM_{\lambda}: \\
\{\sc_{\lambda}\} \setminus \cC(\sc_{\lambda}) 
\subset \bM_{\lambda} \subset \{\sc_{\lambda}\}}}
(-1)^{\sharp \bM^{\rm c}_{\lambda}}
G(\sc_{\lambda}, \bM_{\lambda}).
\label{eqn:LHS2}
\end{equation}

Using only the entries of $\bM_{\lambda}$,
we can define a subcycle $\widehat{\sc_{\lambda}}$ of $\sc_{\lambda}$
uniquely as follows. 
For each $1 \leq j \leq q_{\lambda}$ with 
$c_{\lambda}(j) \in \bM_{\lambda}$, we define
\begin{eqnarray}
\oj &=& \min\{k > j : c_{\lambda}(k) \in \bM_{\lambda} \}
\nonumber\\
\uj &=& \max\{k < j : c_{\lambda}(k) \in \bM_{\lambda} \}. 
\label{eqn:ojuj}
\end{eqnarray} 
Since $\sc_{\lambda}$ is a cyclic permutation,
$\widehat{q_{\lambda}} \equiv \sharp \bM_{\lambda} \geq 1$.
Let $j_1=\min \{1 \leq i \leq q_{\lambda}: c_{\lambda}(j) 
\in \bM_{\lambda}\}$.
If $\widehat{q_{\lambda}} \geq 2$, define
$j_{k+1}=\oj_k, 1 \leq k \leq \widehat{q_{\lambda}}-1$.
Then $\widehat{\sc_{\lambda}} 
=(\widehat{c_{\lambda}}(1) \widehat{c_{\lambda}}(2) 
\cdots \widehat{c_{\lambda}}(\widehat{q_{\lambda}}))
\equiv (c_{\lambda}(j_1) c_{\lambda}(j_2) 
\cdots c(j_{\widehat{q_{\lambda}}}))$.

Moreover, 
we decompose the set $\bM_{\lambda}$ into $M$ subsets,
$\bM_{\lambda}=\bigcup_{m=1}^{M} \bJ_m^{\lambda}$,
by letting 
\begin{equation}
\bJ_m^{\lambda}=\bJ_m^{\lambda}(\sc_{\lambda}, \bM_{\lambda})=
\Big\{ c_{\lambda}(j) \in \bM_{\lambda} :
\uj < ^{\exists}k \leq j, \,
\mbox{s.t.} \, c_{\lambda}(k) \in \I^{(m)} \Big\}, \quad 1 \leq m \leq M.
\label{eqn:bJ1}
\end{equation}
By definition, if $c_{\lambda}(j) \in \bM_{\lambda}$ and
$\uj \leq j-2$, then for all $k$, s.t. 
$\uj < k < j$, 
we see $k \in \bM_{\lambda}^{\rm c} \subset \cC(\sc_{\lambda})$
and thus 
$\tau_{c_{\lambda}(k)} > \tau_{c_{\lambda}(k+1)} \geq \tau_{c_{\lambda}(j)}$.
Then in general
${\bJ_m^{\lambda} \cap \bJ_{m'}^{\lambda}} 
\not= \emptyset,
m \not= m'$, and
$\bJ_1^{\lambda}=\I_{N_1} \cap \bM_{\lambda}=
\I_{N_1} \cap \{\sc_{\lambda}\}$, 
$\bJ_m^{\lambda} \subset \I_{\sum_{k=1}^{m} N_k}$ for $2 \leq m \leq M$,
$\bJ_m^{\lambda} \cap \I^{(k)} \subset \bJ_{k}^{\lambda}$
for $1 \leq k < m \leq M$.

Now we prove the following lemma.
\begin{lem}
\label{thm:G4}
The quantity (\ref{eqn:G}) is equal to
\begin{equation}
\int_{S^{\bM_{\lambda}}} \prod_{j:c_{\lambda}(j) \in \bM_{\lambda}}
\xi(dv_{c_{\lambda}(j)}) \
\rE_{\v} \left[ 
\prod_{m=1}^{M} \prod_{j_m \in \bJ_m^{\lambda}}
\chi_{t_m}(V_{j_m}(t_m))
\prod_{j=1}^{\widehat{q_{\lambda}}}
\cM_{\xi}^{v_{\widehat{c_{\lambda}}(j)}} 
(T, V_{\widehat{c_{\lambda}}(j+1)}(T)) \right].
\label{eqn:G4}
\end{equation}
\end{lem}
\noindent{\it Proof of Lemma \ref{thm:G4}.}
We note that if we set
\begin{eqnarray}
&& F(\{x_{c_{\lambda}(k)}:c_{\lambda}(k) \in \bM_{\lambda}^{\rm c} \}) 
\nonumber\\
&& \qquad
= \int_{S^{\bM_{\lambda}}}
\prod_{j: c_{\lambda}(j) \in \bM_{\lambda}}\Big\{ dx_{c_{\lambda}(j)} 
\ \chi_{{\tau}_{c_{\lambda}(j)}}(x_{c_{\lambda}(j)})
\cG_{{\tau}_{c_{\lambda}(j)}, {\tau}_{c_{\lambda}(j+1)}}
(x_{c_{\lambda}(j)},x_{c_{\lambda}(j+1)}) \Big\}
\nonumber\\
&& \qquad \qquad \times
\prod_{k: c_{\lambda}(k) \in \bM_{\lambda}^{\rm c}}
p( \tau_{c_{\lambda}(k)} - \tau_{c_{\lambda}(k+1)},
x_{c_{\lambda}(k)} | x_{c_{\lambda}(k+1)} ),
\label{eqn:F}
\end{eqnarray}
which is the integral only over $S^{\bM_{\lambda}}$, 
then (\ref{eqn:G}) is obtained by performing the integral of it
also over $S^{\bM^{\rm c}_{\lambda}}=S^{\{\sc_{\lambda}\}}
\setminus S^{\bM_{\lambda}}$, 
\begin{equation}
G(\sc_{\lambda}, \bM_{\lambda}) =
\int_{S^{\bM_{\lambda}^{\rm c}}} 
\prod_{k:c_{\lambda}(k) \in \bM_{\lambda}^{\rm c}} \Big\{ dx_{c_{\lambda}(k)} 
\chi_{{\tau}_{c_{\lambda}(k)}}(x_{c_{\lambda}(k)}) \Big\}
F(\{x_{c_{\lambda}(k)}:c_{\lambda}(k) \in \bM_{\lambda}^{\rm c} \}).
\label{eqn:G2}
\end{equation}
In (\ref{eqn:F}), use the definition 
(\ref{eqn:G1}) for
$\cG_{{\tau}_{c_{\lambda}(j)}, {\tau}_{c_{\lambda}(j+1)}}
(x_{c_{\lambda}(j)},x_{c_{\lambda}(j+1)})$
by putting the integral variables
to be $v=v_{c_{\lambda}(j)}$. We obtain
\begin{eqnarray}
&& 
F(\{x_{c_{\lambda}(k)}:c_{\lambda}(k) \in \bM_{\lambda}^{\rm c}\}) 
\nonumber\\
&=& \int_{S^{ \bM_{\lambda}}}
\prod_{j: c_{\lambda}(j) \in \bM_{\lambda}} \xi (dv_{c_{\lambda}(j)})
\int_{S^{ \bM_{\lambda}}}
\prod_{j: c_{\lambda}(j) \in \bM_{\lambda}} 
\Big\{dx_{c_{\lambda}(j)}
p(\tau_{c_{\lambda}(j)}, x_{c_{\lambda}(j)}|v_{c_{\lambda}(j)})
\chi_{{\tau}_{c_{\lambda}(j)}}(x_{c_{\lambda}(j)}) \Big\}
\nonumber\\
&\times&
\prod_{j:c_{\lambda}(j) \in \bM_{\lambda}} 
\cM_{\xi}^{v_{c_{\lambda}(j)}}({\tau}_{c_{\lambda}(j+1)}, x_{c_{\lambda}(j+1)})
\prod_{k: c_{\lambda}(k) \in \bM_{\lambda}^{\rm c}} 
p({\tau}_{c_{\lambda}(k)}-{\tau}_{c_{\lambda}(k+1)}, 
x_{c_{\lambda}(k)}| x_{c_{\lambda}(k+1)})
\nonumber\\
&=&
\int_{S^{ \bM_{\lambda}}}
\prod_{j: c_{\lambda}(j) \in \bM_{\lambda}} \xi (dv_{c_{\lambda}(j)})
\rE_{\v} \Bigg[ 
\prod_{j:c_{\lambda}(j) \in \bM_{\lambda}} 
\bigg\{ \chi_{{\tau}_{c_{\lambda}(j)}}(V_{c_{\lambda}(j)}({\tau}_{c_{\lambda}(j)}))
\nonumber\\
&& \qquad \qquad \qquad \times
\cM_{\xi}^{v_{c_{\lambda}(j)}}({\tau}_{c_{\lambda}(j+1)}, V_{c_{\lambda}(j+1)}
({\tau}_{c_{\lambda}(j+1)}))^{\1(c_{\lambda}(j+1) \in \bM_{\lambda})}
\nonumber\\
&& \qquad \qquad \qquad  \times 
\cM_{\xi}^{v_{c_{\lambda}(j)}}
({\tau}_{c_{\lambda}(j+1)}, x_{c_{\lambda}(j+1)})^{\1(c_{\lambda}(j+1) 
\in \bM_{\lambda}^{\rm c})}\bigg\}
\nonumber\\
&& \quad \times
\prod_{k: c_{\lambda}(k) \in \bM_{\lambda}^{\rm c}} 
\bigg\{ p({\tau}_{c_{\lambda}(k)}-{\tau}_{c_{\lambda}(k+1)}, 
x_{c_{\lambda}(k)} | V_{c_{\lambda}(k+1)}({\tau}_{c_{\lambda}(k+1)})
)^{\1(c_{\lambda}(k+1) \in \bM_{\lambda})} \nonumber\\
&& \qquad \qquad \qquad \qquad \times 
p({\tau}_{c_{\lambda}(k)}-{\tau}_{c_{\lambda}(k+1)}, 
x_{c_{\lambda}(k)} | x_{c_{\lambda}(k+1)})
^{\1(c_{\lambda}(k+1) \in \bM_{\lambda}^{\rm c})} \bigg\} \Bigg].
\nonumber
\end{eqnarray}
Using Fubini's theorem, (\ref{eqn:G2}) is given by
\begin{eqnarray}
&& \int_{S^{\bM_{\lambda}}}
\prod_{j:c_{\lambda}(i) \in \bM_{\lambda}} \xi(dv_{c_{\lambda}(j)})
\rE_{\v} \Bigg[ 
\prod_{j:c_{\lambda}(j) \in \bM_{\lambda}} 
\chi_{{\tau}_{c_{\lambda}(j)}}(V_{c_{\lambda}(j)}({\tau}_{c_{\lambda}(j)}))
\nonumber\\
&& \qquad \qquad \times 
\prod_{j: c_{\lambda}(j), c_{\lambda}(j+1) \in \bM_{\lambda}}
\cM_{\xi}^{v_{c_{\lambda}(j)}}({\tau}_{c_{\lambda}(j+1)}, 
V_{c_{\lambda}(j+1)}({\tau}_{c_{\lambda}(j+1)}))
\nonumber\\
&& \quad \times 
\int_{S^{\bM_{\lambda}^{\rm c}}} 
\prod_{k: c_{\lambda}(k) \in \bM_{\lambda}^{\rm c}} 
\Big\{ dx_{c_{\lambda}(k)} \chi_{\tau_{c_{\lambda}(k)}}
(x_{c_{\lambda}(k)}) \Big\}
\nonumber\\
&& \qquad \qquad \times 
\prod_{k: c_{\lambda}(k) \in 
\bM_{\lambda}^{\rm c}, c_{\lambda}(k+1) \in \bM_{\lambda}}
p( \tau_{c_{\lambda}(k)} - \tau_{c_{\lambda}(k+1)},
x_{c_{\lambda}(k)} | V_{c_{\lambda}(k+1)}(\tau_{c_{\lambda}(k+1)}))
\nonumber\\
&& \qquad \qquad \times
\prod_{k: c_{\lambda}(k), c_{\lambda}(k+1) \in \bM_{\lambda}^{\rm c}}
p({\tau}_{c_{\lambda}(k)}-{\tau}_{c_{\lambda}(k+1)}, 
x_{c_{\lambda}(k)} | x_{c_{\lambda}(k+1)})
\nonumber\\
&& \qquad \qquad \times
\prod_{j: c_{\lambda}(j) \in \bM_{\lambda}, 
c_{\lambda}(j+1) \in \bM_{\lambda}^{\rm c}}
\cM_{\xi}^{v_{c_{\lambda}(j)}}
(\tau_{c_{\lambda}(j+1)}, x_{c_{\lambda}(j+1)})
\Bigg].
\label{eqn:G3}
\end{eqnarray}
We perform integration over $x_{c_{\lambda}(k)}$'s
for $c_{\lambda}(k) \in \bM^{\rm c}_{\lambda}$
before taking the expectation $\rE_{\v}$.
That is, integrals over $x_{c_{\lambda}(k)}$'s
with indices in intervals $\uj < k < j$
for all $j$, s.t. $c_{\lambda}(j) \in \bM_{\lambda}$
are done.
For each $j$, s.t. $c_{\lambda}(j) \in \bM_{\lambda}$,
if $\uj < j-1$,
\begin{eqnarray}
&& \chi_{{\tau}_{c_{\lambda}(j)}}(V_{c_{\lambda}(j)}({\tau}_{c_{\lambda}(j)}))
\Big\{ \prod_{k=\uj+1}^{j-1}
\int_{S} dx_{c_{\lambda}(k)} \chi_{\tau_{c_{\lambda}(k)}}
(x_{c_{\lambda}(k)}) 
\Big\}
\nonumber\\
&& \qquad \quad \times
p(\tau_{c_{\lambda}(j-1)}-\tau_{c_{\lambda}(j)}, 
x_{c_{\lambda}(j-1)} |V_{c_{\lambda}(j)}(\tau_{c_{\lambda}(j)}))
\nonumber\\
&& \qquad \quad \times
\prod_{l=\uj+2}^{j-1} 
p(\tau_{c_{\lambda}(l-1)}-\tau_{c_{\lambda}(l)},
x_{c_{\lambda}(l-1)} | x_{c_{\lambda}(l)})
\cM_{\xi}^{v_{c_{\lambda}(\, \uj \,)}}
(\tau_{c_{\lambda}(\uj+1)}, x_{c_{\lambda}(\uj+1)}) 
\nonumber
\end{eqnarray}
coincides with the conditional expectation of
\begin{eqnarray}
&&\prod_{k=\uj+1}^{j} \chi_{\tau_{c_{\lambda}(k)}}
(V_{c_{\lambda}(j)}(\tau_{c_{\lambda}(k)}))
\cM_{\xi}^{v_{c_{\lambda}(\, \uj \,)}}
(\tau_{c_{\lambda}(\uj+1)}, V_{c_{\lambda}(j)}(\tau_{c_{\lambda}(\uj+1)}))
\nonumber
\end{eqnarray}
w.r.t. $\rE_{\v}[\, \cdot \,|V_{c_\lambda(j)}(\tau_{c_{\lambda}(j)})]$.
Since 
$$
\prod_{j:c_{\lambda}(j) \in \bM_{\lambda}}
\cM_{\xi}^{v_{c_{\lambda}(\, \uj \,)}} 
(\tau_{c_{\lambda}(\uj+1)} , V_{c_{\lambda}(j)} (\tau_{c_{\lambda}(\uj+1)} ))
= \prod_{j:c_{\lambda}(j) \in \bM_{\lambda}}
\cM_{\xi}^{v_{c_{\lambda}(j)}}
(\tau_{c_{\lambda}(j+1)}, V_{c_{\lambda}(\, \oj \,)}(\tau_{c_{\lambda}(j+1)})),
$$
(\ref{eqn:G3}) is equal to 
\begin{eqnarray}
&& \int_{S^{\bM_{\lambda}}} 
\prod_{j: c_{\lambda}(j) \in \bM_{\lambda}} \xi(dv_{c_{\lambda}(j)})
\nonumber\\
&& \times
\rE_{\v} \left[ 
\prod_{j: c_{\lambda}(j) \in \bM_{\lambda}}
\left\{\prod_{k=\uj+1}^{j}
\chi_{\tau_{c_{\lambda}(k)}}(V_{c_{\lambda}(j)}(\tau_{c_{\lambda}(k)})) 
\cM_{\xi}^{v_{c_{\lambda}(j)}}
(\tau_{c_{\lambda}(j+1)}, V_{c_{\lambda}(\, \oj \,)}(\tau_{c_{\lambda}(j+1)}))
\right\} \right]. 
\nonumber
\end{eqnarray}
Then, by definition (\ref{eqn:bJ1}), 
we arrive at the expression (\ref{eqn:G4})
of $G(\sc_{\lambda}, \bM_{\lambda})$,
if we use the martingale property of $\cM_{\xi}^u$.
\qed
\vskip 0.3cm

Let $\bM \equiv \bigcup_{\lambda=1}^{\ell(\sigma)}
\bM_{\lambda}$. 
Since
$n-\sum_{\lambda=1}^{\ell(\sigma)} \sharp \bM_{\lambda}^{\rm c}
=\sharp \bM$,
the LHS of (\ref{eqn:E=det}), which is written above as 
(\ref{eqn:LHS2}) with Lemma \ref{thm:G4}, becomes now
\begin{eqnarray}
&& \frac{1}{\prod_{m=1}^{M} N_m !} 
\sum_{\sigma \in \cS_n} 
\sum_{\substack{\bM : \\
\I_n \setminus \bigcup_{\lambda=1}^{\ell(\sigma)}
\cC(\sc_{\lambda}) \subset \bM \subset \I_n}}
(-1)^{\sharp \bM -\ell(\sigma)}
\int_{S^{\bM}}
\prod_{\lambda=1}^{\ell(\sigma)}
\prod_{j: c_{\lambda}(j) \in \bM_{\lambda}}
\xi(dv_{c_{\lambda}(j)})
\nonumber\\
&& \qquad \times 
\rE_{\v} \left[ \prod_{\lambda=1}^{\ell(\sigma)}
\left\{ \prod_{m=1}^{M}
\prod_{j_m \in \bJ_m^{\lambda}} 
\chi_{t_m}(V_{j_m}(t_m)) 
\prod_{j=1}^{\widehat{q_{\lambda}}}
\cM_{\xi}^{v_{\widehat{c_{\lambda}}(j)}}
(T, V_{\widehat{c_{\lambda}}(j+1)}(T)) 
\right\} \right].
\label{eqn:E=det2}
\end{eqnarray}

We define 
$$
\widehat{\sigma} \equiv 
\widehat{\sc_1} \widehat{\sc_2} \cdots
\widehat{\sc_{\ell(\sigma)}}
$$
and  
$$
\bJ_m \equiv \bigcup_{\lambda=1}^{\ell(\sigma)} \bJ_m^{\lambda},
\quad 1 \leq m \leq M.
$$
Note that $\ell(\widehat{\sigma})=\ell(\sigma)$.
The obtained $(\bJ_m)_{m=1}^{M}$'s 
form a collection of series of index sets
satisfying the following conditions,
which we write as ${\cal J}(\{N_m\}_{m=1}^{M})$:
\begin{eqnarray}
&&
\bJ_1=\I_{N_1},  \quad
\bJ_m \subset \I_{\sum_{k=1}^{m} N_{k}}
\quad \mbox{for} \quad 2 \leq m \leq M,
\nonumber\\
&& 
\bJ_m \cap \I^{(k)} \subset \bJ_{k} 
\quad \mbox{for} \quad
1 \leq k < m \leq M,  \quad \mbox{and}
\nonumber\\
&& 
\sharp \bJ_m=N_m \quad \mbox{for} \quad 1 \leq m \leq M.
\label{eqn:CJ}
\end{eqnarray}

For each $(\bJ_m)_{m=1}^{M} \in {\cal J}(\{N_m\}_{m=1}^{M})$,
we put 
$$
A_1=0 \quad \mbox{and} \quad
A_m=\sharp \left(\bJ_m \cap \I_{\sum_{k=1}^{m-1}N_{k}} \right)
=\sharp \left( \bJ_m \cap \bigcup_{k=1}^{m-1} \bJ_k \right), \quad 
2 \leq m \leq M.
$$
Then, if we put $\bM=\bigcup_{m=1}^{M} \bJ_m$,
$\sharp \bM=\sum_{m=1}^{M}(N_m-A_m)$,
which means that
from the original index set $\I_n=\bigcup_{m=1}^{M} \I^{(m)}$
with $\sharp \I^{(m)}=N_m, 1 \leq m \leq M$,
we obtain a subset $\bM$ by eliminating
$A_m$ elements at each level $1 \leq m \leq M$.
By this reduction, we obtain
$\widehat{\sigma} \in \cS(\bM)$ from $\sigma \in \cS_n$.
It implies that, for all $\widehat{\sigma} \in \cS(\bM)$,
the number of $\sigma$'s in $\cS_n$ which give
the same $\widehat{\sigma}$ and
$(\bJ_m)_{m=1}^{M}$ by this reduction is given by
$\prod_{m=1}^{M} A_m !$, where $0! \equiv 1$.
Then (\ref{eqn:E=det2}) is equal to 
\begin{eqnarray}
&& 
\sum_{\substack{\bM: \\
\max_m\{N_m\} \leq \sharp \bM \leq N}}
\sum_{\substack{(\bJ_m)_{m=1}^{M} \subset
{\cal J}(\{N_m\}_{m=1}^{M}): \\
\bigcup_{m=1}^{M} \bJ_m=\bM
}}
\frac{\prod_{m=1}^{M} A_m !}{\prod_{m=1}^{M} N_m !}
\sum_{\widehat{\sigma} \in \cS(\bM)}
(-1)^{\sharp \bM-\ell(\widehat{\sigma})}
\nonumber\\
&& \qquad \times
\sharp \bM !
\int_{\W_{\sharp \bM}} \xi^{\otimes \bM} (d \v)
\rE_{\v} \left[
\prod_{m=1}^{M} \prod_{j_m \in \bJ_m}
\chi_{t_m}(V_{j_m}(t_m))
\prod_{\lambda=1}^{\ell(\widehat{\sigma})}
\prod_{j=1}^{\widehat{q_{\lambda}}}
\cM_{\xi}^{v_{\widehat{ c_{\lambda} }(j)}}
(T, V_{\widehat{c_{\lambda}}(j+1)}(T)) \right]
\nonumber\\
&&=
\sum_{\substack{\bM: \\
\max_m\{N_m\} \leq \sharp \bM \leq N}}
\sum_{\substack{(\bJ_m)_{m=1}^{M} \subset
{\cal J}(\{N_m\}_{m=1}^{M}): \\
\bigcup_{m=1}^{M} \bJ_m=\bM
}}
\sharp \bM ! \prod_{m=1}^{M} \frac{A_m!}{N_m!}
\nonumber\\
&& \qquad \times
\int_{\W_{\sharp \bM}} \xi^{\otimes \bM} (d \v)
\rE_{\v} \left[
\prod_{m=1}^{M} \prod_{j_m \in \bJ_m}
\chi_{t_m}(V_{j_m}(t_m))
\det_{j,k \in \bM}
\left[ \cM_{\xi}^{v_k} (T, V_{j}(T)) \right] \right].
\label{eqn:AAA}
\end{eqnarray}
Assume $1 \leq p \leq N$, 
$0 \leq A_m \leq N_m, 2 \leq m \leq M$ and set $A_1 =0$.
Consider
\begin{eqnarray}
&& \Lambda_1 =\left\{(\bJ_m)_{m=1}^{M} \subset
{\cal J}(\{N_{m}\}_{m=1}^{M}): 
\sharp \left(\bigcup_{m=1}^{M} \bJ_m \right) =p,
\right. 
\nonumber\\
&& \hskip 6cm \left.
\sharp \left(\bJ_m \cap \bigcup_{k=1}^{m-1} \bJ_k \right)
=A_m, 2 \leq m \leq M \right\},
\nonumber\\
&& \Lambda_2 = \left\{(\J_m)_{m=1}^{M} :
\sharp \J_m=N_m, 1 \leq m \leq M,
\bigcup_{m=1}^{M} \J_m= \I_{p},
\right. 
\nonumber\\
&& \hskip 6cm \left.
\sharp \left(\J_m \cap \bigcup_{k=1}^{m-1} \J_k \right)
=A_m, 2 \leq m \leq M \right\}.
\nonumber
\end{eqnarray}
Since $V_j(\cdot)$'s are i.i.d. in $\rP_{\v}$, 
the integral in (\ref{eqn:AAA}) has the same value
for all $(\bJ_m)_{m=1}^{M} \in \Lambda_1$
with $\bigcup_{m=1}^{M} \bJ_{m}=\bM$
and it is also equal to
$$
\int_{\W^{\rm A}_{p}} \xi^{\otimes p} (d \v)
\rE_{\v} \left[
\prod_{m=1}^{M} \prod_{j_m \in \J_m}
\chi_{t_m}(V_{j_m}(t_m))
\det_{j,k \in \I_{p}}
\left[ \cM_{\xi}^{v_k} (T, V_j(T)) \right] \right]
$$
for all $(\J_{m})_{m=1}^{M} \in \Lambda_2$.

In $\Lambda_1$, for each $2 \leq m \leq M$,
$A_m$ elements in $\bJ_m$ are chosen from
$\bigcup_{k=1}^{m-1} \bJ_k$,
in which $\sharp (\bigcup_{k=1}^{m-1}
\bJ_{k})=\sum_{k=1}^{m-1}(N_k-A_k)$,
and the remaining $N_m-A_m$ elements in $\bJ_m$
are from $\I^{(m)}$ with $\sharp \I^{(m)}=N_m$.
Then
$$
\sharp \Lambda_1=\prod_{m=2}^{M} 
{\sum_{k=1}^{m-1} (N_k-A_k) \choose A_m}
{N_m \choose N_m-A_m}.
$$

In $\Lambda_2$, on the other hand,
$N_1$ elements in $\J_1$ is chosen from $\I_{p}$,
and then for each $2 \leq m \leq M$,
$A_m$ elements in $\J_m$ are chosen from
$\bigcup_{k=1}^{m-1} \J_k$ with
$\sharp (\bigcup_{k=1}^{m-1} \J_k)=
\sum_{k=1}^{m-1}(N_k-A_k)$
and the remaining $N_m-A_m$ elements in $\J_m$
are from $\I_{p} \setminus \bigcup_{k=1}^{m-1} \J_k$
with $\sharp(\I_{p} \setminus \bigcup_{k=1}^{m-1} \J_k)
=p-\sum_{k=1}^{m-1}(N_k-A_k)$.
Then
$$
\sharp \Lambda_2
={p \choose N_1} 
\prod_{m=2}^{M} 
{\sum_{k=1}^{m-1} (N_k-A_k) \choose A_m}
{p-\sum_{k=1}^{m-1} (N_k-A_k) \choose N_m-A_m}.
$$

Since $\sum_{m=1}^{M}(N_m-A_m)=p$,
we see 
$\sharp \Lambda_2/\sharp \Lambda_1=
p ! \prod_{m=1}^{M} A_m !/N_m!$.
Then (\ref{eqn:AAA}) is equal to the RHS of (\ref{eqn:E=det})
and the proof is completed.
\qed

\clearpage
\SSC{Polynomial Martingales \label{sec:poly_mar}}

For $n \in \N_0$, here we consider the monic
polynomials of degrees $n$ with time-dependent 
coefficients,
\begin{equation}
m_n(t, x)=x^n + \sum_{j=0}^{n-1} c_n^{(j)}(t) x^j, \quad t \geq 0
\label{eqn:poly_mar1}
\end{equation}
satisfying the conditions such that
\begin{equation}
m_n(0,x)=x^n,
\label{eqn:poly_mar2}
\end{equation}
and that, if we replace $x$ by the Markov process
$V(t), t \in \cT$, then they are local martingales.
Such polynomials $\{m_n(t,v)\}_{n \in \N_0}$
are called the {\it polynomial martingales}
associated with the process $V(\cdot)$.

\subsection{Brownian motion (BM) and Hermite polynomials
\label{sec:Hermite}}

Let $V(t)=B(t), t \in \cT=[0, \infty)$,
the one-dimensional standard Brownian motion (BM)
on $S=\R$.
The transition probability density is given by
\begin{equation}
p(t, y|x)
=  \left\{ \begin{array}{ll}
\displaystyle{
\frac{1}{\sqrt{2 \pi t}} e^{-(x-y)^2/2t}},
& \quad t>0, x, y \in \R \cr
& \cr
\delta(y-x),
& \quad t=0, x, y \in \R.
\end{array} \right.
\label{eqn:p_BM}
\end{equation}

For $n \in \N_0$, 
the Hermite polynomials of degrees $n \in \N_0$
are given by
\begin{equation}
H_n(x)=\sum_{j=0}^{[n/2]} (-1)^j 
\frac{n!}{j! (n-2j)!} (2x)^{n-2j}, 
\quad n \in \N_0,
\label{eqn:Hermite1}
\end{equation}
which solve the differential equation
\begin{equation}
y''-2x y'+2n y=0.
\label{eqn:Hermite2}
\end{equation}
The following is proved.

\begin{lem}
\label{thm:Hermite1}
The polynomials of $B(t), t \in [0, \infty)$, 
\begin{equation}
m_n(t, B(t))
=\left( \frac{t}{2} \right)^{n/2} 
H_n \left( \frac{B(t)}{\sqrt{2t}} \right),
\quad t \geq 0, \quad n \in \N_0,
\label{eqn:Hermite_m2}
\end{equation}
are all local martingales.
\end{lem}
\noindent{\it Proof.} \
$m_0(t, B(t))\equiv 1$.
By It\^o's formula, for $n \geq 1$,
\begin{eqnarray}
d m_n(t, B(t)) &=&
\left[ \frac{n}{2} \frac{t^{n/2-1}}{2^{n/2}}
H_n \left( \frac{B(t)}{\sqrt{2t}} \right)
+\left(\frac{t}{2}\right)^{n/2} 
H_n' \left( \frac{B(t)}{\sqrt{2t}}\right)
\left( - \frac{B(t)}{(2t)^{3/2}} \right) \right] dt
\nonumber\\
&& + \left(\frac{t}{2}\right)^{n/2}
H_n' \left( \frac{B(t)}{\sqrt{2t}} \right)
\frac{1}{\sqrt{2t}} dB(t)
+ \frac{1}{2} \left( \frac{t}{2} \right)^{n/2}
H_n'' \left( \frac{B(t)}{\sqrt{2t}} \right)
\frac{1}{2t} dt
\nonumber\\
&=& \frac{1}{2} \left( \frac{t}{2} \right)^{(n-1)/2}
H_n' \left( \frac{B(t)}{\sqrt{2t}} \right) dB(t)
+ A_n(t) dt.
\nonumber
\end{eqnarray}
Here we find
$$
A_n(t)= \frac{t^{n/2-1}}{2^{n/2+2}}
\left[ H_n'' \left( \frac{B(t)}{\sqrt{2t}}\right)
- \sqrt{\frac{2}{t}} B(t)
H_n' \left( \frac{B(t)}{\sqrt{2t}} \right)
+2n H_n \left( \frac{B(t)}{\sqrt{2t}} \right) \right]
= 0,
$$
for (\ref{eqn:Hermite2}).
Then $m_n(t,B(t))$ are given by
stochastic integrals
\begin{equation}
m_n(t, B(t))
= \frac{1}{2^{(n+1)/2}}
\int_0^t s^{(n-1)/2}
H_n' \left( \frac{B(s)}{\sqrt{2s}} \right) dB(s),
\quad t \geq 0, \quad n \geq 1.
\label{eqn:Hermite_m3}
\end{equation}
The proof is thus completed. \qed
\vskip 0.3cm
We call the polynomials
\begin{equation}
m_n(t, x)=
\left( \frac{t}{2} \right)^{n/2} 
H_n \left( \frac{x}{\sqrt{2t}} \right),
\quad t \geq 0, \quad n \in \N_0, 
\label{eqn:Hermite_m1}
\end{equation}
the {\it polynomial martingales
associated with} $B(t), t \in [0, \infty)$.

\vskip 0.3cm
\noindent{\bf Remark 1.} \,
It is obvious that
\begin{equation}
G_{\alpha}^{\rm BM}(t, B(t))
=e^{\alpha B(t)-t \alpha^2/2}
\label{eqn:GBM1}
\end{equation}
is martingale for any $\alpha \in \C$. 
It is known that
\begin{equation}
G_{\alpha}^{\rm BM}(t,x)
=\sum_{n=0}^{\infty} \left( \frac{t}{2} \right)^{n/2}
H_n \left( \frac{x}{\sqrt{2t}} \right) 
\frac{\alpha^n}{n!},
\label{eqn:GBM2}
\end{equation}
where $H_n, n \in \N_0$ are the Hermite polynomials
given by (\ref{eqn:Hermite1}).
Then Lemma \ref{thm:Hermite1} is immediately concluded,
if we confirm that $m_n(t,x), n \in \N_0$ are monic.

\subsection{Squared Bessel Processes (BESQ$^{(\nu)}$)
and Laguerre polynomials
\label{sec:Laguerre}}

Let $V(t)=R^{(\nu)}(t), t \in \cT=[0, \infty)$,
$\nu > -1$,
the squared Bessel process with index $\nu>-1$ (BESQ$^{(\nu)}$)
on $S=\R_+ =\{x \in \R: x \geq 0\}$.
For $\nu=D/2-1, D \in \N$,
$R^{(\nu)}(\cdot)$ can be defined as a sum of squares
of $D$ independent BMs, $B_j(\cdot), 1 \leq j \leq D$ such as
$R^{(\nu)}(t)=\sum_{j=1}^D B_j(t)^2, t \geq 0$.
For general $\nu > -1$, it is given by the solution
of the stochastic differential equation (SDE),
\begin{equation}
R^{(\nu)}(t)=\int_0^t 2 \sqrt{R^{(\nu)}(s)} d B(s)
+2(\nu+1) t, \quad t \geq 0,
\label{eqn:BESQ1}
\end{equation}
where $B(\cdot)$ is a BM, and, 
if $-1 < \nu < 0$, a reflection wall
is put at the origin.
The transition probability density is given by
\begin{equation}
p^{(\nu)}(t, y|x)
=  \left\{ \begin{array}{ll}
\displaystyle{
\frac{1}{2t} \left( \frac{y}{x} \right)^{\nu/2}
\exp \left( - \frac{x+y}{2t} \right)
I_{\nu} \left( \frac{\sqrt{xy}}{t} \right)},
& \quad t>0, x>0, y \in \R_+, \cr
\displaystyle{
\frac{y^{\nu}}{(2t)^{\nu+1} \Gamma(\nu+1)} e^{-y/2t}},
& \quad t >0, x=0, y \in \R_+, \cr
& \cr
\delta(y-x),
& \quad t=0, x, y \in \R_+,
\end{array} \right.
\label{eqn:p_BESQ}
\end{equation}
if $-1 < \nu < 0$, the origin is assumed to be reflecting.
Here $I_{\nu}(x)$ is the modified Bessel function
of the first kind defined by
\begin{equation}
I_{\nu}(x) = \sum_{n=0}^{\infty} 
\frac{1}{\Gamma(n+1) \Gamma(n+1+\nu)}
\left( \frac{x}{2} \right)^{2n+\nu}
\label{eqn:I1}
\end{equation}
with the Gamma function
$\Gamma(z) = \int_{0}^{\infty} e^{-u} u^{z-1} du, 
\ \Re u > 0$.

For $n \in \N_0$, 
the Laguerre polynomial of degree $n \in \N_0$ with 
index $\nu > -1$ 
is given by
\begin{equation}
L^{(\nu)}_n(x)=\sum_{j=0}^{n}(-1)^j 
\frac{\Gamma(n+\nu+1)}{\Gamma(\nu+j+1) (n-j)! j!} x^j,
\quad n \in \N_0,
\label{eqn:Laguerre1}
\end{equation}
which solve the differential equation
\begin{equation}
xy''+(\nu+1-x)y'+ny=0.
\label{eqn:Laguerre2}
\end{equation}

The following is derived.

\begin{lem}
\label{thm:Laguerre1}
For $\nu > -1$, 
\begin{equation}
m_n(t, x)=
(-1)^n n! (2t)^{n} L^{(\nu)}_n
\left( \frac{x}{2t} \right),
\quad t \geq 0, \quad n \in \N_0
\label{eqn:Laguerre_m1}
\end{equation}
are the polynomial martingales
associated with 
the BESQ$^{(\nu)}$, 
$R^{(\nu)}(t), \nu > -1, t \in [0, \infty)$.
\end{lem}
\noindent{\it Proof.} \
$m_0(t, R^{(\nu)}(t))\equiv 1$.
By (\ref{eqn:BESQ1}), the quadratic variation
of BESQ$^{(\nu)}$ is
$\langle R^{(\nu)} \rangle_t = 4 \int_0^t R^{(\nu)}(s) ds, 
t \geq 0, \nu > -1$.
Then, for $n \geq 1$, It\^o's formula gives
\begin{eqnarray}
d m_n(t, R^{(\nu)})
&=& (-1)^n n! \left[
n 2^n t^{n-1} L_n^{(\nu)} \left( \frac{R^{(\nu)}(t)}{2t} \right)
+(2t)^n {L_n^{(\nu)}}' \left( \frac{R^{(\nu)}(t)}{2t} \right)
\left( - \frac{R^{(\nu)}(t)}{2t^2} \right) \right]
\nonumber\\
&& +(-1)^n n! (2t)^{n}
{L_n^{(\nu)}}' \left( \frac{R^{(\nu)}(t)}{2t} \right)
\frac{1}{2t} \left\{
2 \sqrt{R^{(\nu)}(t)} dB(t)+2(\nu+1) dt \right\}
\nonumber\\
&& 
+\frac{1}{2} (-1)^n n! (2t)^n 
{L_n^{(\nu)}}''\left( \frac{R^{(\nu)}(t)}{2t} \right)
\frac{1}{(2t)^2} 4 R^{(\nu)}(t) dt
\nonumber\\
&=& (-1)^n n! 2^n t^{n-1} \sqrt{R^{(\nu)}(t)}
{L_n^{(\nu)}}' \left( \frac{R^{(\nu)}(t)}{2t} \right) dB(t)
+A_n^{(\nu)}(t) dt
\nonumber
\end{eqnarray}
with
\begin{eqnarray}
A_n^{(\nu)}(t) &=& (-1)^n n! 2^n t^{n-1}
\left[ \frac{R^{(\nu)}(t)}{2t}
{L_n^{(\nu)}}'' \left( \frac{R^{(\nu)}(t)}{2t} \right)
\right.
\nonumber\\
&& \qquad \qquad \left.
+ \left( \nu+1-\frac{R^{(\nu)}(t)}{2t} \right)
{L_n^{(\nu)}}' \left( \frac{R^{(\nu)}(t)}{2t} \right)
+ n L_n^{(\nu)} \left( \frac{R^{(\nu)}(t)}{2t} \right) \right].
\nonumber
\end{eqnarray}
For (\ref{eqn:Laguerre2}),
$A_n^{(\nu)}(t) = 0, n \geq 1$.
Then $m_n(t, R^{(\nu)}(t)), \nu > -1$ are given by
stochastic integrals
\begin{equation}
m_n(t, R^{(\nu)}(t))
= (-1)^n n! 2^n 
\int_0^{t} s^{n-1} \sqrt{R^{(\nu)}(s)}
{L_n^{(\nu)}}' \left( \frac{R^{(\nu)}(s)}{2s} \right) dB(s),
\quad t \geq 0, \quad n \geq 1.
\label{eqn:Laguerre_m2}
\end{equation}
The proof is thus completed. \qed

\vskip 0.3cm
\noindent{\bf Remark 2.} \,
For $\alpha \in \C, t \geq 0, x \in \R, \nu > -1$, let
\begin{equation}
G_{\alpha}^{(\nu)}(t,x)
=\frac{e^{\alpha x/(1+2t \alpha)}}{(1+2t \alpha)^{\nu+1}}.
\label{eqn:Gnu1}
\end{equation}
By It\^o's formula, we can see 
\begin{equation}
d G_{\alpha}^{(\nu)}(t, R^{(\nu)}(t))
= \frac{2 \alpha \sqrt{R^{(\nu)}(t)}}{1+2 t \alpha}
G_{\alpha}^{(\nu)}(t, R^{(\nu)}(t)) d B(t),
\label{eqn:Gnu2}
\end{equation}
that is, $G_{\alpha}^{(\nu)}(t,R^{(\nu)}(t))$ is a local martingale
for any $\alpha \in \C$.
It is known that
\begin{equation}
G_{\alpha}^{(\nu)}(t,x)
=\sum_{n=0}^{\infty} (-1)^n n! (2t)^n
L^{(\nu)}_n \left( \frac{x}{2t} \right)
\frac{\alpha^n}{n!}.
\label{eqn:Gnu3}
\end{equation}
Then Lemma \ref{thm:Laguerre1}
is obtained.

\subsection{Random walk (RW)
and Fujita's polynomials
\label{sec:Fujita}}

Let $\Z$ be a set of all integers and $N \in \N \equiv \{1,2, \dots \}$.
Let $V(t), t \in \N_0$ be a one-dimensional, 
simple and symmetric RW on $S=\Z$
starting from 0 at time $t=0$,
\begin{equation}
V(t)=\zeta(1)+\zeta(2)+ \cdots + \zeta(t),
\quad t \in \N,
\label{eqn:RWb1}
\end{equation}
where $\{\zeta(t) : t \in \N\}$ are i.i.d. with
\begin{equation}
\rP[\zeta(1)=1]=\frac{1}{2}, \quad \rP[\zeta(1)=-1]=\frac{1}{2}.
\label{eqn:RWb2}
\end{equation}
The following {\it discrete It\^o's formula}
was given by Fujita \cite{Fuj02,Fuj08,Fuj08b}.

\begin{lem}
\label{sec:Fujita1}
For any $f: \N_0 \times \Z \to \R$ and
any $t \in \N_0$,
\begin{eqnarray}
&& f(t+1,V(t+1))-f(t,V(t))
\nonumber\\
&& \quad = \frac{1}{2} \Big[
f(t+1, V(t)+1)-f(t+1, V(t)-1) \Big]
(V(t+1)-V(t))
\nonumber\\
&& \qquad + \frac{1}{2} \Big[ f(t+1, V(t)+1)- 2 f(t+1, V(t))+ f(t+1, V(t)-1) \Big]
\nonumber\\
&& \qquad +f(t+1, V(t))-f(t, V(t)).
\label{eqn:Ito_Fujita}
\end{eqnarray}
\end{lem}
\vskip 0.3cm

We perform the Esscher transform 
with parameter $\alpha \in \R$,
$V(\cdot) \to \widetilde{V}_{\alpha}(\cdot)$ as
$$
\widetilde{V}_{\alpha}(t)
=\frac{e^{\alpha V(t)}}{\rE[e^{\alpha V(t)}]},
\quad t \in \N_0.
$$
By (\ref{eqn:RWb2}), 
$\rE[e^{\alpha \zeta(1)}]=(e^{\alpha}+e^{-\alpha})/2=\cosh \alpha$,
then we have
\begin{equation}
\widetilde{V}_{\alpha}(t)=
G_{\alpha}(t, V(t))
\label{eqn:tV2}
\end{equation}
with
\begin{equation}
G_{\alpha}(t,x)=\frac{e^{\alpha x}}
{(\cosh \alpha)^t},
\quad t \in \N_0, \quad x \in \Z.
\label{eqn:GG1}
\end{equation}
If we set $f=G_{\alpha}$ in (\ref{eqn:Ito_Fujita}),
the second and third terms in the RHS vanish.
Then
\begin{eqnarray}
&& G_{\alpha}(t+1, V(t+1))-G_{\alpha}(t,V(t))
\nonumber\\
&& \qquad
= \frac{1}{2} \Big[
G_{\alpha}(t+1, V(t)+1) - G_{\alpha}(t+1, V(t)-1) \Big]
\zeta(t+1),
\nonumber
\end{eqnarray}
which implies that $G_{\alpha}(t, V(t))$ is
$\{\zeta(1), \zeta(2), \dots, \zeta(t)\}$-martingale
for any $\alpha \in \R$ \cite{Fuj02,Fuj08,Fuj08b}.
From now on, we simply say
`$G_{\alpha}(t, V(t))$ is martingale'
in such a situation.

Expansion of (\ref{eqn:GG1}) with respect to $\alpha$
around $\alpha=0$ 
\begin{equation}
G_{\alpha}(t,x)=\sum_{n=0}^{\infty} m_n(t,x)
\frac{\alpha^n}{n!},
\label{eqn:GG3}
\end{equation}
determines a series of monic polynomials
of degrees $n$ studied by Fujita in \cite{Fuj02,Fuj08b}
\begin{equation}
m_n(t,x)=x^n + \sum_{j=1}^{n-1} c_n^{(j)} x^j,
\quad n \in \N_0, 
\label{eqn:m1}
\end{equation}
such that
\begin{eqnarray}
&& c_n^{(j)}(0)=0, \quad 1 \leq j \leq n-1, \quad \mbox{and}
\nonumber\\
&&
\mbox{$m_n(t, V(t))$ is martingale}, \quad
t \in \N_0.
\nonumber
\end{eqnarray}
For example,
\begin{eqnarray}
m_0(t,x) &=& 1,
\nonumber\\
m_1(t,x) &=& x,
\nonumber\\
m_2(t,x) &=& x^2-t,
\nonumber\\
m_3(t,x) &=& x^3-3 t x,
\nonumber\\
m_4(t,x) &=& x^4-6 t x^2+(3t+2)t,
\nonumber\\
m_5(t,x) &=& x^5-10 t x^3 + 5 (3t+2) t x.
\nonumber
\end{eqnarray}
They satisfy the recurrence relations
$$
m_n(t,x)=\frac{1}{2}
[ m_n(t+1,x+1)+m_n(t+1,x-1) ],
\quad n \in \N_0.
$$
We call $m_n(t,x), n \in \N$, 
{\it Fujita's polynomials} and 
$m_n(t,V(t)), n \in \N_0$, 
{\it Fujita's polynomial martingales} 
for the simple and symmetric RW \cite{Fuj02,Fuj08b}.

\vskip 0.3cm
\noindent{\bf Remark 3.} \,
The Esscher transform with parameter $\alpha$ 
for BM, $B(t), t \geq 0$ is given by
$$
\widetilde{B}_{\alpha}(t)
= G^{\rm BM}_{\alpha}(t,B(t))
$$
with
\begin{equation}
G^{\rm BM}_{\alpha}(t,x)
= \frac{e^{\alpha x}}{\rE[e^{\alpha B(t)}]}
=\frac{e^{\alpha x}}
{\displaystyle{ \int_{-\infty}^{\infty} dx
e^{\alpha x} p^{\rm BM}(t, x|0)
 }}
= e^{\alpha x - \alpha^2 t/2},
\label{eqn:GBM2b}
\end{equation}
where
$p^{\rm BM}(t,y|x)$ 
is the transition probability density of BM (\ref{eqn:p_BM}).
This is nothing but (\ref{eqn:GBM1}).

\clearpage
\SSC{Integral Transforms 
and Complex-Process Representations (CPR)
for Polynomial Martingales 
\label{sec:integral}}

For each set of polynomial martingales
$\{m_n(t,x): n \in \N_0\}$ associated with
the Markov process $V(t)$, here we want to
determine the integral transform 
of an integrable function $f$
of the form
\begin{equation}
\sfM[f(W) |(t, x)]
=\int_{S} dw \, q(t, w|x) f(w),
\label{eqn:int1}
\end{equation}
such that it
satisfies the equalities
\begin{equation}
m_n(t, x)=\sfM \left[ \left.
(c W)^n \right|(t, x) \right],
\quad \forall n \in \N_0, \quad \forall t \in \cT
\label{eqn:int2}
\end{equation}
with some constant $c \in \C$.
If so, 
given any polynomial $f$, 
\begin{eqnarray}
\label{eqn:int3}
&&
\mbox{$\sfM [f(cW)|(t, V(t))], t \geq 0$ 
is a local martingale, and}
\\
\label{eqn:int4}
&& \sfM [f(cW) |(0, V(0))]
=f(V(0)).
\end{eqnarray}
Note that by setting $f \equiv 1$ in (\ref{eqn:int4}) we have
$\sfM[1|(t, V(t))] \equiv 1, t \geq 0$.

\subsection{BM
\label{sec:int_BM}}

For BM, $V(t)=B(t), t \in [0, \infty)$,
we set
\begin{equation}
c=i.
\label{eqn:c_BM}
\end{equation}
and 
\begin{eqnarray}
q(t,y|x) &=& p(t, y| c^{-1} x)
\nonumber\\
&=&  \left\{ \begin{array}{ll}
\displaystyle{
\frac{1}{\sqrt{2 \pi t}} e^{-(ix+y)^2/2t}},
& \quad t>0, x, y \in \R \cr
& \cr
\delta(y-x),
& \quad t=0, x, y \in \R,
\end{array} \right.
\label{eqn:q_BM}
\end{eqnarray}
Then we can prove the following.

\begin{lem}
\label{thm:int_BM}
With (\ref{eqn:c_BM}) and (\ref{eqn:q_BM}),
(\ref{eqn:int2}) are satisfied.
\end{lem}
{\it Proof.} \quad
Since $p(t,\cdot|x)$ solves the diffusion equation,
$q(t, \cdot|x)=p(t, \cdot|c^{-1}x)$  satisfies
$\partial q/\partial t
=c^{-2} (1/2) \partial^2 q/\partial x^2$.
Then It\^o's formula implies
\begin{eqnarray}
d \sfM[f(W)| (t, B(t))]
&=& \int_{\R} dw \, d q(t, w|B(t)) f(w)
\nonumber\\
&=& \left[ \int_{\R} dw \,
\left\{ \frac{\partial q}{\partial t}(t, w|B(t))
+\frac{1}{2} \frac{\partial^2 q}{\partial x^2}(t,w|B(t)) \right\} f(w) \right]dt
\nonumber\\
&& \quad
+ \left\{ \int_{\R} dw \, \frac{\partial q}{\partial x}(t,w|B(t)) f(w) \right\} dB(t)
\nonumber\\
&=& \left[ \int_{\R} dw \, \frac{\partial q}{\partial x}(t,w|B(t)) f(w) \right] dB(t),
\nonumber
\end{eqnarray}
if $c^{-2}=-1 \Leftrightarrow c=\pm i$.
Therefore, the assignment (\ref{eqn:c_BM}) of the value $c$
guarantees that $\sfM[f(W)|(t,B(t))]$ is a local martingale.
The Hermite polynomials have the following integral representations,
$n \in \N_0, x \in \R$ (for instance, see Eq.(6.1.4) in \cite{AAR99}),
\begin{eqnarray}
H_n(x) 
\label{eqn:Hnint2}
&=& \frac{2^n}{\sqrt{\pi}}
\int_{-\infty}^{\infty} du \,
e^{-(ix+u)^2} (i u)^n \\
\label{eqn:Hnint1}
&=& \frac{2^n}{\sqrt{\pi}} 
\int_{-\infty}^{\infty} du \, 
e^{-u^2} (x+iu)^n.
\end{eqnarray}
The formula (\ref{eqn:Hnint2})
gives (\ref{eqn:int2}) with appropriate
change of variables.
\qed
\vskip 0.5cm
\noindent{\bf Remark 4.} \,
The function $G_{\alpha}^{\rm BM}(t, B(t)), t \geq 0$
given by (\ref{eqn:GBM1}) is martingale
for any $\alpha \in \R$.
Then its Fourier transform with respect to $\alpha$,
\begin{equation}
q(t,w|B(t))=\frac{1}{2\pi} \int_{\R} d \alpha \,
e^{-i \alpha w} G_{\alpha}^{\rm BM}(t, B(t))
\label{eqn:qBM1}
\end{equation}
is also martingale for any $w \in \R$.
We find that
\begin{eqnarray}
q(t,w|x) &=& \frac{1}{2 \pi} \int_{\R} d \alpha \,
e^{-i \alpha w} e^{\alpha x - t \alpha^2/2}
\nonumber\\
&=& \frac{1}{\sqrt{ 2\pi t}}
e^{-(w+ix)^2/2t},
\label{eqn:qBM2}
\end{eqnarray}
which is equal to $p(t,w|c^{-1}x)$ with $c=i$
as mentioned as (\ref{eqn:c_BM}) and (\ref{eqn:q_BM}) above.
By the Fourier reverse transform of (\ref{eqn:qBM1}),
we have
\begin{equation}
G_{\alpha}^{\rm BM}(t,B(t))
=\int_{\R} dw \,
q(t, w|B(t)) e^{i \alpha w}.
\label{eqn:qBM3}
\end{equation}
Expansion of the both sides with respect to $\alpha$ gives
\begin{equation}
\sum_{n=0}^{\infty} m_n(t, B(t)) 
\frac{\alpha^n}{n!}
= \sum_{n=0}^{\infty} \int_{\R} d w \,
q(t,w|B(t)) (iw)^n \frac{\alpha^n}{n!},
\label{eqn:qBM4}
\end{equation}
which implies
\begin{eqnarray}
\sfM[(iW)^n | (t, B(t))]
&\equiv& \int_{\R} dw \,
q(t,w|B(t)) (iw)^n
\nonumber\\
&=& m_n(t, B(t)), \qquad n \in \N_0, 
\label{eqn:qBM5}
\end{eqnarray}
where $m_n(t,x)$ is given by (\ref{eqn:Hermite_m1}).

\subsection{BESQ$^{(\nu)}$ 
\label{sec:int_BESQ}}

For BESQ$^{(\nu)}$, $V(t)=R^{(\nu)}(t), \nu > -1, t \in [0, \infty)$,
we set
\begin{equation}
c=-1.
\label{eqn:c_BESQ}
\end{equation}
and 
\begin{eqnarray}
&& q^{(\nu)}(t,y|x) = p^{(\nu)}(t, y|c^{-1}x)
\nonumber\\
&& \quad = \left\{ \begin{array}{ll}
\displaystyle{
\frac{1}{2t} \left( \frac{y}{x} \right)^{\nu/2}
\exp \left( - \frac{(-x)+y}{2t} \right)
J_{\nu} \left( \frac{\sqrt{xy}}{t} \right)},
& \quad t>0, x>0, y \in \R_+, \cr
\displaystyle{
\frac{y^{\nu}}{(2t)^{\nu+1} \Gamma(\nu+1)} e^{-y/2t}},
& \quad t >0, x=0, y \in \R_+, \cr
& \cr
\delta(y-x),
& \quad t=0, x, y \in \R_+,
\end{array} \right.
\label{eqn:q_BESQ}
\end{eqnarray}
where $J_{\nu}(x)$ is the Bessel function defined by
\begin{equation}
J_{\nu}(z) = \sum_{n=0}^{\infty} 
\frac{(-1)^n}{\Gamma(n+1) \Gamma(n+1+\nu)}
\left( \frac{z}{2} \right)^{2n+\nu}.
\label{eqn:J1}
\end{equation}
As usual we define $z^{\nu}$
to be $\exp(\nu \log z)$, where the argument of $z$
is given by its principal value;
$$
z^{\nu}=\exp \Big[ \nu \Big\{
\log |z| + \sqrt{-1} {\rm arg}(z) \Big\} \Big], \quad
-\pi < {\rm arg}(z) \leq \pi.
$$
In order to obtain (\ref{eqn:q_BESQ}) from (\ref{eqn:p_BESQ}),
we have used the relation
$$
I_{\nu}(z)=\left\{ \begin{array}{ll}
e^{-\nu \pi i/2} J_{\nu}(iz),
& \qquad - \pi < {\rm arg}(z) \leq \pi/2,
\cr
e^{3 \nu \pi i/2} J_{\nu}(iz),
& \qquad \pi/2 < {\rm arg}(z) \leq \pi.
\end{array} \right.
$$

Then we can prove the following.

\begin{lem}
\label{thm:int_BESQ}
With (\ref{eqn:c_BESQ}) and (\ref{eqn:q_BESQ}),
(\ref{eqn:int2}) are satisfied.
\end{lem}
{\it Proof.} \quad
Since $q^{(\nu)}(t, \cdot|x)=p^{(\nu)}(t, \cdot|c^{-1} x)$
satisfies
$$
\frac{\partial q^{(\nu)}}{\partial t}
=c^{-1} \left\{
2x \frac{\partial^2 q^{(\nu)}}{\partial x^2}
+2 (\nu+1) \frac{\partial q^{(\nu)}}{\partial x} \right\},
$$
$\sfM^{(\nu)}[f(W)|(t, R^{(\nu)}(t))]$ is a local martingale,
if $c^{-1}=-1 \Leftrightarrow c=-1$.
Therefore, the assignment (\ref{eqn:c_BESQ}) of the value $c$
guarantees that $\sfM[f(W)|(t,R^{(\nu)}(t))]$ is a local martingale.
The integral representations 
of the Laguerre polynomials in terms of Bessel functions
(for instance, see Eq.(6.2.15) in \cite{AAR99}), 
\begin{equation}
L_n^{(\nu)}(x)
=\frac{1}{n!} \frac{e^x}{x^{\nu/2}}
\int_0^{\infty} du \, e^{-u} u^{n+\nu/2}
J_{\nu} (2 \sqrt{xu}), \quad
n \in \N_0, \nu > -1, x \in \R_+,
\label{eqn:Lnint1}
\end{equation}
give (\ref{eqn:int2}). \qed
\vskip 0.5cm

\noindent{\bf Remark 5.} \,
The following integral formula is established,
\begin{equation}
\frac{e^x}{x^{\nu/2}}
\int_{0}^{\infty} du \, u^{\nu/2} e^{-(1-\alpha) u}
J_{\nu}(2 \sqrt{xu})
= \frac{e^{-x \alpha/(1-\alpha)}}{(1-\alpha)^{\nu+1}}.
\label{eqn:qBESQ1}
\end{equation}
It gives an integral representation
for $G_{\alpha}^{(\nu)}(t,x)$ 
studied in Remark 2,
\begin{equation}
G_{\alpha}^{(\nu)}(t,x)
=\int_0^{\infty} dw \,
q^{(\nu)}(t,w|x) e^{-\alpha w},
\label{eqn:qBESQ2}
\end{equation}
where $q^{(\nu)}$ is given by (\ref{eqn:q_BESQ}).

\subsection{CPR for BM 
\label{sec:cpr_BM}}

The integral formula (\ref{eqn:Hnint1}) implies
\begin{equation}
m_n(t,x)=\check{\rm E}_0 [(x+i W(t))^n],
\quad n \in \N_0, \quad
(t,x) \in [0, \infty) \times \R,
\label{eqn:mnZ1}
\end{equation}
where $\check{\rm E}_0$ denotes the expectation of
BM, $W(t), t \geq 0$. It is independent from
$B(t), t \geq 0$ and started at $W(0)=0$.
Then, if we consider the complex BM, 
\begin{equation}
Z(t)=B(t)+i W(t), \quad t \geq 0,
\label{eqn:cBM1}
\end{equation}
then
\begin{equation}
m_n(t, B(t))=
\sfM [(iW)^n |(t, B(t))]=
\check{\rm E}_0[ Z(t)^n],
\quad t \geq 0, \quad n \in \N.
\label{eqn:mnZ2}
\end{equation}
As a matter of course, the map
$z \to z^n, z \in \C, n \in \N$
are analytic, and then
$Z(t)^n, t \geq 0, n \in \N_0$ are conformal maps
of $Z(t), t \geq 0$. Since the probability distribution of
the complex BM is conformal invariant,
$Z(t)^n, t \geq 0, n \in \N$ are time changes of $Z(t)$.
In other words, $Z(\cdot)^n, n \in \N$ are
{\it conformal local martingales} 
(see Section V.2 of \cite{RY05}).
Since $B(\cdot)=\Re Z(\cdot)$ and 
$W(\cdot)=\Im Z(\cdot)$ are independent
one-dimensional standard BM's, 
$m_n(\cdot, B(\cdot)), n \in \N$, which are
obtained by taking the average over the imaginary parts
of $Z(\cdot)^n$ as (\ref{eqn:mnZ2}) are also
local martingales.

\subsection{CPR for Bessel processes (BES$^{(\nu)}$)
\label{sec:cpr_BES}}

The {\it Bessel process} with index $\nu$
(BES$^{(\nu)}$), $\widetilde{R}^{(\nu)}(t), t \geq 0$, 
is defined by
\begin{equation}
\widetilde{R}^{(\nu)}(t) \equiv 
\sqrt{R^{(\nu)}(t)},
\quad t \geq 0, \quad \nu > -1,
\label{eqn:Bessel}
\end{equation}
where $R^{(\nu)}(t), t \geq 0$ is BESQ$^{(\nu)}$.
It solves the SDE
\begin{equation}
d \widetilde{R}^{(\nu)}(t)
= d B(t) + \frac{2\nu+1}{2} \frac{dt}{\widetilde{R}^{(\nu)}(t)},
\quad t \geq 0.
\label{eqn:Bessel1b}
\end{equation}
The transition probability density 
is obtained from (\ref{eqn:p_BESQ}) as
\begin{eqnarray}
\widetilde{p}^{(\nu)}(t, y|x) &=& p^{(\nu)}(t, y^2|x^2) 2y
\nonumber\\
&=& \left\{ \begin{array}{ll}
\displaystyle{
\frac{1}{t} \frac{y^{\nu+1}}{x^{\nu}}
\exp \left( - \frac{x^2+y^2}{2t} \right)
I_{\nu} \left( \frac{xy}{t} \right)},
& \quad t>0, x>0, y \in \R_+, \cr
\displaystyle{
\frac{y^{2\nu+1}}{2^{\nu} t^{\nu+1} \Gamma(\nu+1)} e^{-y^2/2t}},
& \quad t >0, x=0, y \in \R_+, \cr
& \cr
\delta(y-x),
& \quad t=0, x, y \in \R_+.
\end{array} \right.
\label{eqn:p4}
\end{eqnarray}

Corresponding to this,
the integral transform (\ref{eqn:int1}) 
for BESQ$^{(\nu)}$, which is denoted as 
$\sfM^{(\nu)}[\cdot|\cdot]$,
is converted into that
for BES$^{(\nu)}$ expressed as 
$\widetilde{\sfM}^{(\nu)}[\cdot|\cdot]$
so that the following relation
holds,
\begin{equation}
\sfM^{(\nu)}[f(-W)|(t, R^{(\nu)}(t))]
=\widetilde{\sfM}^{(\nu)}
[\widetilde{f}(iW)|(t, \widetilde{R}^{(\nu)}(t))],
\quad t \geq 0,
\label{eqn:BESBESQ1}
\end{equation}
where $f$ and $\widetilde{f}$ are polynomials
with the relation
$\widetilde{f}(z)=f(z^2), z \in \C$.
For it,  we set
\begin{eqnarray}
\widetilde{q}^{(\nu)}(t, y|x)
 &=& q^{(\nu)}(t, y^2|x^2) 2y
\nonumber\\
&=& \left\{ \begin{array}{ll}
\displaystyle{
\frac{1}{t} \frac{y^{\nu+1}}{x^{\nu}}
\exp \left( - \frac{(-x^2)+y^2}{2t} \right)
J_{\nu} \left( \frac{xy}{t} \right)},
& \quad t>0, x>0, y \in \R_+, \cr
\displaystyle{
\frac{y^{2\nu+1}}{2^{\nu} t^{\nu+1} \Gamma(\nu+1)} e^{-y^2/2t}},
& \quad t >0, x=0, y \in \R_+, \cr
& \cr
\delta(y-x),
& \quad t=0, x, y \in \R_+.
\end{array} \right.
\label{eqn:q4}
\end{eqnarray}
and define
the integral transform for BES$^{(\nu)}, \nu > -1$ by
\begin{equation}
\widetilde{\sfM}^{(\nu)}[f(W)|(t,x)]
=\int_{\R_+} dw \, \widetilde{q}^{(\nu)}(t,w|x) f(w),
\quad (t, x) \in [0, \infty) \times \R_+
\label{eqn:transBES1}
\end{equation}
for an integrable function $f$.
Then the relation (\ref{eqn:BESBESQ1}) is satisfied.

For $m \in \N_0$, the Bessel functions have
the following expansions by the trigonometric functions,
\begin{eqnarray}
J_{2m+1/2}(x)
&=& (-1)^m \sqrt{\frac{2}{\pi x}}
\left[ \sin x \sum_{k=0}^m
\frac{(-1)^k (2m+2k)!}{(2k)!(2m-2k)!} (2x)^{-2k} 
\right.
\nonumber\\
&& \qquad \qquad \quad \left.
+ \cos x \sum_{k=0}^{m-1}
\frac{(-1)^k (2m+2k+1)!}{(2k+1)!(2m-2k-1)!} (2x)^{-(2k+1)} \right],
\nonumber\\
J_{2m+3/2}(x)
&=& (-1)^m \sqrt{\frac{2}{\pi x}}
\left[ - \cos x \sum_{k=0}^m
\frac{(-1)^k (2m+2k+1)!}{(2k)!(2m-2k+1)!} (2x)^{-2k} 
\right.
\nonumber\\
&& \qquad \qquad \quad \left.
+ \sin x \sum_{k=0}^{m}
\frac{(-1)^k (2m+2k+2)!}{(2k+1)!(2m-2k)!} (2x)^{-(2k+1)} \right].
\label{eqn:Jodd}
\end{eqnarray}
They are obtained from Eq.(4.6.12) in \cite{AAR99}.
For example, if we set $m=0$ in (\ref{eqn:Jodd}), we have
$$
J_{1/2}(x)=\sqrt{\frac{2}{\pi x}} \sin x, \quad
J_{3/2}(x)=\sqrt{\frac{2}{\pi x}} 
\left( \frac{\sin x}{x} - \cos x \right).
$$
Assume that $\widetilde{f}(z)$ is a polynomial of $z^2$,
and thus 
\begin{equation}
\widetilde{f}(-z)=\widetilde{f}(z), \quad z \in \C.
\label{eqn:symf}
\end{equation}
Then (\ref{eqn:transBES1}) with (\ref{eqn:q4}) gives 
\begin{eqnarray}
\widetilde{\sfM}^{(1/2)} \left[ \left. \widetilde{f}(iW) 
\right| (t,x) \right]
&=& \int_{\R_+} dw \, \sqrt{\frac{2}{\pi t}} \frac{w}{x}
e^{-(-x^2+w^2)/2t} \sin (xw/t) \widetilde{f}(iw)
\nonumber\\
&=& \frac{1}{\sqrt{2 \pi t}} \frac{1}{i x} 
\int_{\R_+} dw \, w
\left\{ e^{-(w-ix)^2/2t}-e^{-(w+ix)^2/2t} \right\}
\widetilde{f}(iw),
\nonumber
\end{eqnarray}
and
\begin{eqnarray}
&& 
\widetilde{\sfM}^{(3/2)} \left[ \left. \widetilde{f}(iW) 
\right| (t,x) \right]
\nonumber\\
&& \qquad 
= \int_{\R_+} dw \, \sqrt{\frac{2}{\pi t}} \frac{w}{x}
e^{-(-x^2+w^2)/2t} 
\left\{ \frac{t}{xw} \sin (xw/t) - \cos (xw/t) \right\}
\widetilde{f}(iw)
\nonumber\\
&& \qquad
= \frac{1}{\sqrt{2 \pi t}}
\left[ \frac{t}{i x^3} \int_{\R_+} dw \, w
\left\{ e^{-(w-ix)^2/2t}-e^{-(w+ix)^2/2t} \right\}
\right.
\nonumber\\
&& \qquad \qquad \qquad \left.
- \frac{1}{x^2} 
\int_{\R_+} dw \, w^2
\left\{ e^{-(w-ix)^2/2t}+e^{-(w+ix)^2/2t} \right\}
\right] \widetilde{f}(iw).
\nonumber
\end{eqnarray}
By the assumption (\ref{eqn:symf}), they are rewritten as
\begin{eqnarray}
\widetilde{\sfM}^{(1/2)} \left[ \left. \widetilde{f}(iW) 
\right| (t,x) \right]
&=&
\frac{1}{\sqrt{2 \pi t}}
\frac{1}{(-i)x} \int_{\R} dw \, w
e^{-(w+ix)^2/2t} \widetilde{f}(iw)
\nonumber\\
&=& \frac{1}{\sqrt{2 \pi t}}
\int_{-\infty+ix}^{\infty+ix} du \,
\frac{x+iu}{x} e^{-u^2/2t} \widetilde{f}(x+iu),
\nonumber\\
\widetilde{\sfM}^{(3/2)} \left[ \left. \widetilde{f}(iW) 
\right| (t,x) \right]
&=&
\frac{1}{\sqrt{2 \pi t}}
\left[ \frac{t}{(-i) x^3}
\int_{\R} dw \, w e^{-(w+ix)^2/2t} \widetilde{f}(i w)
\right.
\nonumber\\
&& \qquad \qquad  \left.
+ \frac{1}{(-i)^2 x^2}
\int_{\R} dw \, w^2 e^{-(w+ix)^2/2t} \widetilde{f}(iw) \right]
\nonumber\\
&=& \frac{1}{\sqrt{2 \pi t}}
\int_{-\infty+ix}^{\infty+ix} du \,
\left\{ \frac{t(x+iu)}{x^3}+\frac{(x+iu)^2}{x^2} \right\}
e^{-u^2/2t} \widetilde{f}(x+iu),
\nonumber
\end{eqnarray}
where we have changed the integral variables by
$u=w+ix$. Since the integrands are entire,
$\int_{-\infty +ix}^{\infty +ix} du \, (\cdot)$ can be
replaced by $\int_{\R} du \, (\cdot)$.
Then we have the following expressions for the martingales
(\ref{eqn:BESBESQ1}) with $\nu=1/2$ and 3/2 
\begin{eqnarray}
&& \widetilde{\sfM}^{(1/2)} \left[ \left. \widetilde{f}(iW) \right|
(t, \widetilde{R}^{(1/2)}(t)) \right]
=
\check{\rE}_0 \left[ 
\frac{Z^{(1/2)}(t)}{\Re Z^{(1/2)}(t)} 
\widetilde{f}(Z^{(1/2)}(t))
\right], \quad t \geq 0,
\nonumber\\
&&
\widetilde{\sfM}^{(3/2)} \left[ \left. \widetilde{f}(iW) \right|
(t, \widetilde{R}^{(1/2)}(t)) \right]
\nonumber\\
&& \qquad =
\check{\rE}_0 \left[
\left\{ \frac{t Z^{(3/2)}(t)}{(\Re Z^{(3/2)}(t))^3}
+\frac{(Z^{(3/2)}(t))^2}{(\Re Z^{(3/2)}(t))^2} \right\}
\widetilde{f}(Z^{(3/2)}(t)) \right], \quad t \geq 0,
\label{eqn:BES_mart0}
\end{eqnarray}
where
$Z^{(\nu)}(t)=\widetilde{R}^{(\nu)}(t)+iW(t), \nu=1/2, 3/2$.
These calculations are generalized as follows.

\begin{lem}
\label{thm:BES_mart1}
Let $\widetilde{f}(z)$ be a polynomial of $z^2$.
Then 
$\widetilde{\sfM}^{(n+1/2)}[\widetilde{f}(iW)|(\cdot,
\widetilde{R}^{(n+1/2)}(\cdot))]$,
$n \in \N_0$, are local martingales,
and
\begin{equation}
\widetilde{\sfM}^{(n+1/2)} \left[ \left. \widetilde{f}(iW) \right|
(t, \widetilde{R}^{(n+1/2)}(t)) \right]
=\check{\rE}_0 \left[
Q_t^{(n+1/2)}(Z^{(n+1/2)}(t)) 
\widetilde{f}(Z^{(n+1/2)}(t)) \right],
\quad t \geq 0, 
\label{eqn:BES_mart1}
\end{equation}
where
\begin{equation}
Z^{(n+1/2)}(t)=\widetilde{R}^{(n+1/2)}(t)
+ i W(t),
\quad t \geq 0,
\label{eqn:ZZ1}
\end{equation}
and 
\begin{equation}
Q^{(n+1/2)}_t(z)
=\left( \frac{t}{2} \right)^n
\frac{z}{(\Re z)^{2n+1}}
\sum_{k=0}^n \frac{(2n-k)!}{(n-k)! k!}
\left( \frac{2 (\Re z) z}{t} \right)^k, \quad z \in \C.
\label{eqn:Q1}
\end{equation}
\end{lem}

Note that the equalities (\ref{eqn:mnZ2}) with (\ref{eqn:Hermite_m2})
hold even if we replace the complex BM, (\ref{eqn:cBM1}), 
by the present complex diffusion, (\ref{eqn:ZZ1}),
since the imaginary parts are the same; 
\begin{equation}
\check{\rE}_0 
[ (Z^{(n+1/2)}(t))^{k} ]
=m_{k}(t, \widetilde{R}^{(n+1/2)}(t)),
\quad t \geq 0, \quad n, k \in \N_0.
\label{eqn:ZZ2}
\end{equation}
Then for monomials
$\widetilde{f}(z)=z^{2 \ell}$, (\ref{eqn:BES_mart1}) 
gives the following.
For $n, \ell \in \N_0$, 
\begin{eqnarray}
&& \widetilde{\sfM}^{(n+1/2)} 
\Big[ (i W)^{2 \ell} \Big| (t, \widetilde{R}^{(n+1/2)}(t)) \Big]
\nonumber\\
&& \qquad = \left( \frac{t}{2} \right)^n 
\frac{1}{ (\widetilde{R}^{(n+1/2)}(t))^{2n+1}}
\sum_{k=0}^{n} \frac{(2n-k)!}{(n-k)! k!}
\left( \frac{2 \widetilde{R}^{(n+1/2)}(t)}{t} \right)^k 
\check{\rE}_0 [ (Z^{n+1/2}(t))^{2 \ell+k+1} ]
\nonumber\\
&& \qquad = \frac{1}{2^{2n+1}} \left( \frac{t}{2} \right)^{\ell}
\left( \frac{\widetilde{R}^{(n+1/2)}(t)}{\sqrt{2t}} \right)^{-(2n+1)}
\nonumber\\
&& \qquad \qquad \qquad \times
\sum_{k=0}^{n} \frac{(2n-k)!}{(n-k)! k!}
\left( 2 \frac{\widetilde{R}^{(n+1/2)}(t)}{\sqrt{2t}} \right)^{k}
H_{2 \ell+k+1} \left( \frac{\widetilde{R}^{(n+1/2)}(t)}{\sqrt{2t}} \right).
\label{eqn:tildeM}
\end{eqnarray}

\subsection{CPR for RW
\label{sec:cpr_RW}}

For $t >0$, let $\eta_{\ell}(t), \ell \in \N$ be a series of i.i.d.
random variables with the Gamma($t$) probability density
\begin{equation}
\rP^{\Gamma}( \cdot \in dx )=
\frac{1}{\Gamma(t)} x^{t-1} e^{-t} dx,
\quad x > 0,
\label{eqn:Gamma1}
\end{equation}
where $\Gamma(t)$ is the Gamma function.
When $t \in \N$,
\begin{equation}
\eta_{\ell}(t) \d= \varepsilon_{\ell}^{(1)}+ \cdots
+ \varepsilon_{\ell}^{(t)},
\quad \ell \in \N,
\label{eqn:Gamma2}
\end{equation}
where $\varepsilon_{\ell}^{(k)}, k=1,2, \dots, t$ are independent
random variables with standard exponential distribution
${\rm Prob}(\varepsilon_{\ell}^{(k)} \geq x)=e^{-x}, x \geq 0$.
We consider a random variable
\begin{equation}
C(t)=\frac{2}{\pi^2} \sum_{\ell \in \N}
\frac{\eta_{\ell}(t)}{(\ell-1/2)^2},
\label{eqn:Ct1}
\end{equation}
since it is known \cite{BPY01} that 
its Laplace transform is given by
\begin{equation}
\rE^{\Gamma} [ e^{-\lambda C(t)} ]
=\frac{1}{(\cosh \sqrt{2 \lambda})^t},
\quad t >0, 
\label{eqn:Ct2}
\end{equation}
where $\rE^{\Gamma}$ denotes the expectation
with respect to $\eta_{\ell}(t), \ell \in \N$.
In \cite{BPY01}, it is shown that $C(t) \in [0, \infty)$
is infinitely divisible and its probability density
$\mu_{C(t)}(\cdot)$, which is defined for integrable functions $f$ as
\begin{equation}
\rE^{\Gamma}[f(C(t))]
=\int_{0}^{\infty} d c \, \mu_{C(t)}(c) f(c), 
\label{eqn:mu1}
\end{equation}
is explicitly given by
\begin{equation}
\mu_{C(t)}(c)
=\frac{2^t}{\Gamma(t)}
\sum_{\ell=0}^{\infty} (-1)^{\ell}
\frac{\Gamma(\ell+t)}{\Gamma(\ell+1)}
\frac{(2\ell+t)}{\sqrt{2 \pi c^3}}
e^{-(2\ell+t)^2/2c},
\quad t > 0, x > 0.
\label{eqn:mu2}
\end{equation}
The generating function of the polynomials
(\ref{eqn:GG1}) is thus written as
\begin{eqnarray}
G_{\alpha}(t,x)
&=& \rE^{\Gamma}  \Big[e^{\alpha x-\alpha^2 C(t)/2}  \Big]
\nonumber\\
&=& \sum_{n=0}^{\infty}
\frac{\alpha^n}{n!} 
\rE^{\Gamma} \left[
\left( \frac{C(t)}{2} \right)^{n/2}
H_n \left( \frac{x}{\sqrt{2 C(t)}} \right) \right],
\label{eqn:Ct3}
\end{eqnarray}
where we used the formulas (\ref{eqn:GBM1}) and (\ref{eqn:GBM2}).
It gives the relations
\begin{equation}
m_n(t,x)=\rE^{\Gamma} \left[
m_n^{\rm BM}(C(t), x) \right],
\quad n \in \N_0, \quad t \in N_0,
\label{eqn:Ct4}
\end{equation}
that is, let $m_n^{\rm BM}(C(t), x)$ be a random time change
$t \to C(t)$ of the polynomial for BM, 
then its average over $C(t)$ gives Fujita's polynomial for RW.

On the other hand, by (\ref{eqn:mnZ1}), 
if we introduce a one-dimensional BM,
$W(t), t \geq 0$ with $W(0)=0$, and write the expectation
with respect to $W(\cdot)$ as $\check{\rE}$,
we have the following expressions
\begin{equation}
m_n^{\rm BM}(t,x)
=\check{\rE} [(x+i W(t))^n],
\quad n \in \N_0, \quad t \geq 0.
\label{eqn:cp1}
\end{equation}
Combination of (\ref{eqn:Ct4}) and (\ref{eqn:cp1}) gives
\begin{equation}
m_n(t,x)=\rE^{\Gamma} \Big[
\check{\rE} \Big[ (x+i W (C(t)))^n \Big] \Big],
\quad n \in \N_0, \quad t \in \N_0.
\label{eqn:CPR1RW}
\end{equation}

Let
\begin{eqnarray}
\widetilde{W}(t) &\d=& W(C(t))
\d= \sqrt{\frac{C(t)}{t}} W(t),
\quad t \in \N,
\nonumber\\
\widetilde{W}(0) &=& 0.
\label{eqn:Wtilde1}
\end{eqnarray}
We regard $\widetilde{W}(t), t \in \N_0$ as a discrete-time
process and the expectation w.r.t. this process is written as
\begin{equation}
\widetilde{\rE}[f(\widetilde{W}(t))]
=\rE^{\Gamma} \Big[
\check{\rE} [ f(W(C(t))) ] \Big], 
\quad t \in \N_0
\label{eqn:Wtilde2}
\end{equation}
for integrable functions $f$.

With RW, $V(t), t \in \N_0$, we consider a 
discrete-time complex process
\begin{equation}
Z(t)=V(t)+i \widetilde{W}(t), \quad t \in \N_0.
\label{eqn:cp2}
\end{equation}
Note that $\Re Z(t)=V(t) \in \Z$ and
$\Im Z(t)=\widetilde{W}(t) \in \R$.
The above results are summarized as follows \cite{Kat13b}.

\begin{lem}
\label{thm:cpr}
Fujita's polynomial martingales, 
$m_n(t,V(t)), n \in \N_0, t \in \N_0$, for the
simple and symmetric RW have
the following complex-process representations,
\begin{equation}
m_n(t, V(t))
= \widetilde{\rE} [ Z(t)^n ],
\quad n \in \N_0, \quad t \in \N_0.
\label{eqn:cp3}
\end{equation}
\end{lem}

\clearpage
\SSC{Noncolliding Diffusion Processes 
\label{sec:noncolliding}}
\subsection{Map for martingales \label{sec:mapB}}

For a configuration
\begin{equation}
\xi(\cdot)=\sum_{j=1}^N \delta_{u_j}(\cdot) \in \mM_0,
\label{eqn:xiB1}
\end{equation}
we define a polynomial of $x \in \C$ with
a parameter $u \in \C$ as
\begin{equation}
\Phi_{\xi}^{u}(x)=\prod_{r \in \supp \xi \cap \{u\}^{\rm c}}
\frac{x-r}{u-r}.
\label{eqn:Phi1}
\end{equation}
with $\supp \xi =\{u_j : 1 \leq j \leq N\}$.
Note that
\begin{equation}
\Phi_{\xi}^{u_k}(u_j)=\delta_{jk}, \quad
1 \leq j, k \leq N.
\label{eqn:Phi2}
\end{equation}
Then we have the following statement.

\begin{prop}
\label{thm:mapB}
For the Markov process $V(t), t \in \cT$,
assume that the integral transform (\ref{eqn:int1})
satisfying (\ref{eqn:int2}) is obtained.
Then with $\u \in \W_N$ and $\xi=\sum_{j=1}^N \delta_{u_j} \in \mM_0$,
the map (\ref{eqn:map1}) satisfying (\ref{eqn:map2}) 
is given by
\begin{equation}
\cM_{\xi}^u(\cdot, \cdot)
=\sfM[ \Phi_{\xi}^u(cW) |(\cdot, \cdot)].
\label{eqn:mapB1}
\end{equation}
\end{prop}
{\it Proof.} \quad
By definition (\ref{eqn:Phi1}),
$\Phi_{\xi}^u(x)$ is polynomial.
Then $\cM_{\xi}^u(\cdot, V(\cdot))$ is polynomial 
and local martingale by (\ref{eqn:int3}).
By (\ref{eqn:int4}),
$$
\cM_{\xi}^u(0,V(0))
=\sfM[\Phi_{\xi}^u(cW)|(0, V(0))]
=\Phi_{\xi}^u(V(0)).
$$
Then
$$
\cM_{\xi}^{u_k}(0,u_j)=\Phi_{\xi}^{u_k}(u_j)=\delta_{jk},
\quad 1 \leq j, k \leq N
$$
by (\ref{eqn:Phi2}).
The proof is completed. \qed
\vskip 0.3cm

The determinantal martingale (\ref{eqn:D1})
is now given as
\begin{equation}
\cD_{\xi}(t, \V_{\J}(t))
= \det_{j, k \in \J}
\Big[ \sfM[ \Phi_{\xi}^{u_k}(cW)|
(t, V_j(t))] \Big],
\quad \J \subset \I_N, \quad t \in \cT.
\label{eqn:D3}
\end{equation}

The integral transform (\ref{eqn:int1}) is extended
to the linear integral transform of functions
of $\x \in S^N$ such that, if
$F^{(k)}(\x)=\prod_{j=1}^N f_j^{(k)}(x_j)$
with integrable functions $f_j^{(k)}, 1 \leq j \leq N, k=1,2$,
then
\begin{equation}
\sfM \left[ F^{(k)}(\bW) \left| \{(t_{\ell}, x_{\ell})\}_{\ell=1}^N  \right. \right]
= \prod_{j=1}^N \sfM \left[ \left. f^{(k)}_j(W_j) \right| (t_{j}, x_{j}) \right],
\quad k=1,2
\label{eqn:intB1}.
\end{equation}
and
\begin{eqnarray}
&& \sfM \Big[c_1 F^{(1)}(\bW) +c_2 F^{(2)}(\bW) \left| \{(t_{\ell}, x_{\ell})\}_{\ell=1}^N 
\right. \Big]
\nonumber\\
&& \quad
= c_1 \sfM \Big[F^{(1)}(\bW) \left| \{(t_{\ell}, 
x_{\ell})\}_{\ell=1}^N \right. \Big]
+ c_2 \sfM \Big[F^{(2)}(\bW) \left| \{(t_{\ell}, 
x_{\ell})\}_{\ell=1}^N \right. \Big],
\label{eqn:intB2}
\end{eqnarray}
$c_1, c_2 \in \C$, 
for $0 < t_j < \infty, 1 \leq j \leq N$,
where $\bW=(W_1, \dots, W_N) \in S^N$.
In particular, if $t_{\ell}=t, 1 \leq {^{\forall}\ell} \leq N$, we write
$\sfM[\cdot | \{(t_{\ell}, x_{\ell})\}_{\ell=1}^N]$ simply as
$\sfM[\cdot|(t, \x)]$ with $\x=(x_1, \dots, x_N)$.
Then, by the multilinearity of determinant,
(\ref{eqn:D3}) is written as
\begin{equation}
\cD_{\xi}(t, \V_{\J}(t))
= \sfM \left[ \left. \det_{j, k \in \J}
\left[ \Phi_{\xi}^{u_k}(cW_j) \right] \right|
(t, \V_{\J}(t))] \right],
\quad \J \subset \I_N, \quad t \in \cT.
\label{eqn:D4}
\end{equation}

\subsection{Krattenthaler's determinant identity 
and $h$-transform
\label{sec:Kra}}

The following determinant identity was given as
Lemma 2.2 in \cite{Kra90} and as Lemma 3 in \cite{Kra99}
proved by Krattenthaler.

\begin{lem}
\label{thm:Kra1}
Let $X_1, \dots, X_N, A_2, \dots, A_N$ and
$B_2, \dots, B_N$ be indeterminates.
Then there holds
\begin{eqnarray}
&&\det_{1 \leq j, k \leq N}
\Big[ (X_j+A_N)(X_j+A_{N-1}) \cdots (X_j+A_{k+1})
(X_j+B_{k})(X_j+B_{k-1}) \cdots (X_j+B_2) \Big]
\nonumber\\
&& \qquad \qquad \qquad \qquad =
\prod_{1 \leq j < k \leq N} (X_j-X_k)
\prod_{2 \leq j \leq k \leq N} (B_j-A_k).
\label{eqn:Kra1}
\end{eqnarray}
\end{lem}
\vskip 0.3cm

The Vandermonde determinant is given as
\begin{equation}
h(\x)=\det_{1 \leq j, k \leq N}
[x_j^{k-1}] = \prod_{1 \leq j< k \leq N} (x_k-x_j),
\label{eqn:Vand}
\end{equation}
for $\x=(x_1, \dots, x_N) \in S^N$.
As a special case of (\ref{eqn:Kra1}),
we obtain the following determinant identity.

\begin{lem}
\label{thm:det_identity}
Assume that $N \in \N, \x \in \C^N, \u \in \W_N$.
Then
\begin{equation}
\frac{ h(\x)}{h(\u)}
=\det_{1 \leq j, k \leq N} [ \Phi_{\xi}^{u_k}(x_j) ].
\label{eqn:det_identity}
\end{equation}
\end{lem}
{\it Proof.} \quad
In the identity (\ref{eqn:Kra1}), set
\begin{eqnarray}
&& X_j=x_j, \quad 1 \leq j \leq N,
\nonumber\\
&& A_j=-u_j,\quad B_j=-u_{j-1}, \quad 2 \leq j \leq N.
\nonumber
\end{eqnarray}
Then we find
$$
H(\u, \x) \equiv \det_{1 \leq j, k \leq N}
\left[ \prod_{1 \leq \ell \leq N, \ell \not=j}
(u_{\ell}-x_k) \right]
=(-1)^{N(N-1)/2} h(\u) h(\x).
$$
Since 
\begin{eqnarray}
\det_{1 \leq j, k \leq N}
[\Phi_{\xi}^{u_{k}}(x_j)]
&=& \frac{H(\u, \x)}{
\prod_{1 \leq j \leq N} \prod_{1 \leq k \leq N: k \not= j} (u_{k}-u_j)}
\nonumber\\
&=& \frac{H(\u, \x)}{(-1)^{N(N-1)/2} h(\u)^2},
\nonumber
\end{eqnarray}
the identity (\ref{eqn:det_identity}) is obtained
as a special case of (\ref{eqn:Kra1}).
\qed
\vskip 0.3cm

Using the determinant identity (\ref{eqn:det_identity})
with (\ref{eqn:Vand}), we see
\begin{eqnarray}
\cD_{\xi}(t, \V(t)) 
&=& \sfM \left[ \left. \det_{1 \leq j, k \leq N}
[ \Phi_{\xi}^{u_k}(c W_j)] \right| (t, \V(t))] \right]
\nonumber\\
&=& \sfM \left[ \left. 
\frac{h(c\bW)}{h(\u)} \right| (t, \V(t)) \right]
\nonumber\\
&=& \sfM \left[ \left.
\frac{1}{h(\u)} \det_{1 \leq j, k \leq N}
[(cW_j)^{k-1}] \right| (t, \V(t)) \right]
\nonumber\\
&=& \frac{1}{h(\u)}
\det_{1 \leq j, k \leq N}
\Big[ \sfM [(cW_j)^{k-1} | (t, V_j(t))] \Big]
\nonumber\\
&=& \frac{1}{h(\u)} \det_{1 \leq j, k \leq N}
[m_{k-1}(t, V_j(t))].
\label{eqn:D5}
\end{eqnarray}
By multilinearity of determinant, the Vandermonde determinant 
$\det_{1 \leq j, k \leq N} [x_j^{k-1} ] $ does not change
by replacing $x_j^{k-1}$ by any monic polynomial of $x_j$ of 
degree $k-1$, $1 \leq j, k \leq N$.
Since $m_{k-1}(t,x_j)$ is a monic polynomial of $x_j$
of degree $k-1$, (\ref{eqn:D5}) is equal to
\begin{eqnarray}
\cD_{\xi}(t, \V(t))
&=& \frac{1}{h(\u)} \det_{1 \leq j, k \leq N}[(V_j(t))^{k-1}]
\nonumber\\
&=& \frac{h(\V(t))}{h(\u)}.
\label{eqn:D6}
\end{eqnarray}
This is the factor used for the harmonic transform 
($h$-transform).

\subsection{Noncolliding BM and noncolliding BESQ$^{(\nu)}$ 
\label{sec:NC}}

For $N \in \N$, 
we consider $N$-particle systems of BM's,
$\X(t)=(X_1(t), X_2(t), \dots, X_N(t))$, $t \geq 0$, 
and of BESQ$^{(\nu)}$ with
index $\nu > -1$,
$\X^{(\nu)}(t)=(X^{(\nu)}_1(t), X^{(\nu)}_2(t), \dots,
X^{(\nu)}_N(t)), t \geq 0$, 
both {\it conditioned never to collide with each other particle}.
The former process, which is called the {\it noncolliding BM},
solves the following set of SDEs
\begin{equation}
dX_j(t)=dB_j(t)+ 
\sum_{\substack{ 1 \leq k \leq N,\\ k \not= j}}
\frac{dt}{X_j(t)-X_k(t)},
\quad
1 \leq j \leq N, \quad t \geq 0,
\label{eqn:noncollBM}
\end{equation}
with independent one-dimensional standard
BMs, $B_j(t), 1 \leq j \leq N, t \geq 0$
\cite{Dys62,Spo87,Gra99,Joh01,KT10,Osa12,Osa13,Osa13b}.
The latter process,
the {\it noncolliding BESQ$^{(\nu)}$}, does the following set of SDEs
\begin{eqnarray}
dX^{(\nu)}_j(t) &=& 2 \sqrt{X^{(\nu)}_j(t)} 
d \check{B}_j(t) + 2 (\nu+1) dt
\nonumber\\
&& + 4 X^{(\nu)}_j(t)
\sum_{\substack{1 \leq k \leq N, \\ k \not= j}}
\frac{dt}{X^{(\nu)}_j(t)-X^{(\nu)}_k(t)},
\quad 1 \leq j \leq N, \quad t \geq 0,
\label{eqn:noncollBESQ}
\end{eqnarray}
where $\check{B}_j(t), 1 \leq j \leq N, t \geq 0$
are independent one-dimensional standard BMs
different from $B_j(t), 1 \leq j \leq N, t \geq 0$, 
and, if $-1 < \nu < 0$,
the reflection boundary condition is assumed at
the origin \cite{KO01,KT11}. 

Consider subsets of $\R^{N}$,
$\W_N^{\rm A}=\{\x=(x_1, x_2, \dots, x_N) \in \R^N :
x_1 < \cdots < x_N\}$,
and 
$\W_N^{+}=\{\x \in \R_+^N: 
x_1 < \cdots < x_N\}$.
The former is called the Weyl chambers of
types A$_{N-1}$.
If we replace the condition $\x \in \R_+^N$
by $\x \in (0, \infty)^N$ for the
latter, it will be the Weyl chamber of 
type C$_{N}$.
It is proved that, provided $\X(0) \in \W_N^{\rm A}$
and $\X^{(\nu)}(0) \in \W_N^{+}$,
then the SDEs (\ref{eqn:noncollBM})
and (\ref{eqn:noncollBESQ}) guarantee
that with probability one
$\X(t) \in \W_N^{\rm A}$, and
$\X^{(\nu)}(t) \in \W_N^{+}, \forall t > 0$ \cite{GJ13}.
That is, 
in both processes, at any positive time $t > 0$ 
there is no multiple point
at which coincidence of particle positions $X_j(t)=X_{k}(t)$ or
$X^{(\nu)}_j(t) = X^{(\nu)}_k(t)$ 
for $j \not= k$ occurs.
It is the reason why these processes are called
{\it noncolliding diffusion processes} \cite{KT11b}.
In general, however, 
we can consider them starting from initial configurations
with multiple points.
In order to describe a general initial configuration
we express it by a sum of delta measures
in the form 
$\xi(\cdot)=\sum_{j=1}^N \delta_{x_j}(\cdot)$.

Let $\mM$ be the space of nonnegative integer-valued Radon measures 
on $\R$. 
For an element $\xi$ of $\mM$, 
$\xi(\cdot) = \sum_{j \in \I}\delta_{x_j}(\cdot)$
with a countable index set $\I$, 
we introduce the following operations.

\begin{description}
\item[(shift)] with $u \in \R$, 
$\tau_u \xi(\cdot) =\displaystyle{\sum_{i \in \I}} 
\delta_{x_i+u}(\cdot)$,

\item[(dilatation)] with $c>0$,
$c \circ \xi(\cdot)=\displaystyle{\sum_{i \in \I} 
\delta_{c x_i}(\cdot)}$,

\item[(square)]
$\displaystyle{
\xi^{\langle 2 \rangle}(\cdot)
=\sum_{i \in \I} \delta_{x_i^2} (\cdot)}$.
\end{description}

Let $\mM^+=\{(\xi \cap \R_+): \xi \in \mM \}$.
We consider the noncolliding BM and the 
noncolliding BESQ as $\mM$-valued and
$\mM^+$-valued processes and write them as
\begin{equation}
\Xi(t, \cdot)=\sum_{j=1}^N \delta_{X_j(t)}(\cdot),
\quad
\Xi^{(\nu)}(t, \cdot)=\sum_{j=1}^N
\delta_{X^{(\nu)}_j(t)}(\cdot),
\quad t \geq 0,
\label{eqn:Xi1b}
\end{equation}
respectively \cite{KT10,KT11}. 
The probability law of $\Xi(t, \cdot)$
starting from a fixed configuration $\xi \in \mM$
is denoted by $\P_{\xi}$ and that of $\Xi^{(\nu)}(t, \cdot)$
from $\xi \in \mM^+$ by $\P^{(\nu)}_{\xi}$,
and the noncolliding diffusion processes
specified by initial configurations 
are expressed by
$(\Xi(t), \P_{\xi})$ and
$(\Xi^{(\nu)}(t), \P_{\xi}^{(\nu)}), \nu > -1$.
The expectations w.r.t.$\P_{\xi}$ and $\P_{\xi}^{(\nu)}$
are denoted by $\E_{\xi}$ and $\E_{\xi}^{(\nu)}$,
respectively.
The set of $\hmM$-valued
continuous functions defined on $[0,\infty)$ 
is denoted by $\rC([0,\infty) \to \hmM)$
for $\hmM=\mM$ or $\mM^+$.
We introduce a filtration $\{ {\cal F}(t) \}_{t\in [0,\infty)}$
on the space $\rC([0,\infty)\to \hmM)$
defined by ${\cal F}(t) = \sigma (\hXi(s), s \in [0,t])$,
where $\hXi(\cdot)=\Xi(\cdot)$ for $\hmM=\mM$
and $\hXi(\cdot)=\Xi^{(\nu)}(\cdot)$ for $\hmM=\mM^+$. 
Let $\rC_0(S)$ be the set of all continuous
real-valued functions with compact supports on 
$S=\R$ or $\R_+$.
We set 
$
\hmM_{0}= \{ \xi \in \hmM : 
\xi(\{x\})\le 1 \mbox { for any }  x \in S \}$,
which denotes collections of configurations
without any multiple points.

\subsection{$H$-transforms of absorbing processes
\label{sec:harmonic}}

For the noncolliding BM,
$\hXi(\cdot)=\Xi(\cdot)$ 
(resp. BESQ$^{(\nu)}, \nu > -1$, $\hXi(\cdot)=\Xi^{(\nu)}(\cdot)$), we shall set
$\hX(\cdot)=\X(\cdot)$ (resp. $\X^{(\nu)}(\cdot)$),
$\hP_{\xi}=\P_{\xi}$ (resp. $\P^{(\nu)}_{\xi}$),
$\hE_{\xi}=\E_{\xi}$ (resp. $\E^{(\nu)}_{\xi}$),
and 
$S=\R$ (resp. $\R_+$).

Let $0 < t \leq T < \infty$.
For the noncolliding diffusion process $(\hXi(t), \P_{\xi})$,
we consider the expectation
for an $\cF(t)$-measurable bounded function $F$,
$$
\E_{\xi} [F(\hXi(\cdot))].
$$
It is sufficient 
to consider the case that $F$ is given as
$F\left(\hXi(\cdot)\right)
= \prod_{m=1}^M g_m(\hX(t_m))$
for an arbitrary $M \in \N$, $0<t_1< \cdots <t_M \leq T < \infty$, 
with symmetric bounded measurable functions $g_m$ 
on $S^N$, $1 \leq m \leq M$.

We can prove that the noncolliding BM is obtained as 
an $h$-transform of the absorbing BM,
$\bB(t)=(B_1(t), \dots, B_N(t)), t \geq 0$ in the
Weyl chamber $\W_N^{\rm A}$ \cite{Gra99}.
Similarly, the noncolliding BESQ$^{(\nu)}$ 
is realized as an $h$-transform of the absorbing
BESQ$^{(\nu)}(t)$,
$\bR^{(\nu)}(t)=(R^{(\nu)}_1(t), \dots, R^{(\nu)}_N(t))$
in $\W_N^{+}$ \cite{KO01}.
For $\hXi(\cdot)=\Xi(\cdot)$ (resp. 
$\hXi(\cdot)=\Xi^{(\nu)}(\cdot)$), we set
$\V(\cdot)=\bB(\cdot)$ (resp.  $\V(\cdot)=\bR^{(\nu)}(\cdot)$),
and
$\W_N=\W_N^{\rm A}$ (resp. $\W_N^+$).
Put
\begin{equation}
\tau = \inf \{ t > 0 :
\V(t) \notin \W_N \}.
\label{eqn:tau1}
\end{equation}
Then under
$\xi=\sum_{j=1}^N \delta_{u_j}$, the equality
\begin{equation}
\E_{\xi} \left[
\prod_{m=1}^M g_{m}(\hX(t_m)) \right]
=\rE_{\u} \left[
\1(\tau > t_{M})
\prod_{m=1}^{M} g_m(\V(t_m)) 
\frac{h(\V(t_M))}{h(\u)} \right]
\label{eqn:harmonic1}
\end{equation}
is established.

\subsection{DMR for noncolliding diffusion processes
\label{sec:DMR_nonc}}

Now we can prove the following theorem.
\begin{thm}
\label{thm:DMR1}
The noncolliding BM and the noncolliding BESQ$^{(\nu)}$
with $\nu > -1$ have DMR for any $\xi \in \mM_0$,
$\xi(S) < \infty$.
\end{thm}
\noindent{\it Proof.} \quad
We introduce the stopping times
\begin{equation}
\tau_{jk}=\inf\{t> 0 : \hV_j(t)= \hV_k(t)\}, \quad
1\leq j < k \leq N.
\label{eqn:tauij}
\end{equation}
Let $\sigma_{jk} \in \cS_{N}$ be the permutation of $(j,k), 1 \leq j,k \leq N$.
Note that in a configuration $\u'$ if $u'_j=u'_k, j \not=k$, then
$\sigma_{jk}(\u')= \u'$, 
and the processes $\V(t)$ and $\sigma_{jk}(\V(t))$ are identical 
in distribution under the probability measure $\rP_{\u'}$.
By the strong Markov property of the process $\V(t)$
and by the fact that $h$ is anti-symmetric and $g_m, 1 \leq m \leq M$ are symmetric, 
$$
\rE_{\u} \left[
\1(\tau =\tau_{jk} < t_{M})
\prod_{m=1}^{M} g_m(\V(t_m)) 
\frac{h(\V(t_M))}{h(\u)} \right]=0.
$$
Since $\rP_{\u}(\tau_{jk}= \tau_{j' k'})=0$ if $(j,k)\not=(j',k')$,
and 
$$
\tau= \min_{1\leq j<k \leq N} \tau_{jk},
$$
$$
\rE_{\u} \left[
\1(\tau < t_{M})
\prod_{m=1}^{M} g_m(\V(t_m)) 
\frac{h(\V(t_M))}{h(\u)} \right]=0.
$$
Hence, (\ref{eqn:harmonic1}) equals
\begin{equation}
\rE_{\u} \left[
\prod_{m=1}^{M} g_m(\V(t_m)) 
\frac{h(\V(t_M))}{h(\u)} \right].
\label{eqn:harmonic2}
\end{equation}
By the equality (\ref{eqn:D6}), 
the theorem is concluded. \qed
\vskip 0.3cm

Then by Proposition \ref{thm:DM_det}, we will immediately
conclude the following.
\begin{cor}
\label{thm:K_BM}
The noncolliding BM is determinantal for any
$\xi \in \mM_0$, $\xi(\R) < \infty$.
The correlation kernel is given by
\begin{eqnarray}
\mbK_{\xi}(s,x;t,y)
&=&\int_{\R} \xi(dv) p(s,x|v) \cM_{\xi}^v(t,y)
-\1(s>t) p(s-t,x|y),
\nonumber\\
&& \qquad \qquad \qquad 
(s,x), (t,y) \in [0, \infty) \times \R,
\label{eqn:K_BM}
\end{eqnarray}
where
\begin{equation}
\cM_{\xi}^v(t,y)
=\int_{\R} d w \,
\frac{1}{\sqrt{2 \pi t}} e^{-(iy+w)^2/2t}
\Phi_{\xi}^v(iw).
\label{eqn:M_BM}
\end{equation}
\end{cor}
\begin{cor}
\label{thm:K_BESQ}
The noncolliding BESQ$^{(\nu)}$, $\nu>-1$
is determinantal for any
$\xi \in \mM_0^+$, $\xi(\R_+) < \infty$.
The correlation kernel is given by
\begin{eqnarray}
\mbK_{\xi}^{(\nu)}(s,x;t,y)
&=&\int_{\R} \xi(dv) p^{(\nu)}(s,x|v) \cM_{\xi}^v(t,y)
-\1(s>t) p^{(\nu)}(s-t,x|y),
\nonumber\\
&& \qquad \qquad \qquad 
(s,x), (t,y) \in [0, \infty) \times \R_+,
\label{eqn:K_BESQ}
\end{eqnarray}
where
\begin{equation}
\cM_{\xi}^v(t,y)
=\int_{\R_+} d w \,
\frac{1}{2t} \left(\frac{w}{y}\right)^{\nu/2}
e^{(y-w)/2t} J_{\nu} \left( \frac{\sqrt{wy}}{t} \right)
\Phi_{\xi}^v(-w).
\label{eqn:M_BESQ}
\end{equation}
\end{cor}

\subsection{CPR for noncolliding BM
\label{sec:CPR_noncBM}}

Let $W_j(t), t \geq 0, 1 \leq j \leq N$
be a collection of independent BM's on 
a probability space 
$(\check{\Omega}, \check{\cF}, \check{\rP}_{\0})$,
where the expectation w.r.t.$\check{\rP}_{\0}$
is written as $\check{\rE}_{\0}$.
Note that they are independent from
$B_j(t), t \geq 0, 1 \leq j \leq N$, 
and $W_j(0)=0, 1 \leq j \leq N$.
From the equality (\ref{eqn:mnZ2}). we have the following.
Introduce a set of independent complex BM's
\begin{equation}
Z_j(t)=B_j(t)+i W_j(t), \quad
1 \leq j \leq N, \quad t \geq 0, 
\label{eqn:Z1b}
\end{equation}
then
\begin{equation}
\sfM \left[ \left. \prod_{j=1}^N f_j(i W_j) \right|
\{(t_j, V_j(t_j))\}_{j=1}^N \right]
=\check{\rE}_{\0} \left[
\prod_{j=1}^N f_j(Z_j(t_j)) \right]
\label{eqn:cBM0}
\end{equation}
for polynomials $f_j$'s, 
and then the determinantal martingale (\ref{eqn:D1}) is written as
\begin{equation}
\cD_{\xi}^{\u}(T,\V(T))
=\check{\rE}_{\0} \left[
\det_{1 \leq j, k \leq N}
\left[ \Phi_{\xi}^{u_k}(Z_j(T)) \right] \right].
\label{eqn:cBM0b}
\end{equation}
Then Theorem \ref{thm:DMR1} 
is transformed into the following.

Let $Z_j(t), 1 \leq j \leq N, t \geq 0$
be a set of independent complex BM's
given by (\ref{eqn:Z1b}).
If they start at $Z_j(0)=u_j \in \R, 1 \leq j \leq N$,
the probability space is denoted by
$(\Omega, \cF, \bP_{\u})$
with $\u=(u_1, \dots, u_N)$.
The space $(\Omega, \cF, \bP_{\u})$
is a product of two probability spaces
$(\Omega, \cF, \rP_{\u})$ for $B_j(\cdot) =\Re Z_j(\cdot)$
and $(\check{\Omega}, \check{\cF}, \check{\rP}_{\0})$
for $W_j(\cdot)=\Im Z_j(\cdot), 1\leq j \leq N$.
The expectation w.r.t.$\bP_{\u}$ is denoted by $\bE_{\u}$.

\begin{cor}
\label{thm:CPR_noncBM}
The noncolliding BM has CPR (\ref{eqn:CPR1})
with the complex BMs (\ref{eqn:Z1b}) and
\begin{equation}
\varphi^{u}_{\xi}(\cdot)
=\Phi_{\xi}^{u}(\cdot), \quad u \in \C,
\quad \xi \in \mM_0.
\label{eqn:CPR_noncBM}
\end{equation}
\end{cor}

This result was given as Theorem 1.1 in \cite{KT13},
where the present CPR is called the
{\it complex BM representation}.

\subsection{CPR for noncolliding BES$^{(\nu)}$
\label{sec:CPR_noncBES}}

We introduce a set of independent complex diffusions
\begin{equation}
Z^{(\nu)}_j(t)=\widetilde{R}^{(\nu)}_j(t) + i W_j(t),
\quad 1 \leq j \leq N, \quad t \geq 0,
\label{eqn:Z2b}
\end{equation}
where $W_j(t), 1 \leq j \leq N$ are
independent BM's on the probability space
$(\check{\Omega}, \check{\cF}, \check{\rP}_{\0})$.
For BES and BESQ with odd dimensions, $D=2n+3, n \in \N_0$,
the indices are half-odds
$\nu=D/2-1=n+1/2, n \in \N_0$.
Let $\widetilde{f}(z)$ be a polynomial
of $z^2, z \in \C$.
Then Lemma \ref{thm:BES_mart1} gives the following.
For $n \in \N_0$,
\begin{eqnarray}
&&
\widetilde{\sfM}^{(n+1/2)} \left[ \left. \prod_{j=1}^N\widetilde{f}_j( i W_j ) \right|
\left\{(t_j, \widetilde{R}^{(n+1/2)}_j(t_j)) \right\}_{j=1}^N \right]
\nonumber\\
&& \qquad 
=\check{\rE}_{\0} \left[
\prod_{j=1}^N 
Q^{(n+1/2)}_{t_j}(Z^{(n+1/2)}_j (t_j)) \widetilde{f}_j( Z^{(n+1/2)}_j(t_j) ) \right],
\label{eqn:c_process1}
\end{eqnarray}
where $Q^{(n+1/2)}_t$ is given by (\ref{eqn:Q1}). 
Then we have 
\begin{equation}
\cD_{\xi}^{(n+1/2), \u}(T, \V^{(n+1/2)}(T))
=\check{\rE} \left[
\det_{1 \leq j, k \leq N}
\left[ Q_T^{(n+1/2)}(Z^{(n+1/2)}_j(T))
\widetilde{\Phi}_{\xi}^{u_k}
(Z^{(n+1/2)}_j(T)) \right] \right],
\label{eqn:c_process1b}
\end{equation}
where
\begin{equation}
\widetilde{\Phi}_{\xi}^{v}(z)
=\prod_{r \in \supp \xi \cap \{v, -v \}^{\rm c}}
\frac{z^2-r^2}{v^2-r^2},
\quad z, v \in \C.
\label{eqn:Phi2z}
\end{equation}

Then Theorem \ref{thm:DMR1} 
is transformed into the following.
\begin{cor}
\label{thm:CPR_noncBES}
The noncolliding BES$^{(n+1/2)}, n \in \N_0$ has CPR (\ref{eqn:CPR1})
with the complex BMs (\ref{eqn:Z1b}) and
\begin{equation}
\varphi^{u}_{\xi}(\cdot)
=Q_t^{(n+1/2)}(\cdot)
\widetilde{\Phi}_{\xi}^{u}(\cdot), \quad u \in \C,
\quad \xi \in \mM_{+,0}.
\label{eqn:CPR_noncBES}
\end{equation}
\end{cor}

\subsection{Martingales for configurations with multiple points
\label{sec:multi}}

We generalize the function 
(\ref{eqn:Phi1}) as following.
Depending on the transition probability density
of a process, 
$p(s,x|v)$, $0 < s < \infty, x, v \in S$ we put
\begin{equation}
\phi^u_{\xi}((s,x); z, \zeta)=
\frac{p(s,x|\zeta)}{p(s,x|u)}
\frac{1}{z-\zeta} \prod_{r \in \supp \xi}
\left( \frac{z-r}{\zeta-r}\right)^{\xi(\{r\})},
\quad z, \zeta \in \C.
\label{eqn:Phi1b}
\end{equation}
Let $C(\delta_{u})$ be a closed contour on the complex plane $\C$
encircling a point $u$ on $S$
once in the positive direction and set
\begin{eqnarray}
\Phi_{\xi}^{u}((s,x); z) &=& \frac{1}{2 \pi i}
\oint_{C(\delta_{u})} d \zeta \, 
\phi^u_{\xi}((s,x); z, \zeta)
\nonumber\\
&=& {\rm Res} \,\Big[\phi^u_{\xi}((s,x); z, \zeta); \zeta=u \Big].
\label{eqn:Phi1c}
\end{eqnarray}
This function is defined 
for any finite configuration $\xi$, in which
there can be multiple points in general.
(If there is no multiple point, 
(\ref{eqn:Phi1c}) is reduced to (\ref{eqn:Phi1}).) 
Since (\ref{eqn:Phi1c}) is a polynomial with respect to $z$,
we can extend $\cM_{\xi}^{u}(t, y)$ to 
\begin{equation}
\cM_{\xi}^{u}((s, x)|(t,y))
=\sfM \left[\Phi_{\xi}^{u}((s,x); c W) \Big| (t,y) \right],
\quad (s,x), (t,y) \in [0, \infty) \times S.
\label{eqn:detM5}
\end{equation}
Let
\begin{equation}
\xi_{\rm s}(\cdot)=\sum_{u \in \supp \xi} \delta_u (\cdot).
\label{eqn:xidef2}
\end{equation}
Then the correlation kernel (\ref{eqn:K1}) is generalized 
to 
\begin{eqnarray}
\mbK_{\xi}(s,x;t,y)
&=& \int_{S} \xi_{\rm s}(dv)
p(s,x|v) \cM_{\xi}^{v}((s,x) | (t, y))
-\1(s>t) p(s-t, x|y),
\nonumber\\
&& \qquad \qquad (s,x), (t,y) \in [0, \infty) \times S.
\label{eqn:kernel1b}
\end{eqnarray}
By definition of the present martingales 
$\cM_{\xi}^{\cdot}((\cdot, \cdot)|(\cdot, \cdot))$,
it is written by double integral as
\begin{eqnarray}
\mbK_{\xi}(s,x;t,y)
&=& \frac{1}{2 \pi i} \oint_{C(\xi)} d \zeta \,
p(s,x|\zeta) \int_{S} d w \,
p(t,w|c^{-1} y)
\frac{1}{c w-\zeta}
\prod_{r \in \supp \xi}
\left( \frac{c w-r}{\zeta-r} \right)^{\xi(\{r\})} 
\nonumber\\
&& \quad 
-\1(s>t) p(s-t, x|y),
\qquad (s,x), (t,y) \in [0, \infty) \times S,
\label{eqn:kernel1c}
\end{eqnarray}
where $C(\xi)$ denotes a counterclockwise contour
on $\C$ encircling the points in $\supp \xi$ on $S$
but not point $c w, w \in S$;
$C(\xi)=\sum_{v \in \supp \xi} C(\delta_{v})$.

\begin{cor}
\label{thm:K_general}
The noncolliding BM and 
the noncolliding BESQ$^{(\nu)}$ are determinantal for 
any initial configuration with fine number of particles;
$\xi \in \mM$, $\xi(\R) < \infty$,
or $\xi \in \mM^+$, $\xi(\R_+) < \infty$
with correlation kernels (\ref{eqn:kernel1c}).
\end{cor}

In the papers \cite{KT10,KT11}, 
we proved this statement by
deriving the double integral representations (\ref{eqn:kernel1c}) for
the spatio-temporal correlation kernels.
There we used the multiple orthogonal polynomials
\cite{KT09,KT10,KT11} in order to obtain
the expression (\ref{eqn:kernel1c}).
As shown in this lectures, however, 
it is not necessary to use 
multiple orthogonal polynomials
to obtain the results.

As an example, we consider the extreme case such that
all $N$ points are concentrated on an origin,
\begin{equation}
\xi=N \delta_0 \quad
\Longleftrightarrow \quad
\xi_{\rm s}=\delta_0 \quad
\mbox{with} \quad
\xi(\{0\})=N.
\label{eqn:xi01}
\end{equation}
In this case (\ref{eqn:Phi1b}) and (\ref{eqn:Phi1c})
become
\begin{eqnarray}
\phi_{N \delta_0}^0((s,x); z, \zeta)
&=& \frac{p(s,x|\zeta)}{p(s,x|0)}
\frac{1}{z-\zeta} \left( \frac{z}{\zeta} \right)^N
\nonumber\\
&=& \frac{p(s,x|\zeta)}{p(s,x|0)}
\sum_{\ell=0}^{\infty} \frac{z^{N-\ell-1}}{\zeta^{N-\ell}},
\label{eqn:xi02}
\end{eqnarray}
and
\begin{eqnarray}
\Phi_{N \delta_0}^0((s,x);z) &=& 
\frac{1}{p(s,x|0)} \sum_{\ell=0}^{\infty} z^{N-\ell-1}
\frac{1}{2 \pi i} \oint_{C(\delta_0)} d \zeta \,
\frac{p(s,x|\zeta)}{\zeta^{N-\ell}}
\nonumber\\
&=& 
\frac{1}{p(s,x|0)} \sum_{\ell=0}^{N-1} z^{N-\ell-1}
\frac{1}{2 \pi i} \oint_{C(\delta_0)} d \zeta \,
\frac{p(s,x|\zeta)}{\zeta^{N-\ell}},
\label{eqn:xi03}
\end{eqnarray}
since the integrands are holomorphic when $\ell \geq N$,
where we have assumed $\nu > -1$ for BESW$^{(\nu)}$. 

For BM with the transition probability density (\ref{eqn:p_BM}),
(\ref{eqn:xi03}) gives
\begin{eqnarray}
\Phi_{N \delta_0}^0((s,x);z)
&=& \sum_{\ell=0}^{N-1} z^{N-\ell-1}
\frac{1}{2 \pi i}
\oint_{C(\delta_0)} d \zeta \,
\frac{e^{x\zeta/s-\zeta^2/2s}}{\zeta^{N-\ell}}
\nonumber\\
&=& \sum_{\ell=0}^{N-1} 
\left( \frac{z}{\sqrt{2s}} \right)^{N-\ell-1}
\frac{1}{2 \pi i}
\oint_{C(\delta_0)} d \eta \,
\frac{e^{2(x/\sqrt{2s}) \eta-\eta^2}}{\eta^{N-\ell}}
\nonumber\\
&=& \sum_{\ell=0}^{N-1} 
\left( \frac{z}{\sqrt{2s}} \right)^{N-\ell-1}
\frac{1}{(N-\ell-1)!} H_{N-\ell-1} \left(\frac{x}{\sqrt{2s}} \right),
\label{eqn:Hermite1B}
\end{eqnarray}
where we have used the contour integral representation
of the Hermite polynomials \cite{Sze91}
\begin{equation}
H_n(x)=\frac{n!}{2 \pi i} \oint_{C(\delta_0)} d \eta
\frac{e^{2 x \eta-\eta^2}}{\eta^{n+1}},
\quad n \in \N_0, \quad x \in \R.
\label{eqn:Hermite2B}
\end{equation}
Thus its integral transform is calculated as
\begin{eqnarray}
&& \sfM \left[ \left.
\Phi_{N \delta_0}^0((s,x);i W) \right| (t, y) \right]
\nonumber\\
&& \quad 
= \sum_{\ell=0}^{N-1} \frac{1}{(N-\ell-1)!} H_{N-\ell-1}
\left( \frac{x}{\sqrt{2s}} \right)
\frac{1}{(2s)^{(N-\ell-1)/2}}
\sfM[(iW)^{N-\ell-1}|(t,y)]
\nonumber\\
&& \quad 
= \sum_{\ell=0}^{N-1} \frac{1}{(N-\ell-1)!} H_{N-\ell-1}
\left( \frac{x}{\sqrt{2s}} \right)
\frac{1}{(2s)^{(N-\ell-1)/2}} 
m_{N-\ell-1}(t,y)
\nonumber\\
&& \quad 
= \sum_{\ell=0}^{N-1} \frac{1}{(N-\ell-1)! 2^{N-\ell-1}} 
\left( \frac{t}{s} \right)^{(N-\ell-1)/2}
H_{N-\ell-1}
\left( \frac{x}{\sqrt{2s}} \right)
H_{N-\ell-1}
\left( \frac{y}{\sqrt{2t}} \right),
\nonumber
\end{eqnarray}
where we have used Lemma \ref{thm:Hermite1}.
Then we obtain the following,
\begin{eqnarray}
&& \cM_{N \delta_0}^0((s,x)|(t,B(t)))
=\sum_{n=0}^{N-1} \frac{1}{n! 2^n} m_n(s,x) m_n(t,B(t))
\nonumber\\ 
&& \qquad \qquad = \sum_{n=0}^{N-1} \frac{1}{n! 2^n} \left(\frac{t}{s}\right)^{n/2}
H_n \left( \frac{x}{\sqrt{2s}} \right)
H_n \left( \frac{B(t)}{\sqrt{2t}} \right)
\nonumber\\
&& \qquad \qquad =
\sqrt{\pi} e^{x^2/4s+B(t)^2/4t}
\sum_{n=0}^{N-1} 
 \left(\frac{t}{s}\right)^{n/2}
\varphi_n \left( \frac{x}{\sqrt{2s}} \right)
\varphi_n \left( \frac{B(t)}{\sqrt{2t}} \right),
\label{eqn:Hermite3B}
\end{eqnarray}
where
$$
\varphi_n(x)=\frac{1}{\sqrt{ \sqrt{\pi} 2^n n!}}
H_n(x) e^{-x^2/2}, \quad
n \in \N, \quad x \in \R.
$$
Similarly, for BESQ$^{(\nu)}, \nu > -1$
with the transition probability density (\ref{eqn:p_BESQ}),
we obtain
\begin{eqnarray}
\Phi_{N \delta_0}^{(\nu), 0}((s,x); z)
&=& \frac{(2s)^{\nu} \Gamma(\nu+1)}{x^{\nu/2}}
\sum_{\ell=0}^{N-1} z^{N-\ell-1}
\frac{1}{2 \pi i} \oint_{C(\delta_0)} d \zeta \,
\frac{e^{-\zeta/2s}}{\zeta^{N-\ell+\nu/2}}
I_{\nu} \left( \frac{\sqrt{x \zeta}}{s}\right)
\nonumber\\
&=& \Gamma(\nu+1) \sum_{\ell=0}^{N-1}
\left( -\frac{z}{2s} \right)^{N-\ell-1}
\frac{1}{\Gamma(N-\ell+\nu)}
L^{(\nu)}_{N-\ell-1} \left( \frac{x}{2s} \right),
\label{eqn:Laguerre1B}
\end{eqnarray}
where we used the contour integral representation
of the Laguerre polynomials
\begin{equation}
L_n^{(\nu)}(x)
= \frac{\Gamma(n+\nu+1)}{x^{\nu/2}}
\frac{1}{2 \pi i} \oint_{C(\delta_0)}
d \eta \,
\frac{e^{\eta}}{\eta^{n+1+\nu/2}} J_{\nu}(2 \sqrt{\eta x})
\label{eqn:Laguerre2B}
\end{equation}
with the relation
$I_{\nu}(iz)=(-1)^{\nu/2} J_{\nu}(z), - \pi <
{\rm arg}(z) \leq \pi/2$.
By using Lemma \ref{thm:Laguerre1}, we have
\begin{eqnarray}
&& \cM_{N \delta_0}^{(\nu), 0}((s,x)|(t, R^{(\nu)}(t))
\nonumber\\
&& \quad = \sfM^{(\nu)} \left[ \left.
\Phi_{N \delta_0}^{(\nu), 0}((s,x); -W) \right| (t, R^{(\nu)}(t)) \right]
\nonumber\\
&& \quad = \Gamma(\nu+1)
\sum_{n=0}^{N-1} \frac{1}{\Gamma(n+1) \Gamma(n+\nu+1) (2s)^{2n}}
m_n^{(\nu)}(s,x) m_n^{(\nu)} \left(t, R^{(\nu)}(t) \right)
\nonumber\\
&& \quad = \Gamma(\nu+1)
\sum_{n=0}^{N-1} \frac{\Gamma(n+1)}{\Gamma(n+\nu+1)}
\left( \frac{t}{s} \right)^n
L_n^{(\nu)} \left( \frac{x}{2s} \right)
L_n^{(\nu)} \left( \frac{R^{(\nu)}(t)}{2t} \right)
\nonumber\\
&& \quad =
\Gamma(\nu+1) 
\left(\frac{x}{2s} \right)^{-\nu/2}
\left(\frac{R^{(\nu)}(s)}{2s} \right)^{-\nu/2}
e^{x/4s+R^{(\nu)}(t)/4t}
\nonumber\\
&& \qquad \qquad \times
\sum_{n=0}^{N-1} \frac{\Gamma(n+1)}{\Gamma(n+\nu+1)}
\left( \frac{t}{s} \right)^n
\varphi_n^{(\nu)} \left( \frac{x}{2s} \right)
\varphi_n^{(\nu)} \left( \frac{R^{(\nu)}(t)}{2t} \right),
\label{eqn:Laguerre3B}
\end{eqnarray}
where
$$
\varphi_n^{(\nu)}(x)=\sqrt{\frac{\Gamma(n+1)}{\Gamma(n+\nu+1)}}
x^{\nu/2} L_n^{(\nu)}(x) e^{-x/2},
\quad n \in \N_0, \quad x \in \R_+.
$$
The processes (\ref{eqn:Hermite3B}) and (\ref{eqn:Laguerre3B})
are local martingales and
\begin{equation}
\rE_0 \left[ \cM_{N \delta_0}^0((s,x)|(t, \V(t))) \right]
= \rE_0 \left[ \cM_{N \delta_0}^0((s,x)|(0, \V(0))) \right]=1
\label{eqn:HL1}
\end{equation}
for $0 < t \leq T < \infty, 
(s,x) \in [0, T] \times S$.

By the formula (\ref{eqn:kernel1}), we obtain the
correlation kernels as
\begin{eqnarray}
\mbK_{N \delta_0}(s,x;t,y) &=& p(s,x|0) \cM_{N \delta_0}^0((s,x)| (t,y))
\nonumber\\
&=& \frac{e^{-x^2/4s}}{e^{-y^2/4t}}
\bK_{H}(s,x;t,y)
\label{eqn:HermiteB1}
\end{eqnarray}
with
\begin{equation}
\bK_{H}(s,x;t,y)
= \frac{1}{\sqrt{2s}}
\sum_{n=0}^{N-1} \left( \frac{t}{s} \right)^{n/2}
\varphi \left( \frac{x}{\sqrt{2s}} \right)
\varphi \left( \frac{y}{\sqrt{2t}} \right)
-\1(s>t) p(s-t, x|y),
\label{eqn:Hermite5B}
\end{equation}
and
\begin{eqnarray}
\mbK^{(\nu)}_{N \delta_0}(s,x;t,y) 
&=& p^{(\nu)}(s,x|0) \cM_{N \delta_0}^{(\nu),0} ((s,x)| (t,y))
\nonumber\\
&=& 
\frac{(x/2s)^{\nu/2} e^{-x/4s}}{(y/2t)^{\nu/2} e^{-y/4t}}
\bK_{L^{(\nu)}}(s,x;t,y)
\label{eqn:Laguerre5B}
\end{eqnarray}
with
\begin{equation}
\bK_{L^{(\nu)}}(s,x;t,y)= 
\frac{1}{2s}
\sum_{n=0}^{N-1} \left( \frac{t}{s} \right)^{n}
\varphi^{(\nu)} \left( \frac{x}{\sqrt{2s}} \right)
\varphi^{(\nu)} \left( \frac{y}{\sqrt{2t}} \right)
-\1(s>t) p^{(\nu)}(s-t, x|y),
\label{eqn:Laguerre6B}
\end{equation}
where $\bK_{H}(\cdot; \cdot)$ and $\bK_{L^{(\nu)}}(\cdot; \cdot)$
are known as the {\it extended Hermite and Laguerre kernels},
respectively (see, for instance, \cite{For10}).
Here we would like to emphasize the fact that
these kernels have been derived here by not following
any `orthogonal-polynomial arguments'
but by only using proper martingales
associated with the chosen initial configuration (\ref{eqn:xi01}).
In this special case, they are expressed by the
Hermite and Laguerre polynomials in the forms
(\ref{eqn:Hermite3B}) and (\ref{eqn:Laguerre3B}),
respectively.
In the present new approach,
the martingale properties (\ref{eqn:HL1})
and the reducibility (\ref{eqn:reducibility})
coming from the independence of diffusion processes
play essential roles instead of orthogonality
in the theory of orthogonal ensembles in
random matrix theory \cite{Meh04,For10}.

\subsection{Martingales associated with infinite particle systems
\label{sec:IPS2}}
In \cite{KT10} 
we gave useful sufficient conditions of $\xi$ 
so that the noncolliding BM, $(\Xi(t), \P_{\xi})$ is well defined
as a determinantal process even if $\xi(\R)=\infty$.
For $L>0, \alpha>0$ and $\xi\in\mM$ we put
\begin{equation}
M(\xi, L)=\int_{[-L,L]\setminus\{0\}} \frac{\xi(dx)}{x},
\qquad
M_\alpha(\xi, L)
=\left( \int_{[-L,L]\setminus\{0\}} 
\frac{\xi(dx)}{|x|^\alpha}\right)^{1/\alpha},
\label{eqn:defM}
\end{equation}
and
\begin{equation}
M(\xi) = \lim_{L\to\infty}M(\xi, L),
\quad
M_{\alpha}(\xi)= \lim_{L\to\infty}M_\alpha(\xi, L),
\label{eqn:defM2}
\end{equation}
if the limits finitely exist. Then
\vskip 3mm
\noindent ({\bf C.1})
there exists $C_0 > 0$ such that
$|M(\xi,L)|  < C_0$, $L>0$,

\vskip 3mm

\noindent ({\bf C.2}) (i) 
there exist $\alpha\in (1,2)$ and $C_1>0$ such that
$
M_\alpha(\xi) \le C_1,
$ \\
\noindent (ii) 
there exist $\beta >0$ and $C_2 >0$ such that
$$
M_1(\tau_{-a^2} \xi^{\langle 2 \rangle}) \leq C_2
(\max\{|a|, 1\})^{-\beta}
\quad \forall a \in \supp \xi.
$$
\vskip 3mm

\noindent
It was shown that, 
if $\xi \in \mM_0$ satisfies the 
conditions $({\bf C.1})$ and $({\bf C.2})$, 
then for $a \in \R$ and $z \in \C$, 
\begin{equation}
\Phi_{\xi}^{a}(z)\equiv \lim_{L\to\infty}
\Phi_{\xi \cap [a-L, a+L]}^{a}(z) \quad
\mbox{finitely exists},
\label{eqn:CCC1}
\end{equation}
and
\begin{equation}
|\Phi_{\xi}^{a}(z)|\le C \exp \bigg\{c(|a|^\theta +|z|^\theta) \bigg\}
\left| \frac{z}{a}\right|^{\xi(\{0\})} \left|\frac{a}{a-z} \right|,
\quad a\in \supp \xi, \ z \in \C,
\label{eqn:CCC2}
\end{equation}
for some $c, C>0$ and $\theta\in 
(\max\{\alpha, (2-\beta)\},2)$, 
which are determined by the constants $C_0, C_1, C_2$
and the indices $\alpha, \beta$ 
in the conditions (Lemma 4.4 in \cite{KT10}).
We have noted that in the case that $\xi \in \mM_0$ 
satisfies the conditions 
$({\bf C.1})$ and $({\bf C.2})$
with constants $C_0, C_1, C_2$ and indices $\alpha$ and $\beta$,
then $\xi \cap [-L, L], \forall L > 0$ does as well.
Hence we can obtain the convergence of moment generating functions
\begin{equation}
{\Psi}_{\xi \cap [-L, L]}^{\t}[\f] \to
{\Psi}_{\xi}^{\t}[\f]
\quad \mbox{as} \quad L \to \infty,
\label{eqn:CCC3}
\end{equation}
which implies the convergence of
the probability measures 
\begin{equation}
\P_{\xi \cap [-L, L]}
\to \P_{\xi}
\quad \mbox{in $L \to \infty$}
\label{eqn:CCC4}
\end{equation}
in the sense of finite dimensional distributions.
Moreover,
even if $\xi(\R)=\infty$, $\mbK_{\xi}$ 
is well-defined as a correlation kernel and
dynamics of the noncolliding BM with an infinite number of
particles $(\Xi(t), \P_{\xi})$ exists 
as a determinantal process \cite{KT10}.

Similarly, in \cite{KT11}, the following sufficient conditions
for initial configurations $\xi \in \mM^+$ were given so that
the noncolliding BESQ, $(\Xi^{(\nu)}(t), \P^{(\nu)}_{\xi}), \nu >-1$
is well-defined as determinantal processes even if 
$N=\xi(\R_+)= \infty$.

\vskip 0.3cm
\noindent
{\bf (C.A)} (i)
There exists $\alpha \in (1/2,1)$ and $C_1 > 0$
such that 
$M_{\alpha}(\xi) \leq C_1$. 

(ii) There exist $\beta > 0$ and $C_2 > 0$
such that
$$
M_1(\tau_{-a} \xi) 
\leq C_2 (|a| \vee 1)^{-\beta},
\quad \forall a \in \supp \xi.
$$
\vskip 0.3cm

The families of $\xi$ satisfying the conditions
are denoted by $\hat{\mX}=\mX$ for the noncolliding BM
and $\hat{\mX}=\mX^+$ for the noncolliding BESQ$^{(\nu)}$,
respectively.

\begin{prop}
\label{thm:infinite1}
Suppose that $0 < t \leq T < \infty$.
Then the noncolliding BM, $(\Xi(t), \P_{\xi})$, 
started at $\xi \in \mX_0$ and
the noncolliding BESQ$^{(\nu)}$, 
$(\Xi^{(\nu)}(t), \P^{(\nu)}_{\xi}), \nu >-1$, started at $\xi \in \mX^+_0$ 
have DMR for any $\cF(t)$-measurable polynomial function
also in the case with $N=\xi(S)=\infty$.
\end{prop}
\vskip 0.3cm
For $(\Xi(t), \P_{\xi})$, the similar statement
was proved for the complex BM representation
in \cite{KT13} (Corollary 1.3).
Here by the reducibility of the determinantal
martingale given by Lemma \ref{thm:reducibility},
this proposition is readily concluded.

There are two interesting examples
of local martingales for infinite particle systems.
First we consider the configuration
\begin{equation}
\xi_{\Z}(\cdot)=\sum_{j \in \Z} \delta_j(\cdot),
\label{eqn:xiZ}
\end{equation}
that is, the configuration in which every integer point
$\Z$
is occupied by one particle.
It is easy to confirm that $\xi_{\Z} \in \mX_0$
and the noncolliding BM started at $\xi_{\Z}$,
$(\Xi(t), \P_{\xi_{\Z}})$, is a determinantal process
with an infinite number of particles \cite{KT10}.
Since $\prod_{n \in \N}(1-x^2/n^2)=\sin(\pi x)/(\pi x)$,
\begin{eqnarray}
\Phi_{\xi_{\Z}}^{u}(z)
&=& \prod_{r \in \Z, r \not= u}
\frac{z-r}{u-r}
\nonumber\\
&=& \frac{\sin \{\pi(z-u)\}}{\pi (z-u)}
= \frac{1}{2 \pi} \int_{-\pi}^{\pi} d \lambda \,
e^{i \lambda (z-u)},
\quad z, u \in \C.
\label{eqn:PhiZ1}
\end{eqnarray}
Its integral transform is calculated as
\begin{eqnarray}
\sfM \left[ \left. \Phi_{\xi_{\Z}}^{u}(iw) 
\right| (t,x) \right]
&=& \int_{\R} dw \, q(t,w|x) \Phi_{\xi_{\Z}}^{u}(iw)
\nonumber\\
&=& \int_{-\infty}^{\infty} dw \,
\frac{1}{\sqrt{2 \pi t}} e^{-(ix+w)^2/2t}
\Phi_{\xi_{\Z}}^{u}(iw)
\nonumber\\
&=& \frac{1}{2 \pi} \int_{-\pi}^{\pi} d \lambda \,
e^{t \lambda^2/2 + i \lambda(x-u)}.
\label{eqn:mZ0}
\end{eqnarray}
Then we have local martingales
\begin{equation}
\cM_{\xi_{\Z}}^{k}(t, B_j(t))
=\frac{1}{2\pi}
\int_{-\pi}^{\pi} d \lambda \,
\exp \left\{ \frac{\lambda^2}{2} t
+ i \lambda (B_j(t)-k) \right\},
\quad j, k \in \Z, \quad
0 < t \leq T < \infty.
\label{eqn:mZ1}
\end{equation}
We see that
\begin{eqnarray}
\rE_{\xi_{\Z}} \left[
\cM_{\xi_{\Z}}^{k}(t, B_j(t)) \right]
&=& \rE_{\xi_{\Z}} \left[
\cM_{\xi_{\Z}}^{k}(0, B_j(0)) \right]
\nonumber\\
&=& \delta_{j k}, \qquad \qquad
0 < t \leq T < \infty.
\label{eqn:mZ3}
\end{eqnarray}

If $\nu > -1$, the Bessel function $J_{\nu}(z)$ given by
(\ref{eqn:J1}) has an infinite number of pairs of
positive and negative zeros with the same absolute value,
which are all simple. We write the positive zeros
of $J_{\nu}(z)$ arranged in ascending order of the
absolute values as
\begin{equation}
0 < j_{\nu,1} < j_{\nu, 2} < j_{\nu, 3} < \cdots.
\label{eqn:Jzero}
\end{equation}
Then, $J_{\nu}(z)$ has the following infinite product
expression \cite{Wat44},
\begin{equation}
J_{\nu}(z)=\frac{(z/2)^{\nu}}{\Gamma(\nu+1)}
\prod_{j=1}^{\infty} \left( 1 - \frac{z^2}{j_{\nu,j}^2} \right).
\label{eqn:J2}
\end{equation} 
For the noncolliding BESQ$^{(\nu)}$,
we consider the initial configuration
in which every point of the squares of positive zeros
of $J_{\nu}(z)$ is occupied by one particle,
which is denoted as
\begin{equation}
\xi_{J_{\nu}}^{\langle 2 \rangle}(\cdot)
=\sum_{j=1}^{\infty} \delta_{j_{\nu, j}^2}(\cdot).
\label{eqn:xiJ1}
\end{equation}
We can see that $\xi_{J_{\nu}}^{\langle 2 \rangle} \in \mX^+_0$
and thus $(\Xi^{(\nu)}(t), \P_{\xi_{J_{\nu}}^{\langle 2 \rangle}}^{(\nu)})$
is a determinantal process with an infinite number of particles
\cite{KT11}.
For $k \in \N$ we find that 
\begin{eqnarray}
\Phi_{\xi_{J_{\nu}}^{\langle 2 \rangle}}^{(j_{\nu, k})^2}(z)
&=& \left( \frac{ (j_{\nu, k})^2}{z} \right)^{\nu/2}
\frac{1}{(J_{\nu+1}(j_{\nu,k}))^2}
\int_0^1 d \lambda \,
J_{\nu}(\sqrt{\lambda z}) J_{\nu}(\sqrt{\lambda} j_{\nu, k}),
\label{eqn:PhiJ1}
\end{eqnarray}
and their integral transforms gives the martingales,
\begin{eqnarray}
\cM_{\xi_{J_{\nu}}^{\langle 2 \rangle}}^{(j_{\nu, k})^2}
(t, R_j^{(\nu)}(t))
&=& \sfM^{(\nu)} \left[ \left. \Phi_{\xi_{J_{\nu}}^{\langle 2 \rangle}}^{(j_{\nu, k})^2}(-W) 
\right| (t,R_j^{(\nu)}) \right]
\nonumber\\
&=& \left( \frac{(j_{\nu,k})^2}{R^{(\nu)}_j(t)} \right)^{\nu/2}
\frac{1}{(J_{\nu+1}(j_{\nu,k}))^2}
\int_0^1 d \lambda \, e^{\lambda t/2}
J_{\nu}\left(\sqrt{\lambda R_j^{(\nu)}(t)} \right) 
J_{\nu}(\sqrt{\lambda} j_{\nu, k}),
\nonumber\\
&& \qquad
\quad j, k \in \N, \quad 0 < t \leq T \leq \infty.
\label{eqn:mJ1}
\end{eqnarray}
We see that for $0 < t \leq T < \infty$,
\begin{eqnarray}
\rE_{\xi_{J_{\nu}}^{\langle 2 \rangle}}^{(\nu)} \left[
\cM_{\xi_{J_{\nu}}^{\langle 2 \rangle}}^{(j_{\nu, k})^2}
(t, R_j^{(\nu)}(t))
\right]
&=&
\rE_{\xi_{J_{\nu}}^{\langle 2 \rangle}}^{(\nu)} \left[
\cM_{\xi_{J_{\nu}}^{\langle 2 \rangle}}^{(j_{\nu, k})^2}
(0, R_j^{(\nu)}(0))
\right]
\nonumber\\
&=& \delta_{j k}.
\label{eqn:mJ2}
\end{eqnarray}

By the formula (\ref{eqn:K1}),
these martingales determine the correlation kernels,
which are denoted as
$\mbK_{\xi_{\Z}}$ and $\mbK^{(\nu)}_{\xi_{J_{\nu}}^{\langle 2 \rangle}}$.
In the previous papers \cite{KT10,KT11},
we showed 
\begin{eqnarray}
&&
\lim_{\tau \to \infty} \mbK_{\xi_{\Z}}(s+\tau,x;t+\tau,y)
=\bK_{\rm sin}(t-s, y-x),
\nonumber\\
&&
\lim_{\tau \to \infty} \mbK^{(\nu)}_{\xi_{J_{\nu}}^{\langle 2 \rangle}}(s+\tau,x;t+\tau,y)
= \left( \frac{x}{y} \right)^{\nu/2} \bK_{\rm J_{\nu}}(t-s, y|x),
\label{eqn:relax1}
\end{eqnarray}
and proved that the noncolliding BM started 
at (\ref{eqn:xiZ}) and the noncolliding BESQ$^{(\nu)}$
started at (\ref{eqn:xiJ1}) converge in the long-term limit
to the equilibrium determinantal processes governed by
the extended sine kernel 
\begin{eqnarray}
{\bf K}_{\sin}(t, x)
&=& \left\{ \begin{array}{ll} 
\displaystyle{
\int_{0}^{1} d\lambda \, e^{\pi^2 \lambda^2 t/2} 
\cos (\pi \lambda x)},
& \mbox{if $t>0 $} \cr
& \cr
\displaystyle{
\frac{\sin(\pi x)}{\pi x}}
& \mbox{if $t=0$} \cr
& \cr
\displaystyle{
- \int_{1}^{\infty} d\lambda \, 
e^{\pi^2 \lambda^2 t/2} \cos (\pi \lambda x)},
& \mbox{if $t<0$},
\end{array} \right.
\label{eqn:sine-kernel}
\end{eqnarray}
and the extended Bessel kernel \cite{For10}
\begin{equation}
\bK_{J_{\nu}}(t,y|x) = \left\{
   \begin{array}{ll}
\displaystyle{
\frac{1}{4} \int_{0}^{1} d \lambda \,
e^{\lambda t/2} J_{\nu}(\sqrt{\lambda x})
J_{\nu}(\sqrt{\lambda y})
},
& \mbox{if} \quad t >0  \\
& \\
\displaystyle{
\frac{J_{\nu}(\sqrt{x}) \sqrt{y} J_{\nu}'(\sqrt{y})
-\sqrt{x} J_{\nu}'(\sqrt{x}) J_{\nu}(\sqrt{y})}{2(x-y)}
},
& \mbox{if} \quad t=0 \\
& \\
\displaystyle{
- \frac{1}{4}\int_{1}^{\infty} d \lambda \,
e^{\lambda t/2} J_{\nu}(\sqrt{\lambda x})
J_{\nu}(\sqrt{\lambda y})
},
& \mbox{if} \quad t <0,
   \end{array} \right. 
\label{eqn:KBessel1}
\end{equation}
respectively.
These {\it relaxation phenomena} of infinite particle systems
are caused by the following properties of the present martingales,
\begin{eqnarray}
&& \lim_{\tau \to \infty}
\sum_{k \in \Z} p(\tau, x|k) \cM_{\xi_{\Z}}^{k}(t+\tau, B(t))
=\cM_{\xi_{\Z}}^{x}(t, B(t)), \quad x \in \R,
\nonumber\\
&& \lim_{\tau \to \infty} \sum_{k \in \N}
\frac{4 p^{(\nu)}(\tau, x| (j_{\nu,k})^2)}{ (J_{\nu+1}(x))^2}
\cM_{\xi_{J_{\nu}}^{\langle 2 \rangle}}^{(j_{\nu,k})^2}(t+\tau, R^{(\nu)}(t))
= \cM_{\xi_{J_{\nu}}^{\langle 2 \rangle}}^{x}(t, R^{(\nu)}(t))
\nonumber\\
&& \quad
=\frac{x^{\nu/2}}{(J_{\nu+1}(x))^2}
\frac{1}{(R^{(\nu)}(t))^{\nu/2}}
\int_{0}^{1} d \lambda \,
e^{\lambda t/2} J_{\nu}(\sqrt{\lambda x}) 
J_{\nu} \left(\sqrt{\lambda R^{(\nu)}(t)} \right),
\, x \in \R_+,  
\label{eqn:conv_M}
\end{eqnarray}
for $0 < t \leq T < \infty$.

\clearpage
\SSC{Noncolliding Random Walk 
\label{sec:noncollidingRW}}
\subsection{Construction \label{sec:noncRW_const}}

Consider a random walk (RW),  
$\V(t)=(V_1(t), \dots, V_N(t)), t \in \N_0$ on $\Z^N$, 
such that the components $V_j(t), j=1,2,\dots, N$ are independent
simple and symmetric RWs;
\begin{eqnarray}
V_j(0) &=& u_j \in \Z,
\nonumber\\
V_j(t) &=& u_j+\zeta_j(1)+\zeta_j(2)+ \cdots
+ \zeta_j(t), \quad
\quad t \in \N, \quad 1 \leq j \leq N, 
\label{eqn:V1}
\end{eqnarray}
where $\{\zeta_j(t): 1 \leq j \leq N, t \in \N\}$
is a family of i.i.d. random variables binomially
distributed as
\begin{equation}
\rP[\zeta_j(1)=1]=\frac{1}{2}, \quad
\rP[\zeta_j(1)=-1]=\frac{1}{2},
\quad 1 \leq j \leq N.
\label{eqn:P1}
\end{equation}
For each component, $V_j(\cdot), 1 \leq j \leq N$,
the transition probability is given by
\begin{eqnarray}
&& p(t-s,y|x) 
= \rP[V_j(t)=y | V_j(s)=x]
\nonumber\\
&& \quad = \left \{ \begin{array}{l}
\displaystyle{ 
\frac{1}{2^{t-s}} {t-s \choose [(t-s)+(y-x)]/2 }}, \cr
\qquad \qquad
\mbox{if $t \geq s, \, -(t-s) \leq y-x \leq t-s,  \, [(t-s)+(y-x)]/2 \in \Z$}, \cr
0, \cr
\qquad \qquad \mbox{otherwise}.
\end{array} \right.
\label{eqn:tp1}
\end{eqnarray}
Put $\Z^N=\Z^N_{\rm e} \sqcup \Z^N_{\rm o}$ with
\begin{eqnarray}
\Z^N_{\rm e} &=& \{\x=(x_1, \dots, x_N) :
x_j \in 2 \Z, 1 \leq j \leq N\},
\nonumber\\
\Z^N_{\rm o} &=& \{\x=(x_1, \dots, x_N) :
x_j \in 1+2 \Z, 1 \leq j \leq N\}.
\nonumber
\end{eqnarray}
We always take the initial point 
$\u=(u_1, \dots, u_N) =\V(0)$ from $\Z^N_{\rm e}$,
then $\V(t) \in \Z^N_{\rm e}$, if $t$ is even,
and $\V(t) \in \Z^N_{\rm o}$, if $t$ is odd.
The probability space is denoted as $(\Omega, \cF, \rP_{\u})$.
The expectation is written as $\rE_{\u}$,
which is given by the summation over all walks
$\{\V(t) : t \in \N_0 \}$ 
started at $\u$ with the transition probability (\ref{eqn:tp1})
for each component.

Let
$$
\W_N=\{\x=(x_1, \dots, x_N) \in \R^N: x_1 < \cdots < x_N\}
$$
be the Weyl chamber of type A$_{N-1}$.
Define $\tau_{\u}$ be the exit time from the Weyl chamber
of the RW started at $\u \in \Z_{\rm e}^N \cap \W_N$, 
\begin{equation}
\tau_{\u}=\inf \{ t \geq 1: \V(t) \notin \W_N\}.
\label{eqn:tau1b}
\end{equation}
In the present paper, we study the RW 
{\it conditioned to stay in $\W_N$ forever},
that is, $\tau_{\u}=\infty$ is conditioned. 
We call such a conditional RW 
the {\it (simple and symmetric) noncolliding RW}, since when we regard the
$j$-th component $V_j(\cdot)$ as the position of
$j$-th particle on $\Z, 1 \leq j \leq N$,
if $\tau_{\u} < \infty$, 
then at $t=\tau_{\u}$ there is at least
one pair of particles $(j, j+1)$, which collide with each other;
$V_j(\tau_{\u})=V_{j+1}(\tau_{\u}), 1 \leq j \leq N-1$.
Such a conditional RW is also called 
{\it vicious walkers} in statistical physics \cite{Fis84,CK03},
{\it non-intersecting paths}, {\it non-intersecting walks},
and {\it ordered random walks} (see \cite{EK08}).

Let $\mM$ be the space of nonnegative integer-valued Radon measure 
on $\Z$. We consider the noncolliding RW as a process in $\mM$
and represent it by
\begin{equation}
\Xi(t, \cdot)=\sum_{j=1}^N \delta_{X_j(t)}(\cdot),
\quad t \in \N_0,
\label{eqn:Xi1RW}
\end{equation}
where 
\begin{equation}
\X(t)=(X_1(t), \dots, X_N(t)) \in \Z^N \cap \W_N, 
\quad t \in \N_0.
\label{eqn:Xi2}
\end{equation}
The configuration $\Xi(t, \cdot) \in \mM, t \in \N_0$ is
unlabeled, while $\X(t) \in \Z^N \cap \W_N, t \in \N_0$ is labeled.
We write the probability measure for $\Xi(t, \cdot), t \in \N_0$
started at $\xi \in \mM$ as
$\P_{\xi}$ with expectation $\E_{\xi}$,
and introduce a filtration
$\{\cF(t) : t \in \N_0\}$ defined by
$\cF(t)=\sigma(\Xi(s): 0 \leq s \leq t, s \in \N_0)$.
Then the above definition of the noncolliding RW gives the follows: 
Let $\xi=\sum_{j=1}^N \delta_{u_j}$ 
with $\u \in \Z^N_{\rm e} \cap \W_N$,
and $t \in \N$, $t \leq T \in \N$.
For any $\cF(t)$-measurable bounded function $F$,
\begin{equation}
\E_{\xi} \Big[F(\Xi(\cdot)) \Big]
= \lim_{n \to \infty} \rE_{\u} \left[ \left.
F \left(\sum_{j=1}^N \delta_{V_j(\cdot)} \right)
\right| \tau_{\u} > n \right].
\label{eqn:E1}
\end{equation}
The important fact is that, if we write
the Vandermonde determinant as
\begin{equation}
h(\x)=\det_{1 \leq j, k \leq N} [x_j^{k-1}]
=\prod_{1 \leq j < k \leq N} (x_k-x_j),
\label{eqn:Vand1}
\end{equation}
the expectation (\ref{eqn:E1}) is obtained by
an $h$-transform in the sense of Doob
of the form
\begin{equation}
\E_{\xi}\Big[F(\Xi(\cdot))\Big]
=\rE_{\u} \left[ 
F \left(\sum_{j=1}^N \delta_{V_j(\cdot)} \right)
\1(\tau_{\u} > T)
\frac{h(\V(T))}{h(\u)} \right].
\label{eqn:E2}
\end{equation}
See, for instance, Lemma 4.4 in \cite{Koe05}.

The formula (\ref{eqn:E2}) is a discrete analogue
of the construction of noncolliding Brownian motion (BM)
by Grabiner \cite{Gra99} as an $h$-transform of
absorbing BM in $\W_N$.
The noncolliding BM is equivalent to Dyson's BM model
(with parameter $\beta=2$) and the latter is known
as an eigenvalue process of Hermitian matrix-valued BM
\cite{Dys62,Meh04,Spo87,Gra99,Joh01,KT04,KT11b,Tao12,Osa12,Osa13}.
Then the noncolliding RW has been attracted much attention
as a discretization of models associated with
the Gaussian random matrix ensembles
\cite{Bai00,Joh02,NF02,KT03a,Joh05,BS07,For10,Fei12}.

Nagao and Forrester \cite{NF02} studied
a `bridge' of noncolliding RW
started from $\u_0=(2j)_{j=0}^{N-1}$ at $t=0$
and returned to the same configuration $\u_0$
at time $t=2M, M \in \N_0$.
They showed that at time $t=M$
the spatial configuration provides 
a determinantal point process and the correlation 
kernel is expressed by using the symmetric Hahn polynomials.
Johansson \cite{Joh05} generalized the process to
a bridge from $\u_0$ at $t=0$
to $M_2-M_1+\u_0$ at $t=M_1+M_2$, $M_1, M_2 \in \N_0, M_2 > M_1$,
and proved that the process is determinantal.
The dynamical correlation kernel is of
the Eynard-Mehta type and called 
the extended Hahn kernel.
For the noncolliding RW defined 
for infinite time-period $t \in \N_0$ by
(\ref{eqn:E1}) or (\ref{eqn:E2}) \cite{Koe05,EK08},
however, determinantal structure of 
spatio-temporal correlations has not been clarified so far.

\subsection{Determinantal martingales for noncolliding RW
 \label{sec:noncRW_dm}}

Let $\widetilde{W}_j(t), t \in \N_0, 1 \leq j \leq N$
be independent copies of (\ref{eqn:Wtilde1}). 
Set $\widetilde{\bW}(t)=(\widetilde{W}_1(t), \dots, \widetilde{W}_N(t))$
$\in \R^N, t \in \N_0$ in the probability space
$(\widetilde{\Omega}, \widetilde{\cF}, \widetilde{\rP})$.
We consider a complex process
$\bZ(t)=(Z_1(t), \dots, Z_N(t)), t \in \N_0$ with (\ref{eqn:Z1}).
The probability space for (\ref{eqn:Z1}) is a product
of the probability space $(\Omega, \cF, \rP_{\u})$ 
for the RW, $\V(t), t \in \N_0$, 
and $(\widetilde{\Omega}, \widetilde{\cF}, \widetilde{\rP})$
for $\widetilde{\bW}(t), t \in \N_0$.
Let $\bE_{\u}$ be the expectation for the process $\bZ(t), t \in \N_0$ 
with the initial condition $\bZ(0)=\u \in \Z_{\rm e} \cap \W_N$.

By multilinearity of determinant, the Vandermonde determinant 
(\ref{eqn:Vand1}) does not change
in replacing $x_j^{k-1}$ by any monic polynomial of $x_j$ of 
degree $k-1$, $1 \leq j, k \leq N$.
Note that $m_{k-1}(t,x_j)$ is a monic polynomial of $x_j$ of degree $k-1$. 
Then
\begin{eqnarray}
\frac{h(\V(t))}{h(\u)}
&=& \frac{1}{h(\u)} \det_{1 \leq j, k \leq N} [m_{k-1}(t, V_j(t))]
\nonumber\\
&=& \frac{1}{h(\u)} \det_{1 \leq j, k \leq N}
\Big[ \widetilde{\rE} [Z_j(t)^{k-1}  ] \Big]
\nonumber\\
&=& \widetilde{\rE} \left[ \frac{1}{h(\u)}
\det_{1 \leq j, k \leq N}
[Z_j(t)^{k-1}] \right],
\nonumber
\end{eqnarray}
where we have used the multilinearity of determinant
and independence of $Z_j(t)$'s.
Therefore, we have obtained the equality,
\begin{equation}
\frac{h(\V(t))}{h(\u)}
= \widetilde{\rE} \left[
\frac{h(\bZ(t))}{h(\u)} \right], \quad t \in \N_0.
\label{eqn:Vand2}
\end{equation}

Now we consider the determinant identity \cite{KT13},
\begin{equation}
\frac{h(\z)}{h(\u)}
=\det_{1 \leq j, k \leq N} \Big[
\Phi_{\xi}^{u_k}(z_j) \Big],
\label{eqn:Vand3}
\end{equation}
where $\xi=\sum_{j=1}^N \delta_{u_j}, \u =(u_1, \dots, u_N) \in \W_N$, and
$\Phi_{\xi}^{u_k}(z)$ is given by
\begin{eqnarray}
\Phi_{\xi}^{u_k}(z)
&=& \prod_{\substack{1 \leq j \leq N,  \cr j \not= k}}
\frac{z-u_j}{u_k-u_j},
\quad 1 \leq k \leq N.
\label{eqn:Phi_RW1}
\end{eqnarray}
Let
\begin{equation}
\cM_{\xi}^{u_k}(t, V_j(t))
\equiv \widetilde{\rE}
\Big[ \Phi_{\xi}^{u_k}(Z_j(t)) \Big],
\quad t \in \N_0, \quad1 \leq j \leq N. 
\label{eqn:mart1}
\end{equation}
Since $\Phi_{\xi}^{u_k}(z)$ is a polynomial of $z$
of degree $N-1$, 
$\cM_{\xi}^{u_k}(t, V_j(t))$ is expressed by a linear combination 
of the polynomial martingales $\{m_n(t,V_j(t))\ : 0 \leq n \leq N-1\}$. 
Then $\cM_{\xi}^{u_k}(t, V_j(t)), 1 \leq j \leq N,
t \in \N_0$
are independent martingales
and 
\begin{eqnarray}
\rE_{\u}[\cM_{\xi}^{u_k}(t, V_j(t))]
&=& \rE_{\u}[\cM_{\xi}^{u_k}(0, V_j(0))]
\nonumber\\
&=& \cM_{\xi}^{u_k}(0, u_j)
\nonumber\\
&=& \Phi_{\xi}^{u_k}(u_j)=\delta_{jk},
\quad 1 \leq j, k \leq N.
\label{eqn:cM2}
\end{eqnarray} 
Using the identity (\ref{eqn:Vand3})
for $h(\bZ(t))/h(\u)$ in (\ref{eqn:Vand2}), 
we have
\begin{eqnarray}
\frac{h(\V(t))}{h(\u)}
&=& \widetilde{\rE} \left[
\det_{1 \leq j, k \leq N}
[ \Phi_{\xi}^{u_k}(Z_j(t)) ] \right]
\nonumber\\
&=& \det_{1 \leq j, k \leq N} \Big[
\widetilde{\rE} [ \Phi_{\xi}^{u_k}(Z_j(t)) ] \Big],
\nonumber
\end{eqnarray}
where independence of $Z_j(t)$'s is used.
Let
\begin{equation}
\cD_{\xi}(t,\V(t))
=\det_{1 \leq j, k \leq N}
[\cM_{\xi}^{u_k}(t, V_j(t))],
\quad t \in \N_0.
\label{eqn:mart2}
\end{equation}
We obtain the equality 
\begin{equation}
\frac{h(\V(t))}{h(\u)}
=\cD_{\xi}(t, \V(t)), \quad t \in \N_0. 
\label{eqn:dm2}
\end{equation}

In other words, for the complex processes
(\ref{eqn:Z1}),
$\Phi_{\xi}^{u_k}(Z_j(\cdot)), 1 \leq j, k \leq N$ are 
discrete-time complex martingales
(see Section V.2 of \cite{RY05}) such that,
for any $t \in \N_0$, 
\begin{equation}
\bE_{\u} [\Phi_{\xi}^{u_k}(Z_j(t))]
=\bE_{\u}[\Phi_{\xi}^{u_k}(Z_j(0))] =\delta_{j k},
\quad 1 \leq j, k \leq N.
\label{eqn:conformal1}
\end{equation}
\vskip 0.3cm

\subsection{DMR and CPR for noncolliding RW
 \label{sec:noncRW_DMR}}

Since we consider the noncolliding RW as a process
represented by an unlabeled configuration (\ref{eqn:Xi1}),
measurable functions of $\Xi(\cdot)$ are only
symmetric functions of $N$ variables, 
$X_j(\cdot), 1 \leq j \leq N$.
Then by the equality (\ref{eqn:dm2}), we obtain the
following DMR and CPR
for the present noncolliding RW \cite{Kat13b}.

\begin{thm}
\label{thm:DMR_RW}
Suppose that $N \in \N$ and
$\xi=\sum_{j=1}^{N} \delta_{u_j}$
with $\u=(u_1, \dots, u_N) \in \Z_{\rm e}^N \cap \W_N$.
Let $t  \in \N, t \leq T \in \N$.
For any ${\cal F}(t)$-measurable bounded function $F$
we have
\begin{eqnarray}
\E_{\xi} \left[ F \left(\Xi(\cdot) \right) \right]
&=& \rE_{\u} \left[F \left( \sum_{j=1}^{N} \delta_{V_j(\cdot)} \right)
\cD_{\xi}( T, \V(T)) \right]
\nonumber\\
&=& \bE_{\u} \left[F \left( \sum_{j=1}^{N} \delta_{\Re Z_j(\cdot)} \right)
\det_{1 \leq j, k \leq N}
[\Phi_{\xi}^{u_k}(Z_j(T))] \right]. 
\label{eqn:DMR1}
\end{eqnarray}
\end{thm}
\vskip 0.3cm
\noindent{\it Proof.} \,
It is sufficient 
to consider the case that $F$ is given as
$F(\Xi(\cdot))
= \prod_{m=1}^M g_m(\X(t_m))$
for $M \in \N$, $t_m \in \N, 1 \leq m \leq M$, 
$t_1< \cdots <t_M \leq T \in \N$, 
with symmetric bounded measurable functions $g_m$ on 
$\Z^N$, $1 \leq m \leq M$.
Here we prove the equalities
\begin{eqnarray}
\E_{\xi}\left[ \prod_{m=1}^M g_m(\X(t_m)) \right]
&=& \rE_{\u} \left[ \prod_{m=1}^M g_m(\V(t_m))
\cD_{\xi}( T, \V(T)) \right]
\nonumber\\
&=& \bE_{\u} \left[ \prod_{m=1}^M g_m(\V(t_m))
\det_{1 \leq j, k \leq N}
[\Phi_{\xi}^{u_k}(Z_j(T)) \right].
\label{eqn:prA1}
\end{eqnarray}
By (\ref{eqn:E2}), the LHS of (\ref{eqn:prA1}) is given by
\begin{equation}
\rE_{\u} \left[ 
\prod_{m=1}^M g_m(\V(t_m))
\1(\tau_{\u} > t_M)
\frac{h(\V(t_M))}{h(\u)} \right],
\label{eqn:prA2}
\end{equation}
where we used the fact that $h(\V(\cdot))/h(\u)$ is martingale.
At time $t=\tau_{\u}$, there are at least one pair $(j, j+1)$ 
such that $V_j(\tau_{\u})=V_{j+1}(\tau_{\u}), 1 \leq j \leq N-1$.
We choose the minimal $j$.
Let $\sigma_{j, j+1}$ be the permutation of
the indices $j$ and $j+1$ and for
$\v=(v_1, \dots, v_N) \in \Z^N$ we put
$\sigma_{j,j+1}(\v)=(v_{\sigma_{j,j+1}(k)})_{k=1}^N
=(v_1, \dots, v_{j+1}, v_j, \dots, v_N)$.
Let $\u'$ be the labeled configuration of the process
at time $t=\tau_{\u}$.
Since $u'_j=u'_{j+1}$ by the above setting,
under the probability law $\rP_{\u'}$
the processes $\V(t), t > \tau_{\u}$ and 
$\sigma_{j,j+1}(\V(t)), t > \tau_{\u}$
are identical in distribution.
Since $g_m, 1 \leq m \leq M$ are symmetric,
but $h$ is antisymmetric, the Markov property
of the process $\V(\cdot)$ gives
$$
\rE_{\u} 
\left[ 
\prod_{m=1}^M g_m(\V(t_m))
\1(\tau_{\u} \leq t_M)
\frac{h(\V(t_M))}{h(\u)} \right]=0.
$$
Therefore, (\ref{eqn:prA2}) is equal to
$$
\rE_{\u} \left[ 
\prod_{m=1}^M g_m(\V(t_m))
\frac{h(\V(t_M))}{h(\u)} \right].
$$
By the equality (\ref{eqn:dm2}) and the
martingale property of $\cD_{\xi}(\cdot, \V(\cdot))$, 
we obtain the first line of (\ref{eqn:prA1}). 
By definitions of $\bE_{\u}$ and $\cD_{\xi}$, the second line is valid.
Then the proof is completed.
\qed
\vskip 0.5cm

Note that the CPR in the second line of (\ref{eqn:DMR1})
may correspond to the complex BM representation
reported in \cite{KT13} for the noncolliding BM.

Then by Proposition \ref{thm:DM_det}, we will immediately
conclude the following \cite{Kat13b}.
Let
\begin{equation}
\mbK_{\xi}(s,x;t,y)
= \left\{ \begin{array}{l}
\displaystyle{
\sum_{j=1}^N p(s, x|u_j) \cM_{\xi}^{u_j}(t,y)
- \1(s>t) p(s-t,x|y),
}
\cr
\qquad \qquad \qquad \mbox{if} \quad
(s,x), (t,y) \in \N_0 \times \Z, \quad s+x, t+y \in 2 \Z,
\cr
0,
\qquad \qquad \quad \mbox{otherwise},
\end{array} \right.
\label{eqn:K_RW1}
\end{equation}
where $p$ is the transition probability for RW (\ref{eqn:tp1}).
Here 
\begin{equation}
\cM_{\xi}^{u_j}(t,y)
= \widetilde{\rE}[
\Phi_{\xi}^{u_j}(y+i \widetilde{W}(t))]
\label{eqn:mart1b}
\end{equation}
is a functional of initial configuration
$\xi=\sum_{j=1}^N \delta_{u_j}$ through (\ref{eqn:Phi_RW1}).

\begin{cor} 
\label{thm:main_RW}
For any initial configuration $\xi \in \mM$
with $\xi(\Z^N_{\rm e})=N \in \N$, 
the noncolliding RW, $(\Xi(t), t \in \N_0, \P_{\xi})$
is determinantal with the kernel
(\ref{eqn:K_RW1}) with (\ref{eqn:mart1b}) 
in the sense that
the moment generating function (\ref{eqn:GF1})
is given by Fredholm determinant
\begin{equation}
\Psi_{\xi}^{\t}[\f]
=\mathop{{\rm Det}}_
{\substack{
(s,t)\in \{t_1,t_2,\dots, t_M\}^2, \\
(x,y)\in \Z^2}
}
 \Big[\delta_{st} \delta_x(y)
+ \mbK_{\xi}(s,x;t,y) \chi_{t}(y) \Big],
\label{eqn:Fred}
\end{equation}
and then all spatio-temporal correlation functions
are given by determinants as
\begin{eqnarray}
&& \rho_{\xi} \Big(t_1,\x^{(1)}_{N_1}; \dots;t_M,\x^{(M)}_{N_M} \Big) 
=
\left\{ \begin{array}{l}
\displaystyle{
\det_{\substack
{1 \leq j \leq N_{m}, 1 \leq k \leq N_{n}, \\
1 \leq m, n \leq M}
}
\Bigg[
\mbK_{\xi}(t_m, x_{j}^{(m)}; t_n, x_{k}^{(n)} )
\Bigg],
}
\cr
\qquad \mbox{if} \quad
\x^{(m)}_{N_m} \in \Z_{\rm e}^{N_m} \cap \W_{N_m}, 
t_m=\mbox{even}, \cr
\qquad  \mbox{or} \quad
\x^{(m)}_{N_m} \in \Z_{\rm o}^{N_m} \cap \W_{N_m}, 
t_m=\mbox{odd},  \quad 1 \leq m \leq M, \cr
0, \quad \mbox{otherwise},
\end{array} \right.
\nonumber\\
\label{eqn:main1}
\end{eqnarray}
$t_m \in \N, 1 \leq m \leq M$, 
$t_1 < \cdots < t_M$, and
$0 \leq N_m \leq N, 1 \leq m \leq M$.
\end{cor}

\clearpage
\SSC{DMR in O'Connell Process 
\label{sec:DMR_OP}}
\subsection{Quantum Toda lattice and Whittaker function \label{sec:QToda}}

Let $a >0$.
The Hamiltonian of the 
GL($N, \R$)-quantum Toda lattice is given by
\begin{equation}
\cH_N^a=-\frac{1}{2} \Delta+ \frac{1}{a^2} V_N(\x/a),
\quad \x=(x_1, x_2, \dots, x_N) \in \R^N
\label{eqn:TodaH1}
\end{equation}
with the Laplacian 
$\Delta=\sum_{j=1}^N \partial^2/\partial x_j^2$
and the potential
\begin{equation}
V_N(\x)=\sum_{j=1}^{N-1} e^{-(x_{j+1}-x_j)}.
\label{eqn:VN1}
\end{equation}
For $\vnu=(\nu_1, \nu_2, \dots, \nu_N) \in \W_N^{\rm A}$,
the eigenfunction problem
\begin{equation}
\cH_N^a \psi^{(N)}_{\vnu}(\x)=\lambda(\vnu) \psi^{(N)}_{\vnu}(\x)
\label{eqn:Toda_ev1}
\end{equation}
for the eigenvalue 
\begin{equation}
\lambda(\vnu)=-|\vnu|^2/2
\label{eqn:Toda_ev2}
\end{equation}
is uniquely solved
under the condition that
\begin{equation}
e^{-\vnu \cdot \x} \psi^{(N)}_{\vnu}(\x)
\quad \mbox{is bounded}
\label{eqn:Toda_ev3}
\end{equation}
and 
\begin{equation}
\lim_{\x \to \infty, \x \in \W_N}
e^{-\vnu \cdot \x} \psi^{(N)}_{\vnu}(\x)
=\prod_{1 \leq j < k \leq N} \Gamma(\nu_k-\nu_j),
\label{eqn:Toda_ev4}
\end{equation}
where $\x \to \infty, \x \in \W_N$ means 
$x_{j+1}-x_j \to \infty, 1 \leq j \leq N-1$, 
and $\Gamma(\cdot)$ denotes the Gamma function.
The eigenfunction $\psi^{(N)}_{\vnu}(\, \cdot \,)$
is called the class-one Whittaker function \cite{BO11,OCo12b}.

The class-one Whittaker function $\psi^{(N)}_{\vnu}(\x)$
has several integral representations, one of which 
was given by Givental \cite{Giv97},
\begin{equation}
\psi_{\vnu}^{(N)}(\x)
= \int_{\mbT_N(\x)} 
\exp\left(\cF_{\vnu}^{(N)}(\T) \right) d \T.
\label{eqn:Givental1}
\end{equation}
Here  the integral is performed 
over the space $\mbT_N(\x)$ of all real lower
triangular arrays with size $N$,
$\T=(T_{j,k}, 1 \leq k \leq j \leq N)$, 
with $T_{N,k}=x_k, 1 \leq k \leq N$,
and 
\begin{equation}
\cF_{\vnu}^{(N)}(\T)
=
\sum_{j=1}^{N} \nu_{j}
\left( \sum_{k=1}^{j} T_{j, k}
-\sum_{k=1}^{j-1} T_{j-1, k} \right)
- \sum_{j=1}^{N-1} \sum_{k=1}^j
\Big\{ e^{-(T_{j,k}-T_{j+1,k})}
+e^{-(T_{j+1, k+1}-T_{j,k})} \Big\}.
\label{eqn:Givental2}
\end{equation}
We can prove that \cite{OCo12a,COSZ11}
\begin{equation}
\lim_{a \to 0} a^{N(N-1)/2}
\psi^{(N)}_{a \vnu}(\x/a)
=\frac{\displaystyle{
\det_{1 \leq j, \ell \leq N}
[e^{x_j \nu_{\ell}}]}}{h(\vnu)}.
\label{eqn:Toda4}
\end{equation}
where $h(\vnu)$ is the Vandermonde determinant
(\ref{eqn:Vand}).

The following orthogonality relation is proved
for the class-one Whittaker functions
\cite{STS94,Wal92},
\begin{equation}
\int_{\R^N} \psi^{(N)}_{-i \k}(\x)
\psi^{(N)}_{i \k'}(\x) d \x
=\frac{1}{s_N(\k) N!}
\sum_{\sigma \in \cS_N} \delta(\k-\sigma(\k')),
\label{eqn:orth1}
\end{equation}
for $\k, \k' \in \R^N$,
where $s_N(\cdot)$ is the density function of
the Sklyanin measure \cite{Skl85}
\begin{eqnarray}
s_N(\vmu)
&=& \frac{1}{(2 \pi)^N N!}
\prod_{1 \leq j < \ell \leq N}
|\Gamma(i(\mu_{\ell}-\mu_j))|^{-2}
\nonumber\\
&=& \frac{1}{(2 \pi)^N N!}
\prod_{1 \leq j < \ell \leq N}
\left\{ (\mu_{\ell}-\mu_j)
\frac{\sinh \pi (\mu_{\ell}-\mu_j)}{\pi} \right\},
\quad \vmu \in \R^N.
\label{eqn:sN1}
\end{eqnarray}
Borodin and Corwin proved that for a class of test functions,
the orthogonality relation (\ref{eqn:orth1}) can be
extended for any $\k, \k' \in \C^N$ \cite{BC11}.
Moreover, the following recurrence relations
with respect to $\vnu$ 
are established \cite{KL01,BC11};
for $1 \leq r \leq N-1, \vnu \in \C^N$,
\begin{equation}
\sum_{\substack{I \subset \{1, \dots, N\},\\
|I|=r}}
\prod_{\substack{j \in I, \\ 
k \in \{1, 2, \dots, N\} \setminus I}}
\frac{1}{i(\nu_k-\nu_j)}
\psi^{(N)}_{i(\vnu+i \e_I)}(\x)
= \exp \left( - \sum_{j=1}^r x_j \right)
\psi^{(N)}_{i \vnu}(\x),
\label{eqn:iden1}
\end{equation}
where $\e_I$ is the vector with ones in the slots
of label $I$ and zeros otherwise;
$$
(\e_I)_j= \left\{ \begin{array}{ll}
1, \quad & j \in I, \cr
0, \quad & j \in \{1, \dots, N\} \setminus I.
\end{array} \right.
$$
In particular, for $r=1$,
\begin{equation}
\sum_{j=1}^{N} \prod_{1 \leq k \leq N: k \not=j}
\frac{1}{i(\nu_k-\nu_j)}
\psi^{(N)}_{i(\vnu+i \e_{\{j\}})}(\x)
=e^{-x_1} \psi^{(N)}_{i \vnu}(\x),
\label{eqn:iden2}
\end{equation}
where the $\ell$-th component of the 
vector $\e_{\{j\}}$ is
$(\e_{\{j\}})_{\ell}=\delta_{j \ell}, 1 \leq j, \ell \leq N$.
As fully discussed by Borodin and Corwin \cite{BC11},
the recurrence relations (\ref{eqn:iden1}) are
derived as the $q \to 1$ limit of
the eigenfunction equations associated to
the Macdonald difference operators
in the theory of symmetric functions \cite{Mac99}.
For more details on Whittaker functions, 
see \cite{KL01,BO11,Kat11,OCo12a,Kat12a,BC11} and references therein.

\subsection{O'Connell process \label{sec:OP}}

O'Connell introduced an $N$-component diffusion process,
$N \geq 2$, which can be regarded as a stochastic version of 
a quantum open Toda-lattice \cite{OCo12a}.
Let $a > 0$.
The infinitesimal generator of the O'Connell process is given by
\begin{eqnarray}
\cL^{\nu,a}_N &=&
- (\psi^{(N)}_{\vnu}(\x/a))^{-1} 
\left(\cH_N^a +\frac{1}{2}|\vnu|^2 \right) \psi^{(N)}_{\vnu}(\x/a)
\nonumber\\
&=& \frac{1}{2} \Delta
+\nabla \log \psi^{(N)}_{\vnu}(\x/a) \cdot \nabla,
\label{eqn:generator}
\end{eqnarray}
where $\nabla=(\partial/\partial x_1, \dots, \partial/\partial x_N)$.
This multivariate diffusion process
is an extension of a one-dimensional diffusion studied by Matsumoto and Yor 
\cite{MY00,MY05}. 
(The Matsumoto-Yor process describes 
time-evolution of the relative coordinate of the $N=2$ case.)

We can show that the O'Connell process is realized
as the following mutually killing BMs 
{\it conditioned that all particles survive forever},
if the particle position of the $j$-th BM
is identified with the $j$-th component of the O'Connell process,
$1 \leq j \leq N$ \cite{Kat11,Kat12a,Kat12b}.
Let $B_j(t), 1 \leq j \leq N$ be independent 
one-dimensional standard BMs started at
$B_j(0)=x_j \in \R$,
and for $\vnu=(\nu_1, \nu_2, \dots, \nu_N) \in \R^N$,
\begin{equation}
B_j^{\nu_j}(t)=B_j(t)+ \nu_j t, \qquad 1 \leq j \leq N
\label{eqn:driftBM}
\end{equation}
be drifted BMs.
We consider an $N$-particle system of
BMs with drift vector $\vnu$,
$\B^{\nu}(t)=(B_1^{\nu_1}(t), \dots, B_N^{\nu_N}(t)), t \geq 0$,
such that the probability $P_N^a(t|\{\B^{\nu}(s)\}_{0 \leq s \leq t})$
that all $N$ particles survive
up to time $t$ conditioned on a path $\{\B^{\nu}(s)\}_{0 \leq s \leq t}$
decays following the equation
\begin{equation}
\frac{d}{dt}
P_N^a(t|\{\B^{\nu}(s)\}_{0 \leq s \leq t})
=- \frac{1}{a^2} V_N(\B^{\nu}(t)/a) 
P_N^a(t|\{\B^{\nu}(s)\}_{0 \leq s \leq t}),
\quad t \geq 0.
\label{eqn:dPdt}
\end{equation}
It is a system of mutually killing BMs, in which the Toda-lattice
potential (\ref{eqn:VN1}) determines the
decay rate of the survival probability depending on
a configuration $\B^{\nu}(t)$ \cite{Kat12a}.
With the initial condition $\B^{\nu}(0)=\x \in \W_N$,
the survival probability
is obtained by averaging over all paths of BMs started at $\x$ as
\begin{equation}
P_N^a(t; \x, \vnu) = \rE_{\x}\Big[ P_N^a(t|\{\B^{\nu}(s)\}_{0 \leq s \leq t}) \Big], 
\label{eqn:surP_OP}
\end{equation}
and
we can show that \cite{Kat12a,OCo12b,Kat12b}
\begin{eqnarray}
\lim_{t \to \infty} 
P_N^a(t; \x, \vnu) &=& c^a_1(N, \vnu) e^{-\vnu \cdot \x/a}
\psi^{(N)}_{\vnu}(\x/a),
\quad \mbox{if $\vnu \in \W_N$, $\vnu \not=0$},
\nonumber\\
P_N^{a}(t; \x, 0) &\sim& c^a_2(N) t^{-N(N-1)/4} \psi^{(N)}_0(\x/a)
\quad \mbox{as $t \to \infty$},
\label{eqn:PNasym}
\end{eqnarray}
where $c^a_1(N,\vnu)$ and $c^a_2(N)$ are independent of
$\x$ and $t$.
Then, conditionally on surviving of all $N$ particles,
the equivalence of this vicious BM,
which has a killing term given by the Toda-lattice potential, 
with the O'Connell process is proved.
We note that the parameter $a>0$ in the killing rate (\ref{eqn:dPdt})
with (\ref{eqn:VN1}) indicates the characteristic range
of interaction to kill neighboring particles
as well as the characteristic length in which neighboring particles
can exchange their order in $\R$.
It implies that if we take the limit $a \to 0$,
the O'Connell process is reduced to the noncolliding BM.
(The original vicious Brownian motion is a system of BMs such that
if pair of particles collide they are annihilated immediately.
The noncolliding BM is the vicious BM conditioned never to
collide with each other, and thus all particles
survive forever.)

The transition probability density 
for the O'Connell process with $\vnu$ is given by
\cite{Kat12a}
\begin{equation}
P_N^{\nu, \, a}(t, \y|\x)
=e^{-t|\vnu|^2/2 a^2}
\frac{\psi^{(N)}_{\vnu}(\y/a)}{\psi^{(N)}_{\vnu}(\x/a)}
Q_N^{a}(t, \y|\x),
\quad \x, \y \in \R^N, t \geq 0,
\label{eqn:tpdA1}
\end{equation}
with
\begin{equation}
Q_N^{a}(t, \y|\x)
=\int_{\R^N} e^{-t|\k|^2/2}
\psi_{i a \k}^{(N)}(\x/a) \psi_{-i a \k}^{(N)}(\y/a)
s_N(a \k) d \k.
\label{eqn:QN1}
\end{equation}
(See also Proof of Proposition 4.1.32 in \cite{BC11}.)
As a matter of fact, we can confirm that 
$u(t, \x) \equiv P_N^{\nu, \, a}(t, \y|\x)$ satisfies
the Kolmogorov backward equation
associated with the infinitesimal generator
$\cL_N^{\nu,a}$ given by (\ref{eqn:generator}),
\begin{eqnarray}
\frac{\partial u(t, \x)}{\partial t}
&=& \cL_N^{\nu,a} u(t, \x)
\nonumber\\
&=& \frac{1}{2} \sum_{j=1}^{N}
\frac{\partial^2 u(t, \x)}{\partial x_j^2}
+\sum_{j=1}^N
\frac{\partial \log \psi_{\vnu}^{(N)}(\x/a)}{\partial x_j}
\frac{\partial u(t, \x)}{\partial x_j}, 
\label{eqn:OConnell1}
\end{eqnarray}
$\x \in \R^N, t \geq 0$, under the condition
$u(0, \x)=\delta(\x-\y) 
\equiv \prod_{j=1}^N \delta(x_j-y_j), \y \in \R^N$.
We denote the O'Connell process by
\begin{equation}
\X^a(t)=(X^a_1(t), X^a_2(t), \dots, X^a_N(t)), \quad t \geq 0.
\label{eqn:Xa1}
\end{equation}
It is defined as an $N$-particle diffusion process
in $\R$ such that its backward Kolmogorov equation
is given by (\ref{eqn:OConnell1}).
Therefore, (\ref{eqn:Xa1}) is a unique solution of the
following stochastic differential equation
for given initial configuration $\X^a(0)=\x \in \R^N$,
\begin{equation}
d X^a_j(t)= dB_j(t)+ \Big[\F_{N}^{\nu, a}(\X^a(t))\Big]_j dt, 
\quad 1 \leq j \leq N, t \geq 0
\label{eqn:SDE1}
\end{equation}
with
\begin{equation}
\F_{N}^{\nu, a}(\x)=\nabla \log \psi^{(N)}_{\vnu}(\x/a),
\label{eqn:SDE2}
\end{equation}
where $\{B_j(t)\}_{j=1}^N$ are
independent one-dimensional standard BMs
and $[\V]_j$ denotes the $j$-th coordinate of a vector $\V$.

\subsection{Special entrance law \label{sec:entrance}}

Let $N \in \N$, and define
\begin{equation}
\vrho =
\left( -\frac{N-1}{2}, -\frac{N-1}{2}+1, \dots,
\frac{N-1}{2}-1, \frac{N-1}{2} \right).
\label{eqn:rho_OP}
\end{equation}
O'Connell considered the process starting from
$\x=-M \vrho$ and let $M \to \infty$ \cite{OCo12a}.
It was claimed in \cite{OCo12a} (see also \cite{BO11}) that
\begin{equation}
\psi^{(N)}_{\vnu}(-M \vrho)
\sim C e^{-N(N-1) M/8}
\exp \Big( e^{M/2} \cF_{\0}(\T^0) \Big)
\label{eqn:asym1}
\end{equation}
as $M \to \infty$, where the coefficient $C$ 
and the critical point $\T^0$ are
independent of $\vnu$.
Then as a limit of (\ref{eqn:tpdA1}) with (\ref{eqn:QN1}),
we have a probability density function
\begin{eqnarray}
\cP_N^{\nu,a}(t, \x)
&\equiv&
\lim_{M \to \infty} P_N^{\nu, a}(t, \x| -M\vrho)
\nonumber\\
&=& e^{-t|\vnu|^2/2a^2}
\psi^{(N)}_{\vnu}(\x/a) \vartheta_N^{a}(t,\x)
\label{eqn:tpdB1}
\end{eqnarray}
with
\begin{equation}
\vartheta_N^a(t,\x)
=\int_{\R^N} e^{-t|\k|^2/2}
\psi^{(N)}_{-i a \k}(\x/a)
s_N(a \k) d \k
\label{eqn:theta1}
\end{equation}
for any $t > 0$.
Since we have taken the limit $M \to \infty$
for the state $-M \vrho$, we cannot speak of
initial configurations any longer, but for
an arbitrary series of increasing times,
$0 < t_1 < t_2 < \dots < t_M < \infty$, 
the probability density function of the multi-time
joint distributions is given by
\begin{equation}
\cP_N^{\nu, a}(t_1, \x^{(1)}; t_2, \x^{(2)};
\dots; t_M, \x^{(M)}) = 
\prod_{m=1}^{M-1} P^{\nu, a}_N(t_{m+1}-t_m, \x^{(m+1)}|\x^{(m)})
\cP_N^{\nu,a}(t_1, \x^{(1)})
\label{eqn:multiA1}
\end{equation} 
for $\x^{(m)} \in \R^N, 1 \leq m \leq M$.
We can call the probability measure
$\cP_N^{\nu,a}(t, \x) d \x$ with (\ref{eqn:tpdB1})
and $d \x=\prod_{j=1}^N dx_j$
an {\it entrance law coming from} ``$-\infty \vrho$" \cite{OCo12a}
(see, for instance, Section XII.4 of \cite{RY05} for entrance laws).
We note that, by (\ref{eqn:tpdA1}), (\ref{eqn:multiA1}) is
written as
\begin{eqnarray}
&& \cP_N^{\nu, a}(t_1, \x^{(1)}; t_2, \x^{(2)};
\dots; t_M, \x^{(M)})
\nonumber\\
&& \quad = 
e^{-t_M |\vnu|^2/2a^2}
\psi^{(N)}_{\vnu}(\x^{(M)}/a)
\prod_{m=1}^{M-1} Q^{a}_N(t_{m+1}-t_m, \x^{(m+1)}|\x^{(m)})
\vartheta^a_N(t_1, \x^{(1)}).
\nonumber
\end{eqnarray}

The expectation with respect to the distribution 
of the present process started according to the
special entrance law (\ref{eqn:tpdB1}) is denoted by
$\E^{\nu, a}$.
For measurable functions 
$f^{(m)}, 1 \leq m \leq M$,
\begin{eqnarray}
&& \E^{\nu, a} \left[ \prod_{m=1}^{M} f^{(m)}(\X^a(t_m)) \right]
\nonumber\\
&& =e^{-t_M|\vnu|^2/2a^2}
\left\{\prod_{m=1}^M \int_{\R^N} d\x^{(m)} \right\}
f^{(M)}(\x^{(M)}) \psi^{(N)}_{\vnu}(\x^{(M)}/a)
Q^a_N(t_M-t_{M-1}, \x^{(M)}|\x^{(M-1)})
\nonumber\\
&& \qquad \times \prod_{m=2}^{M-1} f^{(m)}(\x^{(m)}) 
Q^a_N(t_m-t_{m-1}, \x^{(m)}|\x^{(m-1)})
f^{(1)}(\x^{(1)}) \vartheta^a_N(t_1, \x^{(1)}),
\label{eqn:Exp2}
\end{eqnarray}
$0 < t_1 < \dots < t_M < \infty$,
where $d\x^{(m)}=\prod_{j=1}^N dx^{(m)}_j, 1 \leq m \leq M$.

The present special entrance law (\ref{eqn:tpdB1}) 
is called a Whittaker measure 
by Borodin and Corwin \cite{BC11}
and denoted by ${\bf WM}_{(\nu; t)}(\x)$.
Note that in the notation of \cite{BC11},
a Whittaker process is a `triangular array extension'
of the Whittaker measure and is not the same as
the O'Connell process.

When $M=1$, for $t > 0$, (\ref{eqn:Exp2}) gives
\begin{eqnarray}
&& \E^{\nu, a}[f(\X^a(t))]
\nonumber\\
&& \quad = e^{-t|\vnu|^2/2a^2}
\int_{\R^N} d \x f(\x) \psi^{(N)}_{\vnu}(\x/a)
\vartheta_N^a(t, \x)
\nonumber\\
&& \quad = e^{-t|\vnu|^2/2a^2}
\int_{\R^N} d \x f(\x) \psi^{(N)}_{\vnu}(\x/a)
\int_{\R^N} d \k e^{-t|\k|^2/2}
\psi^{(N)}_{-i a \k}(\x/a) s_N(a\k).
\label{eqn:SA1}
\end{eqnarray}

\subsection{Combinatorial limit $a \to 0$ \label{sec:a_0}}

The transition probability density
of the absorbing BM in $\W_N^{\rm A}$
is given by the Karlin-McGregor determinant
of (\ref{eqn:p_BM}),
\begin{equation}
q_N(t,\y|\x)=
\det_{1 \leq j, k \leq N} [p(t, y_j|x_k)],
\quad \x, \y \in \W_N^{\rm A}, t \geq 0.
\label{eqn:qN1}
\end{equation}

Consider the drift transform of (\ref{eqn:qN1}),
$$
q_N^{\nu}(t, \y|\x)
=\exp \left\{ -\frac{t}{2}|\vnu|^2
+\vnu \cdot (\y-\x) \right\}
q_N(t, \y|\x).
$$
Then, if 
$\vnu \in \overline{\W}_N^{\rm A}
=\{ \x \in \R^N : x_1 \leq x_2 \leq \cdots \leq x_N \}$,
the transition probability density
of the noncolliding BM
with drift $\vnu$ is given by
\cite{BBO05}
\begin{equation}
p_N^{\nu}(t, \y|\x)
=e^{-t|\vnu|^2/2}
\frac{\displaystyle{\det_{1 \leq j, k \leq N}[e^{\nu_j y_k}]}}
{\displaystyle{\det_{1 \leq j, k \leq N}[e^{\nu_j x_k}]}}
q_N(t, \y|\x),
\quad \x, \y \in \W_N, \quad t \geq 0.
\label{eqn:pNnu1}
\end{equation}
In the limit $\nu_j \to 0, 1 \leq j \leq N$, (\ref{eqn:pNnu1})
becomes
\begin{equation}
p_N(t, \y|\x)=\frac{h(\y)}{h(\x)} q_N(t, \y|\x),
\quad \x, \y \in \W_N^{\rm A}, t \geq 0.
\label{eqn:pN1}
\end{equation}

We prove the following.
(The superscript $a \nu$ is used for the processes
with drift vector $a \vnu=(a \nu_1, \dots, a \nu_N)$.)
\begin{lem}
\label{thm:a0lim}
For $\vnu \in \overline{\W}_N$,
\begin{eqnarray}
\lim_{a \to 0} \cP_N^{a \nu, a}(t, \x) d\x
&=& p_N(t^{-1}, \x/t |\vnu) d (\x/t)
\nonumber\\
&=& p_N^{\nu}(t, \x|\0) d \x,
\quad t > 0.
\label{eqn:a0lim}
\end{eqnarray}
\end{lem}
\vskip 0.3cm

\noindent{\it Proof} \,
By the asymptotics (\ref{eqn:Toda4}) and the
definition (\ref{eqn:qN1}) of $q_N$, we have
\begin{equation}
\lim_{a \to 0} a^{N(N-1)/2} e^{-t|\vnu|^2/2} 
\psi^{(N)}_{a \vnu}(\x/a)
=\left( \frac{2 \pi}{t} \right)^{N/2} e^{|\x|^2/2t}
\frac{q_N(t^{-1}, \x/t|\vnu)}{h(\vnu)}.
\label{eqn:EqA1}
\end{equation}

For $\vartheta_N^a$ defined by the integral (\ref{eqn:theta1}),
we can show that the Whittaker function with purely imaginary index
multiplied by the Sklyanin density,
$\psi^{(N)}_{-i a \k}(\, \cdot \,) s_N(a \k)$,
is uniformly integrable in $a >0$
with respect to the Gaussian measure
$e^{-t |\k|^2/2} d \k, t > 0$.
Then the integral and the limit $a \to 0$ is
interchangeable.
Since 
$$
\psi^{(N)}_{-i a \k}(\x/a)
\sim (-i a)^{-N(N-1)/2}
\frac{\displaystyle{\det_{1 \leq j, \ell \leq N}
[e^{-ix_j k_{\ell}}]}}{h(\k)}, \quad
\mbox{as $a \to 0$}
$$
by (\ref{eqn:Toda4}), and (\ref{eqn:sN1}) gives
$s_N(a \k) \sim a^{N(N-1)} 
(h(\k))^2/\{(2 \pi)^N N!\}$, as $a \to 0$, 
we have
\begin{eqnarray}
&& \lim_{a \to 0} a^{-N(N-1)/2} \vartheta^a_N(t, \x)
\nonumber\\
&& \quad = 
\frac{1}{(2\pi)^N N!} \int_{\R^N} d \k
e^{-t|\k|^2/2} \det_{1 \leq j, \ell \leq N}
[e^{-i x_j k_{\ell}}] h(i \k)
\nonumber\\
&& \quad = \frac{t^{-N(N+1)/4}}{(2 \pi)^{N/2}}
e^{-|\x|^2/2t}
\frac{1}{N!} \int_{\R^N} d (\sqrt{t} \k)
\det_{1 \leq j, \ell \leq N}
\left[ \frac{
e^{-(\sqrt{t} k_{\ell}+i x_j/\sqrt{t})^2/2}
}{\sqrt{2\pi}}
\prod_{m=1}^{\ell-1} (i \sqrt{t} k_{\ell}
-i \sqrt{t} k_m) \right].
\nonumber
\end{eqnarray}
By multi-linearity of determinant,
\begin{eqnarray}
&& 
\frac{1}{N!} \int_{\R^N} d (\sqrt{t} \k)
\det_{1 \leq j, \ell \leq N}
\left[ \frac{
e^{-(\sqrt{t} k_{\ell}+i x_j/\sqrt{t})^2/2}
}{\sqrt{2\pi}}
\prod_{m=1}^{\ell-1} (i \sqrt{t} k_{\ell}
-i \sqrt{t} k_m) \right]
\nonumber\\
&& \quad =
\det_{1 \leq j, \ell \leq N} \left[
\int_{\R} d(\sqrt{t} k)
\frac{e^{-(\sqrt{t} k+i x_j/\sqrt{t})^2/2}}
{\sqrt{2 \pi}} \prod_{m=1}^{\ell-1}
(i \sqrt{t} k-i \sqrt{t} k_m) \right]
\nonumber\\
&& \quad =
\det_{1 \leq j, \ell \leq N} \left[
\int_{\R} d u
\frac{e^{-(u+i x_j/\sqrt{t})^2/2}}
{\sqrt{2 \pi}} \prod_{m=1}^{\ell-1}
(i u-i \sqrt{t} k_m) \right].
\label{eqn:EqA3}
\end{eqnarray}
The integral in the determinant (\ref{eqn:EqA3})
can be identified with an integral representation
given by Bleher and Kuijlaars \cite{BK05,KT10}
for the multiple Hermite polynomial of type II,
$$
P_{\xi_{\ell-1}}(x_j/\sqrt{t})
\quad \mbox{with} \quad
\xi_{\ell-1}(\cdot)=\sum_{m=1}^{\ell-1} \delta_{i \sqrt{t} k_m}(\cdot).
$$
(We set $\xi_0(\cdot) \equiv 0$ and
$\prod_{m=1}^0 (\cdot) \equiv 1$.)
It is a monic polynomial of $x_j/\sqrt{t}$ with 
degree $\ell-1$.
Then (\ref{eqn:EqA3}) is equal to the Vandermonde determinant
$$
h(\x/\sqrt{t})=t^{N(N-1)/4} h(\x/t).
$$
Therefore, we obtain
\begin{equation}
\lim_{a \to 0} a^{-N(N-1)/2} 
\vartheta^a_N(t, \x)
=\frac{1}{(2 \pi t)^{N/2}}
e^{-|\x|^2/2t} h(\x/t).
\label{eqn:EqA5}
\end{equation}
Combining (\ref{eqn:EqA1}) and (\ref{eqn:EqA5}), 
we obtain the equality
\begin{equation}
\lim_{a \to 0} \cP_N^{a \nu, a}(t, \x)
=\frac{h(\x/t)}{h(\vnu)}
q_N(t^{-1}, \x/t | \vnu) t^{-N},
\label{eqn:EqA6}
\end{equation}
which gives the first equality of 
(\ref{eqn:a0lim}) by the formula (\ref{eqn:pN1}).
The second equality is concluded by the
reciprocal relation proved as Theorem 2.1
in \cite{Kat12b}
(see (\ref{eqn:reciprocal2}) below).
The proof is then completed. \qed
\vskip 0.3cm
Moreover, if we take the limit $\vnu \to 0$ in (\ref{eqn:a0lim}),
we have the following
\begin{eqnarray}
\lim_{\vnu \to 0} \lim_{a \to 0}
\cP_N^{a \nu, a}(t, \x)
&=& p_N(t, \x|\0)
\nonumber\\
&=& \frac{t^{-N^2/2}}{(2 \pi)^{N/2}
\prod_{j=1}^N \Gamma(j)}
e^{-|\x|^2/2t} (h_N(\x))^2.
\label{eqn:GUE}
\end{eqnarray}
This is the probability density of 
the eigenvalue distribution of the 
Gaussian unitary ensemble (GUE)
with variance $\sigma^2=t$
of random matrix theory.
It implies that {\it a geometric lifting} of the
GUE-eigenvalue distribution is 
the $\vnu \to 0$ limit of the
entrance law coming from ``$-\infty \vrho$",
\begin{eqnarray}
\cP_N^{a}(t, \x)
&\equiv& \lim_{\vnu \to 0} \cP_N^{\nu, a}(t, \x)
\nonumber\\
&=& \psi^{(N)}_0(\x/a) \vartheta_N^a(t, \x)
\nonumber\\
&=& \psi^{(N)}_0(\x/a) \int_{\R^N} e^{-t|\k|^2/2}
\psi^{(N)}_{-i a \k}(\x/a) s_N(a \k) d \k.
\label{eqn:geoGUE}
\end{eqnarray}

\subsection{Determinantal formula of Borodin and Corwin \label{sec:BC}}

For $x \in \R, a >0$, set
\begin{equation}
\Theta^a(x)=\exp(-e^{-x/a}).
\label{eqn:Theta}
\end{equation}
Note that $\displaystyle{\lim_{a \to 0} \Theta^a(x)=\1(x>0)}$,
that is, (\ref{eqn:Theta}) is 
a softening of an indicator function $\1(x >0)$.

Let $\widetilde{\delta}=\sup \{|\nu_j|: 1 \leq j \leq N\}$
and choose $0 < \delta < 1$ so that $\widetilde{\delta} < \delta/2$.
Borodin and Corwin \cite{BC11} proved that
$\E^{\nu, a}[\Theta^a(X^a_1(t)-h)], 
h \in \R$ is given by a Fredholm determinant
of a kernel $K_{e^{h/a}}$
for the contour integrals on $C((-1) \circ \nu)$,
$(-1) \circ \nu(\cdot)=\sum_{j=1}^N \delta_{-\nu_j}(\cdot)$;
\begin{equation}
\E^{\nu, a} \Big[
\Theta^a(X_1^a(t)-h) \Big]
=\mathop{\Det}_{(v,v')\in C((-1) \circ \nu)^2}
\Big[ \delta(v-v')+K_{e^{h/a}}(v,v') \Big],
\label{eqn:BC0}
\end{equation}
where
\begin{equation}
K_u(v,v')=\int_{-i \infty+\delta}^{i \infty+\delta}
\frac{ds}{2 \pi i}
\Gamma(-s) \Gamma(1+s) \prod_{\ell=1}^N
\frac{\Gamma(v+\nu_{\ell})}{\Gamma(s+v+\nu_{\ell})}
\frac{u^s e^{tvs/a^2+ts^2/2a^2}}{v+s-v'}, \quad u >0.
\label{eqn:Ku1}
\end{equation}
Here the Fredholm determinant is defined by
the sum of infinite series of multiple 
contour-integrals
\begin{equation}
\mathop{\Det}_{(v,v')\in C((-1) \circ \nu)^2}
\Big[ \delta(v-v')+K_{u}(v,v') \Big]=
\sum_{L=0}^{\infty} \frac{1}{L!}
\prod_{j=1}^L \oint_{C((-1) \circ \nu)} \frac{d v_j}{2 \pi i}
\det_{1 \leq j, k \leq L} [K_{u}(v_j, v_k)],
\label{eqn:FredholmDet}
\end{equation}
where the term for $L=0$ is assumed to be 1.
Note that (\ref{eqn:Ku1}) depends on
$\nu, a$ and $t$;
$K_u(\cdot,\cdot)=K_u(\cdot, \cdot; \nu, a, t)$.

The Fredholm determinant formula (\ref{eqn:BC0}) discovered by
Borodin and Corwin \cite{BC11} is surprising, since
the O'Connell process is not determinantal as mentioned above.
We would like to understand the origin of such 
a determinantal structure surviving 
in the geometric lifting from the noncolliding BM to
the O'Connell process.

\subsection{Variations of CPR and DMR \label{sec:DMR_OP2}}

Let
\begin{equation}
\hnu(\cdot)=\sum_{j=1}^N \delta_{\hat{\nu}_j}(\cdot) \in \mM_0.
\label{eqn:hnu1}
\end{equation}
For such $\hnu$, define a function of $x \in \C$
with parameters $\u \in \C, a >0$ by
\begin{equation}
\Phi_{\hnu}^{u,a}(x)
=\Gamma(1-a(u-x)) 
\prod_{r \in \supp \hnu \cap \{u\}^{\rm c}}
\frac{\Gamma(a(r-u))}{\Gamma(a(r-x))}.
\label{eqn:Phia1}
\end{equation}
This function has simple poles at
\begin{equation}
x_n=-\frac{n}{a}+u, \quad n \in \N
\label{eqn:Phia2}
\end{equation}
and satisfies 
\begin{equation}
\Phi_{\hnu}^{\hnu_k, a}(\hnu_j)
=\delta_{jk},
\quad 1 \leq j, k \leq N.
\label{eqn:Phia3}
\end{equation}

Since 
\begin{equation}
\Gamma(a x) \sim \frac{1}{ax}, \quad x \in \C, \quad
\mbox{as $a \to 0$},
\label{eqn:a_0_b1}
\end{equation}
\begin{equation}
\lim_{a \to 0} 
\Phi_{\hnu}^{u,a}(x)
= \Phi_{\hnu}^u(x) = \prod_{r \in \supp \hnu \cap \{u\}^{\rm c}}
\frac{x-r}{u-r}, \quad x, u \in \C.
\label{eqn:a_0_b2}
\end{equation}
We say that $\Phi_{\hnu}^u(\cdot)$ is
the {\it combinatorial limit} of
$\Phi_{\hnu}^{u,a}(\cdot)$,
and that
$\Phi_{\hnu}^{u,a}(x)$ is the
{\it geometric lifting} of
$\Phi_{\hnu}^u(\cdot)$.
All poles (\ref{eqn:Phia2}) 
go to infinity in the limit $a \to 0$ and
the function becomes entire in the combinatorial limit.

Let $Z_j(t), 1 \leq j \leq N, t \geq 0$ be a set of
independent complex BM's given by (\ref{eqn:cBM1}).
In \cite{Kat12c}, we showed that the determinant
formula (\ref{eqn:BC0}) with
(\ref{eqn:FredholmDet}) of Borodin and Corwin
is rewritten as follows.
\begin{eqnarray}
&& \E^{a \circ \hnu, a}
[\Theta^a(X^a_1(t)-h)]
\nonumber\\
&& \quad
= \bE_{\widehat{\vnu}} \left[
\det_{1 \leq j, k \leq N}
\Big[\delta_{jk} 
-\Phi_{\hnu}^{\hnu_k,a}(Z_j(1/t))
\1(\Re Z_j(1/t) < h/t) \Big] \right]
\label{eqn:OP_det1}
\end{eqnarray}
for $a>0, h \in \R, t>0$.

The functions $\Phi_{\hnu}^{\hnu_k,a}(\cdot), 1 \leq k \leq N$
are not entire, but they are holomorphic 
in $\C \setminus \{-n/a+\hnu_{k} : n \in \N\}$.
Then $\Phi_{\hnu}^{\hnu_k,a}(\Z_j(1/t)), 1 \leq j, k \leq N$
are conformal martingales.
Since (\ref{eqn:Phia3}) holds, for each $t > 0$,
the RHS of (\ref{eqn:OP_det1}) is equal to
\begin{eqnarray}
&& \bE_{\widehat{\vnu}} \left[
\det_{1 \leq j, k \leq N}
\Big[\Phi_{\hnu}^{\hnu_k,a}(Z_j(1/t))
-\Phi_{\hnu}^{\hnu_k,a}(Z_j(1/t))
\1(\Re Z_j(1/t) < h/t) \Big] \right]
\nonumber\\
&& \quad 
= \bE_{\widehat{\vnu}} \left[
\det_{1 \leq j, k \leq N}
\Big[ (1- \1(\Re Z_j(1/t)< h/t))
\Phi_{\hnu}^{\hnu_k,a}(Z_j(1/t))
\Big] \right]
\nonumber\\
&& \quad 
= \bE_{\widehat{\vnu}} \left[
\prod_{j=1}^N \1(\Re Z_j(1/t) \geq h/t)
\det_{1 \leq j, k \leq N}
\Big[
\Phi_{\hnu}^{\hnu_k,a}(Z_j(1/t))
\Big] \right].
\nonumber
\end{eqnarray}
That is, we have the equality
\begin{eqnarray}
&& \E^{a \circ \hnu, a}
[\Theta^a(X^a_1(t)-h)]
\nonumber\\
&& \quad =
\bE_{\widehat{\vnu}} \left[
\prod_{j=1}^N \1(\Re Z_j(1/t) \geq h/t)
\det_{1 \leq j, k \leq N}
\Big[
\Phi_{\hnu}^{\hnu_k,a}(Z_j(1/t))
\Big] \right],
\quad a > 0, h \in \R.
\label{eqn:CPR_OP1}
\end{eqnarray}
The observable $\Theta^a(X^a_1(t)-h), h \in \R$ in the LHS
is a softening of the indicator 
$\1(X_1(t) \geq h)$.
Its expectation
for the O'Connell process started according to
the entrance law coming from ``$-\infty \vrho$"
has the representation as given by the RHS,
in which the `sharp' indicators
$\1(V_j(1/t) \geq h/t), 1 \leq j \leq N$ are
observed, but the complex weight on paths 
is `softened'(geometrically lifted).
We regard (\ref{eqn:CPR_OP1}) as
a `variation' of CPR.

\begin{prop}
\label{thm:CPR_OP1}
The O'Connell process has a variation of CPR
for $\Theta^a(X^a_1(\cdot)-h), a>0, h \in \R$
as given by (\ref{eqn:CPR_OP1}).
\end{prop}

For the noncolliding BM
$\X(t)=(X_1(t), \dots, X_N(t)), t \geq 0$
and $\vnu=(\nu_1, \dots, \nu_N) \in \W_N^{\rm A}$,
we consider the noncolliding BM with drift vector
$\vnu$, 
\begin{eqnarray}
&&\X^{\nu}(t)=(X^{\nu}_1(t), \dots, X^{\nu}_N(t)),
\nonumber\\
&& \mbox{with} \quad
X^{\nu}_j(t)=X_j(t)+\nu_j t, \quad 1 \leq j \leq N,
\quad t \geq 0.
\label{eqn:reciprocal1}
\end{eqnarray}
Put
$\nu(\cdot)=\sum_{j=1}^{N} \delta_{\nu_j}(\cdot) \in \mM_0$.
In \cite{Kat12b}, we proved the equality
\begin{equation}
\E_{\nu} \left[ \prod_{j=1}^N \1(X_j(t) \geq h) \right]
=\E_{N \delta_0} \left[ \prod_{j=1}^N
\1(X_j^{\nu}(1/t) \geq h/t) \right],
\quad t >0.
\label{eqn:reciprocal2}
\end{equation}
We call such a relation
{\it reciprocal time relation} \cite{Kat12b}.
By applying Corollary \ref{thm:CPR_noncBM},
the LHS of (\ref{eqn:reciprocal2}) is 
given by CPR and we have the equality
\begin{eqnarray}
&& \E_{N \delta_0} \left[ \prod_{j=1}^N
\1(X_j^{\nu_j}(1/t) \geq h/t) \right]
\nonumber\\
&& \quad =
\bE_{\widehat{\vnu}} \left[
\prod_{j=1}^N \1(\Re Z_j(1/t) \geq h/t)
\det_{1 \leq j, k \leq N}
\Big[
\Phi_{\nu}^{\nu_k}(Z_j(1/t))
\Big] \right],
\quad h \in \R.
\label{eqn:reciprocal3}
\end{eqnarray}
The determinant formula (\ref{eqn:CPR_OP1})
can be regarded as a geometric lifting
of the expression (\ref{eqn:reciprocal3})
for the noncolliding BM with drift.

The functions
$\Phi_{\hnu}^{\nu_k,a}(x), 1 \leq k \leq N$
are not polynomial of $x$, but here we apply
the equality (\ref{eqn:cBM0}).
By multilinearity of determinant, the RHS of (\ref{eqn:CPR_OP1})
becomes
\begin{eqnarray}
&& \rE_{\hnu} \left[
\prod_{j=1}^N \1(\Re Z_j(1/t) \geq h/t)
\det_{1 \leq j, k \leq N}
\Big[ \check{E}_0[
\Phi_{\hnu}^{\hnu_k,a}(Z_j(1/t))]
\Big] \right]
\nonumber\\
&& \quad =
\rE_{\hnu} \left[
\prod_{j=1}^N \1(\Re Z_j(1/t) \geq h/t)
\det_{1 \leq j, k \leq N}
\Big[ \sfM[
\Phi_{\hnu}^{\hnu_k,a}(iW)|(1/t,B_j(1/t))]
\Big] \right].
\nonumber
\end{eqnarray}
Then we put
\begin{equation}
\cM_{\hnu}^{u,a}(\cdot, \cdot)
=\sfM[ \Phi_{\hnu}^{u,a}(iW)|(\cdot, \cdot)],
\quad a>0, u \in \C,
\label{eqn:martingale_OP}
\end{equation}
and 
\begin{equation}
\cD_{\hnu}^a(1/t, \B(1/t))
=\det_{1 \leq j, k \leq N}
\Big[ \sfM [ \Phi_{\hnu}^{u,a}(iW)|(1/t, B_j(1/t))] \Big],
\quad a>0, \quad t > 0.
\label{eqn:DM_OP}
\end{equation}
Then (\ref{eqn:CPR_OP1}) is written as follows.
\begin{eqnarray}
&& \E^{a \circ \hnu, a}
[\Theta^a(X^a_1(t)-h)]
\nonumber\\
&& \quad =
\rE_{\hnu} \left[
\prod_{j=1}^N \1(B_j(1/t) \geq h/t)
\cD_{\hnu}^a(1/T, \B(1/T))
 \right],
 \nonumber\\
&& \qquad \qquad 
a > 0, \quad h \in \R, \quad 0 < t \leq T \leq \infty.
\label{eqn:DMR_OP1}
\end{eqnarray}

\begin{prop}
\label{thm:DMP_OP1}
The O'Connell process has a variation of DMP
for $\Theta^a(X^a_1(\cdot)-h), a>0, h \in \R$
as given by (\ref{eqn:DMR_OP1}).
Then the expectation 
of $\Theta^a(X^a_1(\cdot)-h), a>0, h \in \R$
is F-determinantal.
\end{prop}

\clearpage
\vskip 1cm
\noindent{\bf Acknowledgements} \quad
The present author would like to thank 
Hirofumi Osada for giving him
an opportunity to give lectures
at Faculty of Mathematics, Kyushu University.
He thanks
Tomoyuki Shirai, Hideki Tanemura, and Syota Esaki 
for useful discussion.



\begin{thebibliography}{99}
\bibitem{AAR99}
Andrews, G. E., Askey, R., Roy, R.:
{\it Special functions}. Cambridge: 
Cambridge University Press, 1999

\bibitem{Bai00}
Baik, J.: 
Random vicious walks and random matrices. 
Commun. Pure Appl. Math. {\bf 53}, 1385-1410 (2000)

\bibitem{BS07}
Baik, J., Suidan, T. M.:
Random matrix central limit theorems for
nonintersecting random walks.
Ann. Probab. {\bf 35}, 1807-1834 (2007)

\bibitem{BO11}
Baudoin, F., O'Connell, N.:
Exponential functionals of Brownian motion
and class-one Whittaker functions.
Ann. Inst. H. Poincar\'e, 
B {\bf 47}, 1096-1120 

\bibitem{BBO05}
Biane, P., Bougerol, P., O'Connell, N.:
Littelmann paths and Brownian paths.
Duke Math. J. {\bf 130}, 127-167 (2005)

\bibitem{BBO09}
Biane, P., Bougerol, P., O'Connell, N.:
Continuous crystal and Duistermaat-Heckman measure
for Coxeter groups.
Adv. Math. {\bf 221}, 1522-1583 (2009)

\bibitem{BPY01}
Biane, P., Pitman, J., Yor, M.:
Probability laws related to the Jacobi theta
and Riemann zeta functions,
and Brownian excursions.
Bull. Amer. Math. Soc.
{\bf 38}, 435-465 (2001)

\bibitem{BK05}
Bleher, P. M., Kuijlaars, A. B. J.: 
Integral representations for multiple Hermite
and multiple Laguerre polynomials.
Ann. Inst. Fourier.
{\bf 55}, 2001-2014 (2005)

\bibitem{BC11}
Borodin, A., Corwin, I.:
Macdonald processes.
to appear in 
Probab. Theory Relat. Fields; 
{\sf arXiv:math.PR/1111.4408}

\bibitem{BR05}
Borodin, A., Rains, E. M.:
Eynard-Mehta theorem, Schur process and their Pffaffian analogs.
J. Stat. Phys. {\bf 121}, 291-317 (2008)

\bibitem{CK03}
Cardy, J., Katori, M.:
Families of vicious walkers. 
J. Phys. A {\bf 36}, 609-629 (2003)

\bibitem{COSZ11}
Corwin, I., O'Connell, N., Sepp\"al\"ainen, T.,
Zygouras, N.:
Tropical combinatorics and Whittaker functions.
{\sf arXiv:math.PR/1110.3489}

\bibitem{Dys62}
Dyson, F. J. :
A Brownian-motion model for the eigenvalues of a random matrix.
J. Math. Phys. {\bf 3}, 1191-1198 (1962)

\bibitem{EK08}
Eichelsbacher, P., K\"onig, W.:
Ordered random walks.
Electron. J. Probab. {\bf 13}, no.46, 1307-1336 (2008)

\bibitem{EM98}
Eynard, B., Mehta, M. L. :
Matrices coupled in a chain: I.
Eigenvalue correlations.
J. Phys. A {\bf 31}, 4449-4456 (1998)

\bibitem{Fei12}
Feierl, T.:
The height of watermelons with wall. 
J. Phys. A {\bf 45}, 095003 (2012)

\bibitem{Fis84}
Fisher, M.E.:
Walks, walls, wetting, and melting. 
J. Stat. Phys. {\bf 34}, 667-729, (1984)

\bibitem{For10}
Forrester, P. J.:
{\it Log-gases and Random Matrices}.
London Mathematical Society Monographs, Princeton:
Princeton University Press, 2010

\bibitem{FNH99}
Forrester, P. J., Nagao, T., Honner, G.:
Correlations for the orthogonal-unitary and
symplectic-unitary transitions at the
hard and soft edges. 
Nucl. Phys. {\bf B553[PM]}, 601-643 (1999)

\bibitem{Fuj02}
Fujita, T. :
{\it Stochastic Calculus for Finance}
(in Japanese).
Tokyo: Kodansha, 2002

\bibitem{Fuj08b}
Fujita, T.:
{\it Random Walks and Stochastic Calculus}
(in Japanese).
Tokyo: Nihon-Hyoron-Sha, 2008

\bibitem{Fuj08}
Fujita, T., Kawanishi, Y.:
A proof of Ito's formula using a discrete Ito's formula.
Stud. Sci. Math. Hungr. {\bf 45}, 125-134 (2008)

\bibitem{Giv97}
Givental, A.:
Stationary phase integrals, quantum Toda lattices,
flag manifolds and the mirror conjecture.
In: Topics in Singular Theory,
AMS Trans. Ser. 2, vol. 180, pp.103-115,
AMS, Rhode Island (1997)

\bibitem{Gra99}
Grabiner, D. J.:
Brownian motion in a Weyl chamber,
non-colliding particles, and random matrices. 
Ann. Inst. Henri Poincar\'e,
Probab. Stat. {\bf 35}, 177-204 (1999)

\bibitem{GJ13}
Graczyk, P., Ma{\l}ecki, J.:
Multidimensional Yamada-Watanabe theorem 
and its applications to particle systems. 
J. Math. Phys. {\bf 54}, 021503/1-15 (2013)

\bibitem{Joh01}
Johansson, K.:
Universality of the local spacing distribution
in certain ensembles of Hermitian Wigner matrices.
Commun. Math. Phys. {\bf 215}, 683-705 (2001)

\bibitem{Joh02}
Johansson, K.:
Non-intersecting paths, random tilings and random matrices. 
Probab. Th. Rel. Fields {\bf 123}, 225-280 (2002)

\bibitem{Joh05}
Johansson, K.:
Non-intersecting, simple, symmetric random walks
and the extended Hahn kernel.
Ann. Inst. Fourier {\bf 55}, 2129-2145 (2005)

\bibitem{KS91}
Karatzas, I., Shreve, S. E.:
{\it Brownian Motion and Stochastic Calculus}.
2nd edition, New York:
Springer, 1991

\bibitem{Kat11}
Katori, M.:
O'Connell's process as a vicious Brownian motion.
Phys. Rev. E {\bf 84}, 061144 (2011)

\bibitem{Kat12a}
Katori, M.:
Survival probability of mutually killing Brownian motion
and the O'Connell process.
J. Stat. Phys.{\bf 147}, 206-223 (2012) 

\bibitem{Kat12b}
Katori, M.:
Reciprocal time relation of noncolliding Brownian motion
with drift.
J. Stat. Phys. {\bf 148}, 38-52 (2012).

\bibitem{Kat12c}
Katori, M.:
System of complex Brownian motions associated with
the O'Connell process.
J. Stat. Phys. {\bf 149}, 411-431 (2012)

\bibitem{Kat13}
Katori, M.:
Determinantal martingales and noncolliding diffusion processes.\\
{\sf arXiv:math.PR/1305.4412}

\bibitem{Kat13b}
Katori, M. Determinantal martingales and correlations
of noncolliding random walks.
{\sf arXiv:math.PR/1307.1856}

\bibitem{KT03a}
Katori, M., Tanemura, H.:
Functional central limit theorems for vicious walkers.
Stoch. Stoch. Rep. {\bf 75}, 369-390 (2003)

\bibitem{KT04}
Katori, M., Tanemura, H.:
Symmetry of matrix-valued stochastic processes and
noncolliding diffusion particle systems. 
J. Math. Phys. {\bf 45}, 3058-3085 (2004)

\bibitem{KT07a}
Katori, M., Tanemura, H.:
Infinite systems of noncolliding generalized meanders
and Riemann-Liouville differintegrals.
Probab. Theory Relat. Fields {\bf 138}, 113-156 (2007)

\bibitem{KT09}
Katori, M., Tanemura, H.:
Zeros of Airy function and relaxation process.
J. Stat. Phys. {\bf 136}, 1177-1204 (2009)

\bibitem{KT10}
Katori, M., Tanemura, H.:
Non-equilibrium dynamics of Dyson's model
with an infinite number of particles.
Commun. Math. Phys. {\bf 293}, 469-497 (2010)

\bibitem{KT11}
Katori, M., Tanemura, H.:
Noncolliding squared Bessel processes.
J. Stat. Phys. {\bf 142}, 592-615 (2011)

\bibitem{KT11b}
Katori, M. and Tanemura, H.:
Noncolliding processes, matrix-valued processes
and determinantal processes. 
Sugaku Expositions (AMS) 
{\bf 24}, 263-289 (2011)

\bibitem{KT13}
Katori, M., Tanemura, H.:
Complex Brownian motion representation
of the Dyson model.
Electron. Commun. Probab. {\bf 18}, no.4, 1-16 (2013)

\bibitem{KL01}
Kharchev, S., Lebedev, D.:
Integral representations for the eigenfunctions
of quantum open and periodic Toda chains
from the QISM formalism.
J. Phys. A: Math. Gen.{\bf 34}, 2247-2258 (2001)

\bibitem{Koe05}
K\"onig, W.:
Orthogonal polynomial ensembles in probability theory.
Probab. Surveys {\bf 2}, 385-447 (2005)

\bibitem{KO01}
K\"onig, W., O'Connell, N.:
Eigenvalues of the Laguerre process as non-colliding 
squared Bessel process.
Elec. Commun. Probab. {\bf 6}, 107-114 (2001)

\bibitem{Kra90}
Krattenthaler, C.:
Generating functions for plane partitions
of a given shape.
Manuscripta Math. {\bf 69}, 173-202 (1990)

\bibitem{Kra99}
Krattenthaler, C.:
Advanced determinant calculus.
S\'eminaire Lotharingien Combin. 42 (1999) (The Andrews Festschrift), 
paper B42q, 67 pp

\bibitem{Lev96}
Levin, B. Ya.:
{\it Lectures on Entire Functions}. 
Translations of Mathematical Monographs, {\bf 150},
Province R. I. :
Amer. Math. Soc.,1996

\bibitem{Mac99}
Macdonald, I. G.:
{\it Symmetric Functions and Hall Polynomials}. 
2nd edition, New Yort: Oxford University Press, 1999

\bibitem{MY00}
Matsumoto, H., Yor, M.:
An analogue of Pitman's $2M-X$ theorem
for exponential Wiener functionals,
Part I: A time-inversion approach.
Nagoya Math. J. {\bf 159}, 125-166 (2000)

\bibitem{MY05}
Matsumoto, H., Yor, M.:
Exponential functionals of Brownian motion
I: Probability laws at fixed time.
Probab. Surveys {\bf 2}, 312-347 (2005)

\bibitem{Meh04}
Mehta, M. L.:
{\it Random Matrices}. 3rd edition, 
Amsterdam: Elsevier, 2004

\bibitem{NF98}
Nagao, T. and Forrester, P.:
Multilevel dynamical correlation functions for Dyson's
Brownian motion model of random matrices.
Phys. Lett. {\bf A247}, 42-46 (1998)

\bibitem{NF02}
Nagao, T., Forrester, P.J.:
Vicious random walkers and a discretization of
Gaussian random matrix ensembles.
Nucl. Phys. B {\bf 620} [FS], 551-565 (2002)

\bibitem{OCo12a}
O'Connell, N.:
Directed polymers and the quantum Toda lattice.
Ann. Probab. {\bf 40}, 437-458 (2012)

\bibitem{OCo12b}
O'Connell, N.:
Whittaker functions and related stochastic processes.
to appear in proceedings of Fall 2010 MSRI semester 
`Random matrices, interacting particle systems and integrable systems',
{\sf arXiv:math.PR/1201.4849}

\bibitem{Osa12}
Osada, H.:
Infinite-dimensional stochastic differential equations
related to random matrices.
Probab. Theory Relat. Fields {\bf 153}, 471-509 (2012)

\bibitem{Osa13}
Osada, H.:
Interacting Brownian motions in infinite dimensions
with logarithmic potentials.
Ann. Probab. {\bf 41}, 1-49 (2013)

\bibitem{Osa13b}
Osada, H.:
Interacting Brownian motions in infinite dimensions
with logarithmic potentials II:
Airy random point field.
Stoch. Proc. Appl. {\bf 123}, 813-838 (2013)

\bibitem{RY05}
Revuz, D., Yor, M.:
{\it Continuous Martingales and Brownian Motion.} 3rd edition,
New York: Springer, 2005

\bibitem{Sch00}
Schoutens, W.:
{\it Stochastic Processes and Orthogonal Polynomials.}
Lecture Notes in Statistics {\bf 146}, 
New York: Springer, 2000

\bibitem{STS94}
Semenov-Tian-Shansky, M. A.:
Quantization of open Toda lattices.
In: Dynamical Systems VII:
Integrable Systems, Nonholonomic Dynamical Systems. 
Edited by V. I. Arnol'd and S. P. Novikov.
Encyclopaedia of Mathematical Sciences,
vol.16. pp.226-259, Springer, Berlin (1994)

\bibitem{ST03}
Shirai, T., Takahashi, Y.: 
Random point fields associated with certain
Fredholm determinants I:
fermion, Poisson and boson point process.
J. Funct. Anal.
{\bf 205}, 414-463 (2003)

\bibitem{Skl85}
Sklyanin, E. K.:
The quantum Toda chain.
In: {\it Non-linear Equations in Classical
and Quantum Field Theory}, 
Lect. Notes in Physics,
{\bf 226}, pp. 195-233, Berlin: Springer, 1985

\bibitem{Sos00}
Soshnikov, A. : 
Determinantal random point fields.
Russian Math. Surveys {\bf 55}, 923-975 (2000)

\bibitem{Spo87}
Spohn, H.:
Interacting Brownian particles:
a study of Dyson's model.
In:
{\it Hydrodynamic Behavior and Interacting Particle Systems}. 
G. Papanicolaou (ed),  
IMA Volumes in Mathematics and its Applications, {\bf 9}, Berlin: 
Springer-Verlag, 1987, pp.151-179

\bibitem{Sze91}
Szeg\"o, G.:
{\it Orthogonal Polynomials}.
Providence:
Amer. Math. Soc., 1991

\bibitem{Tao12}
Tao, T.:
{\it Topics in Random Matrix Theory}.
Providence:
Amer. Math. Soc., 2012

\bibitem{TW04}
Tracy, C. A., Widom, H.:
Differential equations for Dyson processes.
Commun. Math. Phys. {\bf 252}, 7-41 (2004)

\bibitem{TW07}
Tracy, C. A., Widom, H.:
Nonintersecting Brownian excursions,
Ann. Appl. Probab. {\bf 17}, 953-979 (2007)

\bibitem{Wal92}
Wallach, N. R.:
{\it Real Reductive Groups II}.
San Diego CA: Academic Press,  1992

\bibitem{Wat44}
Watson, G. N.:
{\it A Treatise on the Theory of Bessel Functions.}
2nd edition, 
Cambridge: Cambridge University Press, 1944

\end{thebibliography}
\end{document}